\documentclass{article}

\usepackage[french,british]{babel}

\usepackage[T1]{fontenc}

\usepackage{amsmath}

\usepackage{amsfonts}
\usepackage{amsbsy}
\usepackage{amssymb}

\usepackage{amsthm}

\theoremstyle{plain}

\usepackage{latexsym}
\usepackage{amsmath}
\usepackage{amsfonts}
\usepackage{amssymb}
\usepackage{psfrag}
\usepackage{graphicx}

\usepackage{amsmath,amssymb}
\usepackage{graphicx}

\usepackage{amsmath,amssymb}
\usepackage{graphicx}

\newtheorem{ozn}{Definition}[section]
\newtheorem{thm}{Theorem}[section]

\newtheorem{nas}{Corollary}[section]
\newtheorem{zau}{Remark}[section]
\newtheorem{lema}{Lemma}[section]
\newtheorem{pry}{Example}[section]

\newcommand{\eps}{\varepsilon}
\newcommand{\me}{\mathbf}
\newcommand{\mr}{\mathbb}
\newcommand{\mt}{\mathsf}
\newcommand{\md}{\mathcal}
\newcommand{\ld}{\left}
\newcommand{\rd}{\right}

\newcommand{\ip}{\int_{-\pi}^{\pi}}

\newcommand{\be}{\begin{equation}}
\newcommand{\ee}{\end{equation}}

\newcommand{\bem}{\begin{multline}}
\newcommand{\eem}{\end{multline}}

\newcommand{\bml}{\begin{multline*}}
\newcommand{\eml}{\end{multline*}}

\newcommand{\beg}{\begin{gather}}
\newcommand{\eeg}{\end{gather}}

\begin{document}

\title{Minimax-robust forecasting of sequences with periodically stationary long memory multiple seasonal increments}

\author{
Maksym Luz\thanks {maksim-luz@ukr.net},
Mikhail Moklyachuk\thanks
{Department of Probability Theory, Statistics and Actuarial
Mathematics, Taras Shevchenko National University of Kyiv, Kyiv 01601, Ukraine, moklyachuk@gmail.com}
}

\date{\today}

\maketitle

\renewcommand{\abstractname}{Abstract}
 \begin{abstract}
We introduce stochastic sequences $\zeta(k)$ with periodically stationary generalized multiple increments of fractional order which combines cyclostationary, multi-seasonal, integrated and fractionally integrated patterns. We solve the problem of optimal estimation of linear functionals
constructed from unobserved values of stochastic sequences $\zeta(k)$  based on their observations at points $ k<0$.
For sequences with known spectral densities, we obtain formulas for calculating values of the mean square errors and the spectral characteristics of the optimal estimates of functionals.
Formulas that determine the least favorable spectral densities and minimax (robust) spectral
characteristics of the optimal linear estimates of functionals
are proposed in the case where spectral densities of sequences
are not exactly known while some sets of admissible spectral densities are given.
\end{abstract}

\vspace{2ex}

\maketitle

\vspace{2ex}

\section{Introduction}

A variety of non-stationary and long memory time series models are introduced and investigated by researchers in the past decade (see, for example, papers by Dudek and Hurd \cite{Dudek-Hurd}, Johansen and Nielsen  \cite{Johansen}, Reisen et al.\cite{Reisen2014}).
Such models are used when analyzing data which arise in different field of economics, finance, climatology, air pollution, signal processing.

Since the book by Box and Jenkins (1970), autoregressive moving average (ARMA) models integrated of order $d$ are standard models used for time series
analysis. These models are described by the equation
\be
 \psi (B) (1-B) ^ dx_t = \theta (B) \eps_t,
 \label{ARIMA_model} \ee
 where $ \eps_t $, $ t \in \mr Z $, is a sequence of zero mean i.i.d. random variables,  $ \psi (z) $, $ \theta ( z) $ are polynomials of $ p $ and $ q $ degrees
 respectively with roots outside the unit circle.
 This integrated ARIMA model is generalized by adding a seasonal component. A new model is described by the equation (see new edition of the book by Box and Jenkins \cite{Box_Jenkins} for detailes)
 \be
 \Psi (B ^ s) (1-B ^ s) ^ Dx_t = \Theta (B ^ s) \eps_t, \label{seasonal_2_model} \ee
where $ \Psi (z) $ and $ \Theta (z) $ are polynomials of degrees of $ P $ and $ Q $ respectively which have roots outside the unit circle.

When an ARIMA  sequence determined by equation (\ref{ARIMA_model}) is inserted in relation (\ref{seasonal_2_model}) instead of $ \eps_t$ we have
general multiplicative model
\be
 \Psi (B ^ s) \psi (B) (1-B) ^ d (1-B ^ s) ^ Dx_t = \Theta (B ^ s) \theta (B) \eps_t \label {seasonal_3_model} \ee
with parameters $ (p, d, q) \times (P, D, Q)_s $, $d,D\in\mr N^*$, called SARIMA $(p, d, q)\times(P, D, Q)_s$ model.

A good performance is shown by models which include fractional integration, that is when parameters $d$ and $D$ are fractional.
When $|d+D|<1/2$ and $|D|<1/2$, a process described by equation (\ref{seasonal_3_model}) is stationary and invertible.
We refer to the paper by Porter-Hudak \cite{Porter-Hudak} who studied  seasonal ARFIMA models and applied them to the monetary aggregates used by U.S. Federal Reserve.

 Closely related to fractionally integrated ARMA and GARMA processes described by equation
 \be
 (1-2uB+B^2) ^ dx_t = \eps_t,\quad |u|\leq1. \label{GARMA_model} \ee
  is SARFIMA process.
These processes were introduced and studied by  Granger and Joyeux \cite{Granger_Joyeux}, Hosking \cite{Hosking}, Andel \cite{Andel}, Gray et al. \cite{Gray} in order to model long-memory stationary time series.
Fractionally integrated models are a powerful tool for studying a variety of real world processes.
For the resent works dedicated to statistical inference for  seasonal long-memory sequences,
we refer to Arteche and   Robinson \cite{Arteche}, who applied the log-periodogram and Gaussian or Whittle methods of memory parameters estimation for seasonal/cyclical asymmetric long memory processes with application to UK inflation data,  and also Tsai,  Rachinger and  Lin \cite{Tsai}, who developed methods of estimation of parameters in case of measurement errors.
Baillie,   Kongcharoen and  Kapetanios \cite{Baillie} compared MLE and semiparametric estimation procedures for prediction problems based on ARFIMA models.
 Based on simulation study, they indicate better performance of  MLE predictor than the one  based on two-step local Whittle estimation.
Hassler and  Pohle \cite{Hassler} (see also Hassler \cite{Hassler_book}) assess a predictive  performance of various methods of  forecasting  of inflation and return volatility time series and show strong evidences for models with a fractional integration component.

Another type of non-stationarity is described by stochastic processes with time-dependent spectrum.
A wide class of processes with time-dependent spectrum is formed by  periodically correlated, or cyclostationary, processes introduced by Gladyshev \cite{Gladyshev}.
These processes are widely used in signal processing and communications (see Napolitano \cite{Napolitano} for a review of recent works on cyclostationarity and its applications).
Periodic time series may be considered as an extension of a SARIMA model (see Lund \cite{Lund} for a test assessing if a PARMA model is preferable to a SARMA one) and are suitable for forecasting stream flows with quarterly, monthly or weekly cycles (see Osborn \cite{Osborn}).

Baek, Davis and Pipiras \cite{Baek} introduced a periodic dynamic factor model (PDFM) with periodic vector autoregressive (PVAR) factors, in contrast to seasonal VARIMA factors.

The models mentioned above are used in estimation of model's parameters and forecast issues.
Meanwhile a direct application of the developed results to real data may lead to significant increasing of errors of estimates due to presence of outliers, measurement errors, incomplete information about the spectral, or model, structure etc. From this point of view, we see an increasing interest to robust methods of estimation that are reasonable in such cases.
For example, Reisen,  et al. \cite{Reisen2018} proposed a semiparametric robust estimator for fractional parameters in the SARFIMA model and illustrated its application
to forecast of sulfur dioxide $SO_2$ pollutant concentrations. Solci at al. \cite{Solci} proposed robust estimates of periodic  autoregressive (PAR) model.

Robust approaches  are successfully  applied  to the problem of estimation of linear functionals from unobserved values of stochastic processes.
The paper by Grenander \cite{Grenander} should be marked as the first one where the minimax extrapolation problem for stationary processes was
formulated as a game of two players and solved.
 Hosoya \cite{Hosoya}, Kassam  \cite{Kas1982}, Franke \cite{Franke1985}, Vastola and  Poor \cite{VastPoor1984}, Moklyachuk \cite{Moklyachuk,Moklyachuk2015} studied minimax  extrapolation (forecasting), interpolation (missing values estimation) and filtering (smoothing) problems for stationary sequences and processes. Recent results of minimax extrapolation problems for stationary vector-valued processes and periodically correlated processes belong to  Moklyachuk and Masyutka \cite{Mokl_Mas_extr,Mokl_Mas_pred} and  Moklyachuk and  Golichenko (Dubovetska) \cite{Dubovetska6} respectively. Processes with stationary increments are investigated by Moklyachuk and Luz  \cite{Luz_Mokl_extra,Luz_Mokl_book}.
 We also mention works  by Moklyachuk and Sidei~\cite{Sidei_extr,Sidei_book}, who derive minimax estimates of stationary processes from observations with missed values.
 Moklyachuk and  Kozak \cite{Kozak_Mokl} studied interpolation problem for stochastic sequences with periodically stationary increments.

In this article, we present results of investigation of stochastic sequences with periodically stationary long memory multiple seasonal increments
motivated by articles by Dudek \cite{Dudek}, Gould et al. \cite{Gould} and Reisen et al. \cite{Reisen2014}, who considered models with multiple seasonal patterns for inference and forecasting, and Hurd and Piparas \cite{Hurd}, who introduced two models of periodic autoregressive time series with multiple periodic coefficients.

In Section \ref{spectral_ theory}, we give definition of
generalized multiple (GM)  increment sequence $\chi_{\overline{\mu},\overline{s}}^{(d)}(\vec{\xi}(m))$ and introduce stochastic
sequences $\zeta(m)$ with periodically stationary (periodically correlated, cyclostationary) GM increments.
Such kind of non-stationary stochastic sequence combines  periodic structure of the covariation function of the sequences as well as multiple seasonal factors,
including the integrating one. The section also contains a short review of the spectral theory of vector-valued GM increment sequences.
Section \ref{classical_extrapolation} deals with the classical estimation problem for linear functionals $A\zeta$ and $A_N\zeta $ which are constructed from unobserved values of the sequence $\zeta(m)$ when the spectral structure of the sequence $\zeta(m)$ is known.
Estimates are obtained by representing the sequence $\zeta(m)$ as a vector sequence $\vec \xi(m)$ with stationary GM  increments and applying the Hilbert space projection technique.
An approach to forecasting in the presence of non-stationary  fractional integration is discussed in Section~\ref{fractional_extrapolation}.
Section \ref{examples} contains examples of forecasting of particular models of time series.
In Section \ref{minimax_extrapolation}, we derive the minimax (robust) estimates in cases, where spectral densities of sequences are not exactly known
while some sets of admissible spectral densities are specified which are generalizations of the corresponding sets of admissible spectral densities described in a survey article by
Kassam and Poor \cite{Kassam} for stationary stochastic processes.

\section{Stochastic sequences with periodically stationary generalized multiple increments}\label{spectral_ theory}

\subsection{Definition and spectral representation of a periodically stationary GM increment}

In this section, we present definition, justification and a brief review of the spectral theory of stochastic sequences with periodically stationary multiple seasonal increments.
This type of stochastic sequences will allow us to deal with a wide range of non-stationarity in time series analysis.

 Consider a stochastic sequence $\{\eta(m),m\in\mathbb Z\}$. By $B_{\mu}$ denote a backward shift operator with the step $\mu\in
\mathbb Z$, such that $B_{\mu}\eta(m)=\eta(m-\mu)$; $B:=B_1$. Recall the following definition \cite{Luz_Mokl_book,Pinsker,Yaglom}.

\begin{ozn}\label{def_Pryrist}
For a given stochastic sequence $\{\eta(m),m\in\mathbb Z\}$, the
sequence
\begin{equation}
\label{Pryrist}
\eta^{(n)}(m,\mu)=(1-B_{\mu})^n\eta(m)=\sum_{l=0}^n(-1)^l{n \choose l}\eta(m-l\mu),
\end{equation}
where ${n \choose l}=\frac{n!}{l!(n-l)!}$,  is called
stochastic $n$th increment sequence with a step $\mu\in\mathbb Z$.
\end{ozn}

Note that in Definition \ref{def_Pryrist},  the step parameter $\mu$ is not fixed and   varies over the set $\mr Z$.
The introduced increment (\ref{Pryrist}) is applicable for describing the integrated stochastic sequence (\ref{ARIMA_model}).
The varying step $\mu$ provides a flexibility of the integrated processes.
For instance, let a sequence $x_m$ satisfy the equation  $x_m=x_{m-1}+\eps_m+a\eps_{m-1}$.
Then the $\mu$-step increment $x_{m}-x_{m-\mu}=\sum_{k=0}^{\mu-1}(x_{m-k}-x_{m-k-1})$ is stationary as a sum of stationary $1$-step increments.
To deal with seasonal time series (\ref{seasonal_2_model}) we need to extend definition of stochastic increment sequence as follows.

\begin{ozn}\label{def_seasonal_Pryrist}
For a given stochastic sequence $\{\eta(m),m\in\mathbb Z\}$, the
sequence
\begin{equation}
\label{seasonal_Pryrist}
\eta^{(n)}_s(m,\mu)=(1-B_{\mu}^s)^n\eta(m)=\sum_{l=0}^n(-1)^l{n \choose l}\eta(m-l\mu s)
\end{equation}
is called stochastic seasonal increment sequence with a fixed seasonal parameter $s\in \mr N^*=\mr N\setminus\{0\}$ and a varying step $\mu\in\mathbb Z$.
\end{ozn}

\begin{zau}
For $s=1$, under the seasonal increment $\eta^{(n)}_1(m,\mu)$ we  understand the increment $\eta^{(n)}(m,\mu)$ from Definition \ref{def_Pryrist}.
\end{zau}

We mention the following properties of the  seasonal increment sequence $\eta^{(n)}_s(m,\mu)$, which will be used for proving Theorem \ref{thm1}:
\begin{eqnarray}\label{tot1}
    \eta^{(n)}_s(m,-\mu)& = &(-1)^{n}\eta^{(n)}_s(m+n\mu s,\mu),
\\
\label{tot2}
    \eta^{(n)}_s(m,\mu)& = &\sum_{l=0}^{(\mu-1)n}A_{l}\eta^{(n)}(m-ls,1), \quad \mu>0,
\end{eqnarray}
  where  $\{A_l,l=0,1,2,\ldots,(\mu-1)n\}$ are coefficients
from the representation
$$(1+x+\ldots+x^{\mu-1})^n=\sum_{l=0}^{(\mu-1)n}A_lx^l.$$

General multiplicative model (\ref{seasonal_3_model}) \cite{Box_Jenkins} indicates the necessity of dealing with increments of different seasonal parameters.
Moreover, for each season factor at each differencing order it is possible to make different steps by applying the operator
$(1-B_{\mu_1}^s)\cdot\ldots\cdot(1-B_{\mu_n}^s)$ instead of $(1-B_{\mu}^s)^n$. Thus, the following generalization is reasonable.

\begin{ozn}\label{def_multiplicative_Pryrist}
For a given stochastic sequence $\{\eta(m),m\in\mathbb Z\}$, the
sequence
\begin{gather}
\nonumber
\chi_{\overline{\mu},\overline{s}}^{(d)}(\eta(m)):=\chi_{\overline{\mu},\overline{s}}^{(d)}(B)\eta(m)
=(1-B_{\mu_1}^{s_1})^{d_1}(1-B_{\mu_2}^{s_2})^{d_2}\cdot\ldots\cdot(1-B_{\mu_r}^{s_r})^{d_r}\eta(m)
\\=\sum_{l_1=0}^{d_1}\ldots \sum_{l_r=0}^{d_r}(-1)^{l_1+\ldots+ l_r}{d_1 \choose l_1}\cdot\ldots\cdot{d_r \choose l_r}\eta(m-\mu_1s_1l_1-\cdots-\mu_rs_rl_r)
\label{GM_Pryrist}
\end{gather}
is called stochastic  generalized multiple (GM)  increment sequence of differentiation   order
$d:=d_1+d_2+\ldots+d_r$, $\overline{d}=(d_1,d_2,\ldots,d_r)\in (\mr N^*)^r$,
with a fixed seasonal  vector $\overline{s}=(s_1,s_2,\ldots,s_r)\in (\mr N^*)^r$
and a varying step $\overline{\mu}=(\mu_1,\mu_2,\ldots,\mu_r)\in (\mr N^*)^r$ or $\in (\mr Z\setminus\mr N)^r$.
\end{ozn}

\begin{pry}
Seasonal autoregressive integrated moving average  (SARIMA) model $\{x_m, m\in \mr Z\}$ with multiple period is defined by the difference equation
\[
\phi(B)(1-B)^d\prod_{i=1}^r\Phi_i(B^{s_i})(1-B^{s_i})^{d_i}x_m
=\theta(B) \prod_{i=1}^r\Theta_i(B^{s_i}) \eps_m,
\]
where all  roots of polynomials $\phi(z)$, $\theta(z)$, $\Phi_i(z)$, $\Theta_i(z)$ lie outside the unit circle, $1<s_1<\ldots<s_r$. The sequence $y_m=(1-B)^d\prod_{i=1}^r(1-B^{s_i})^{d_i}x_m$ is stationary in this case and we can define a GM increment sequence $\chi_{\overline{\mu},\overline{s}}^{(d)}(x_m)$ such that $\chi_{\overline{1},\overline{s}}^{(d)}(x_m)=y_m$, $m\in \mr Z$.
\end{pry}

Let $\gamma$ denotes a triple $(\overline{\mu},\overline{s},\overline{d})$.
For $i=1,2,\ldots, r$, $j\in \mr Z$ define coefficients $M_i^j:=\ld[\frac{j}{\mu_is_i}\rd]$ and  $I_i^j:=\mr I\{j\mod \mu_is_i=0\}$, where $\mr I\{\cdot\}$ is the indicator function,
and notations $n_i:=\mu_is_id_i$, $\langle\overline{s},\overline{\mu},\overline{d}\rangle_k:=\sum_{i=1}^k\mu_is_id_i=\sum_{i=1}^kn_i$, $n(\gamma):=\langle\overline{s},\overline{\mu},\overline{d}\rangle_r$. Denote the maximun of two numbers as $x \vee y$ and the minimum as $x \wedge y$.

\begin{lema} \label{lem_multiplicative_Pryrist}
 The  multiplicative increment operator $\chi_{\overline{\mu},\overline{s}}^{(d)}(B)$ admits a representation
\[
\chi_{\overline{\mu},\overline{s}}^{(d)}(B)
=\prod_{i=1}^r(1-B_{\mu_i}^{s_i})^{d_i}
=\sum_{k=0}^{n(\gamma)}e_{\gamma}(k)B^k,
\]
\begin{gather*}
e_{\gamma}(k)=\sum_{k_{r-1}=0\vee k_r-n_r}^{\langle\overline{s},\overline{\mu},\overline{d}\rangle_{r-1}\wedge k_r}\
\sum_{k_{r-2}=0\vee k_{r-1}-n_{r-1}}^{\langle\overline{s},\overline{\mu},\overline{d}\rangle_{r-2}\wedge k_{r-1}}\ldots
\sum_{k_1=0\vee k_2-n_2}^{\langle\overline{s},\overline{\mu},\overline{d}\rangle_{1}\wedge k_2}\ld((-1)^{\sum_{i=1}^rM_i^{k_i-k_{i-1}}}\rd.
\\
\ld.\times \prod_{i=1}^rI_i^{k_i-k_{i-1}}
\prod_{i=1}^r{d_i \choose M_i^{k_i-k_{i-1}}}\rd),
\end{gather*}
where $k_0:=0$, $k_r:=k$.
\end{lema}
\begin{proof}
See Appendix.
\end{proof}

\begin{ozn}
\label{oznStPryrostu}
A stochastic GM increment sequence $\chi_{\overline{\mu},\overline{s}}^{(d)}(\eta(m))$ generated
by a stochastic sequence $\{\eta(m),m\in\mathbb Z\}$ is wide sense
stationary if the mathematical expectations
\begin{eqnarray*}
\mt E\chi_{\overline{\mu},\overline{s}}^{(d)}(\eta(m_0))& = &c^{(d)}_{\overline{s}}(\overline{\mu}),
\\
\mt E\chi_{\overline{\mu}_1,\overline{s}}^{(d)}(\eta(m_0+m))\chi_{\overline{\mu}_2,\overline{s}}^{(d)}(\eta(m_0))
& = & D^{(d)}_{\overline{s}}(m;\overline{\mu}_1,\overline{\mu}_2)
\end{eqnarray*}
exist for all $m_0,m,\overline{\mu},\overline{\mu}_1,\overline{\mu}_2$ and do not depend on $m_0$.
The function $c^{(d)}_{\overline{s}}(\overline{\mu})$ is called mean value  and the function $D^{(d)}_{\overline{s}}(m;\overline{\mu}_1,\overline{\mu}_2)$ is
called structural function of the stationary GM increment sequence.
\\
The stochastic sequence $\{\eta(m),m\in\mathbb   Z\}$
determining the stationary GM increment sequence
$\chi_{\overline{\mu},\overline{s}}^{(d)}(\eta(m))$ by   \eqref{GM_Pryrist} is called stochastic
sequence with stationary GM increments (or GM increment sequence of order $d$).
\end{ozn}

The following theorem is a generalization of the corresponding theorem for stochastic increment sequence $\eta^{(d)}(m,\mu)$ \cite{Luz_Mokl_book,Yaglom:1955}.

\begin{thm}\label{thm1}
The mean value  and the structural function
  of the stochastic stationary GM sequence $\chi_{\overline{\mu},\overline{s}}^{(d)}(\eta(m))$ can be represented in the
 forms
\begin{eqnarray}
\label{serFnaR}
c^{(d)}_{ \overline{s}}(\overline{\mu})& = &c\prod_{i=1}^r\mu_i^{d_i},
\\
\label{strFnaR}
 D^{(d)}_{\overline{s}}(m;\overline{\mu}_1,\overline{\mu}_2)& = &\int_{-\pi}^{\pi}e^{i\lambda
m} \chi_{\overline{\mu}_1}^{(d)}(e^{-i\lambda})\chi_{\overline{\mu}_2}^{(d)}(e^{i\lambda})\frac{1}
{|\beta^{(d)}(i\lambda)|^2}dF(\lambda),
\end{eqnarray}
where
\[\chi_{\overline{\mu}}^{(d)}(e^{-i\lambda})=\prod_{j=1}^r(1-e^{-i\lambda\mu_js_j})^{d_j}, \quad \beta^{(d)}(i\lambda)= \prod_{j=1}^r\prod_{k_j=-[s_j/2]}^{[s_j/2]}(i\lambda-2\pi i k_j/s_j)^{d_j},
\]
 $c$ is a constant, $F(\lambda)$ is a left-continuous
nondecreasing bounded function. The constant $c$
and the function $F(\lambda)$ are determined uniquely by the GM
increment sequence $ \chi_{\overline{\mu},\overline{s}}^{(d)}(\eta(m))$.
\\
On the other hand, a function $c^{(d)}_{ \overline{s}}(\mu)$ which has form
$(\ref{serFnaR})$ with a constant $c$ and a function
$D^{(d)}_{\overline{s}}(m;\overline{\mu}_1,\overline{\mu}_2)$ which has form $(\ref{strFnaR})$ with a
function $F(\lambda)$ satisfying the indicated conditions are
the mean value and the structural function of a stationary GM
increment sequence $\chi_{\overline{\mu},\overline{s}}^{(d)}(\eta(m))$.
\end{thm}
\begin{proof}
See Appendix.
\end{proof}

Note that  by the spectral function and the spectral density of a stochastic sequence with stationary GM increments, we will call the spectral function and the spectral density of the corresponding stationary GM increment sequence. Representation $(\ref{strFnaR})$ and the Karhunen theorem
 \cite{Karhunen,Gihman_Skorohod} imply the spectral representation of
the stationary GM increment sequence $\chi_{\overline{\mu},\overline{s}}^{(d)}(\eta(m))$:
\begin{equation}
\label{predZnaR}
\chi_{\overline{\mu},\overline{s}}^{(d)}(\eta(m))=\int_{-\pi}^{\pi}
e^{im\lambda}\chi_{\overline{\mu}}^{(d)}(e^{-i\lambda})\frac{1}{\beta^{(d)}(i\lambda)}dZ_{\eta^{(d)}}(\lambda),
\end{equation}
where $Z_{\eta^{(d)}}(\lambda)$ is a stochastic process with uncorrelated increments on $[-\pi,\pi)$ connected with the spectral function $F(\lambda)$ by
the relation
\[
 \mt E|Z_{\eta^{(d)}}(\lambda_2)-Z_{\eta^{(d)}}(\lambda_1)|^2=F(\lambda_2)-F(\lambda_1)<\infty, \quad -\pi\leq \lambda_1<\lambda_2<\pi.
\]

Finally, we are ready to give a definition of  periodically stationary  GM increment sequence.

\begin{ozn}
\label{OznPeriodProc}
A stochastic sequence $\{\zeta(m),m\in\mathbb Z\}$ is called  stochastic
sequence with periodically stationary (periodically correlated) GM increments with period $T$ if the mathematical expectations
\begin{eqnarray*}
\mt E\chi_{\overline{\mu},T\overline{s}}^{(d)}(\zeta(m+T)) & = & \mt E\chi_{\overline{\mu},T\overline{s}}^{(d)}(\zeta(m))=c^{(d)}_{T\overline{s}}(m,\overline{\mu}),
\\
\mt E\chi_{\overline{\mu}_1,T\overline{s}}^{(d)}(\zeta(m+T))\chi_{\overline{\mu}_2,T\overline{s}}^{(d)}(\zeta(k+T))
& = & D^{(d)}_{T\overline{s}}(m+T,k+T;\overline{\mu}_1,\overline{\mu}_2)
\\
& = & D^{(d)}_{T\overline{s}}(m,k;\overline{\mu}_1,\overline{\mu}_2)
\end{eqnarray*}
exist for every  $m,k,\overline{\mu}_1,\overline{\mu}_2$ and  $T>0$ is the least integer for which these equalities hold.
\end{ozn}

It follows from  Definition \ref{OznPeriodProc} that the sequence
\begin{equation}
\label{PerehidXi}
\xi_{p}(m)=\zeta(mT+p-1), \quad p=1,2,\dots,T; \quad m\in\mathbb Z
\end{equation}
forms a vector-valued sequence
$\vec{\xi}(m)=\left\{\xi_{p}(m)\right\}_{p=1,2,\dots,T}, m\in\mathbb Z$
with stationary GM increments as follows:
\[
\chi_{\overline{\mu},\overline{s}}^{(d)}(\xi_p(m))=\sum_{l_1=0}^{d_1}\ldots \sum_{l_r=0}^{d_r}(-1)^{l_1+\ldots+ l_r}{d_1 \choose l_1}\cdot\ldots\cdot{d_r \choose l_r}\xi_p(m-\mu_1s_1l_1-\ldots-\mu_rs_rl_r)\]
\[
=\sum_{l_1=0}^{d_1}\ldots \sum_{l_r=0}^{d_r}(-1)^{l_1+\ldots+ l_r}{d_1 \choose l_1}\cdot\ldots\cdot{d_r \choose l_r}\zeta((m-\mu_1s_1l_1-\ldots-\mu_rs_rl_r)T+p-1)\]
\[=\chi_{\overline{\mu},T\overline{s}}^{(d)}(\zeta(mT+p-1)),\quad p=1,2,\dots,T,\]
where $\chi_{\overline{\mu},\overline{s}}^{(d)}(\xi_p(m))$ is the GM increment of the $p$-th component of the vector-valued sequence $\vec{\xi}(m)$.
The following theorem describes a spectral representation of the sequence $\vec{\xi}(m)$.

\begin{pry}
Define a periodic seasonal autoregressive integrated moving average model (PSARIMA) $\{X_m,m\in \mr Z\}$,  with multiple seasonal patterns by relation
\[
\phi_m(B)(1-B^T)^d\prod_{i=1}^r\Phi_{i,m}(B)(1-B^{Ts_i})^{d_i}X_m
=\theta_m(B) \prod_{i=1}^r\Theta_{i,m}(B) \eps_m,
\]
where all  polynomials $\phi_m(z)$, $\theta_m(z)$, $\Phi_{i,m}(z)$, $\Theta_{i,m}(z)$ are $T$-periodic by parameter  $m$ functions, $1<s_1<\ldots<s_r$.
Define
\begin{eqnarray*}
\Phi_m(z)& := &\phi_m(z)\prod_{i=1}^r\Phi_{i,m}(z)=\sum_{k=0}^{q_1}\Phi_m(k)z^k,
\\
\Theta_m(z)& := &\theta_m(z) \prod_{i=1}^r\Theta_{i,m}(z)=\sum_{k=0}^{q_2}\Theta_m(k)z^k
\end{eqnarray*}
 and put $\Phi_m(k)=0$ for $k>q_1$, $\Theta_m(k)=0$ for $k>q_2$. Then the increment sequence $$Y_m=(1-B^T)^d\prod_{i=1}^r(1-B^{Ts_i})^{d_i}X_m$$
 is periodically stationary and allows a stationary vector representation
$$\me Y_m=(1-B)^d\prod_{i=1}^r(1-B^{s_i})^{d_i}\me X_m$$ with
$$\me Y_m=(Y_{mT},Y_{mT+1},\ldots, Y_{mT+T-1})^{\top}, \quad\me X_m=(X_{mT},X_{mT+1},\ldots, X_{mT+T-1})^{\top},$$ $\boldsymbol{\eps}_m=(\eps_{mT},\eps_{mT+1},\ldots, \eps_{mT+T-1})^{\top}$.
We can write the relation
\[
\me  \Pi \me Y_m+\sum_{l=1}^{q_1^*}\me \Pi_l\me Y_{m-l}=\me \Xi \boldsymbol{\eps}_m+\sum_{l=1}^{q_2^*}\me \Xi_l\boldsymbol{\eps}_{m-l},
\]
where $\me \Pi(k,j)=\Phi_k(k-j)$, $\me \Xi(k,j)=\Theta_k(k-j)$ for $k\geq j$, $\me \Pi(k,j)=0$, $\me \Xi(k,j)=0$ otherwise.
$\me \Pi_l(k,j)=\Phi_k(lT+k-j)$, $\me \Xi_l(k,j)=\Theta_k(lT+k-j)$ \cite{Dudek-Hurd}, provided
$det(\me \Pi  +\sum_{l=1}^{q_1^*}\me \Pi_lz^l)\neq 0$ for $|z|\leq 1$ \cite{Hannan}.
A GM increment sequence is defined as
$$\chi_{\overline{\mu},\overline{s}}^{(d)}(\me X_m)=(1-B^{\mu_0})^d\prod_{i=1}^r(1-B^{s_i\mu_i})^{d_i}\me X_m,\,\, m\in \mr Z.$$
\end{pry}

\begin{thm}\label{thm1}
The structural function
$D^{(d)}_{\overline{s}}(m;\overline{\mu}_1,\overline{\mu}_2)$ of the vector-valued stochastic stationary
GM increment sequence $\chi_{\overline{\mu},\overline{s}}^{(d)}(\vec{\xi}(m))$ can be represented in the form
\begin{equation}
\label{strFnaR_vec}
 D^{(d)}_{\overline{s}}(m;\overline{\mu}_1,\overline{\mu}_2)=\int_{-\pi}^{\pi}e^{i\lambda
m} \chi_{\overline{\mu}_1}^{(d)}(e^{-i\lambda})\chi_{\overline{\mu}_2}^{(d)}(e^{i\lambda})\frac{1}
{|\beta^{(d)}(i\lambda)|^2}dF(\lambda),
\end{equation}
where $F(\lambda)$ is the matrix-valued spectral function of the stationary stochastic sequence $\chi_{\overline{\mu},\overline{s}}^{(d)}(\vec{\xi}(m))$.
The stationary GM increment sequence $\chi_{\overline{\mu},\overline{s}}^{(d)}(\vec{\xi}(m))$ admits the spectral representation
\begin{equation}
\label{SpectrPred_vec}
\chi_{\overline{\mu},\overline{s}}^{(d)}(\vec{\xi}(m))
=\int_{-\pi}^{\pi}e^{im\lambda}\chi_{\overline{\mu}}^{(d)}(e^{-i\lambda})\frac{1}{\beta^{(d)}(i\lambda)}d\vec{Z}_{\xi^{(d)}}(\lambda),
\end{equation}
where $\vec{Z}_{\xi^{(d)}}(\lambda)=\{Z_{ p}(\lambda)\}_{p=1}^{T}$ is a (vector-valued) stochastic process with uncorrelated increments on $[-\pi,\pi)$ connected with the spectral function $F(\lambda)$ by
the relation
\[
 \mt E(Z_{p}(\lambda_2)-Z_{p}(\lambda_1))(\overline{ Z_{q}(\lambda_2)-Z_{q}(\lambda_1)})
 =F_{pq}(\lambda_2)-F_{pq}(\lambda_1),\]
 \[  -\pi\leq \lambda_1<\lambda_2<\pi,\quad p,q=1,2,\dots,T.
\]
\end{thm}

\subsection{Moving average representation of  periodically stationary GM increment}
\label{moving_average}

Denote by $ H=L_2(\Omega, \cal F, \mt P)$ the Hilbert space of random variables $\zeta$ with zero first moment, ${\mt E}{\zeta}=0$, finite second moment, ${\mt E}|{\zeta}|^2<\infty$, endowed with the inner product $\langle \zeta,\eta\rangle={\mt E}{\zeta\overline{\eta}}$.
Denote by $H(\vec{\xi}^{(d)})$ the closed linear subspace of the space $H$ generated by
components
$\{\chi_{\overline{\mu},\overline{s}}^{(d)}(\xi_p(m)), p=1,\dots,T;\, m\in \mathbb Z \}$
of the stationary stochastic GM increment sequence $\vec{\xi}^{(d)}=\{\chi_{\overline{\mu},\overline{s}}^{(d)}(\xi_p(l))\}_{p=1}^{T}$, $\overline{\mu}>\overline{0}$,
and denote by $H^{q}(\vec{\xi}^{(d)})$ the closed linear subspace generated by components  $\{\chi_{\overline{\mu},\overline{s}}^{(d)}(\xi_p(m)),\, p=1,\dots,T;\, m\leqslant q \}$, $q\in \mathbb Z$.
Define a subspace
\[ S(\vec{\xi}^{(d)})=\bigcap_{q\in \mathbb{Z}} H^{q}(\vec{\xi}^{(d)})
\]
  of the Hilbert space $H(\vec{\xi}^{(d)})$. Then the  space $H(\vec{\xi}^{(d)})$ admits a decomposition
$$ H(\vec{\xi}^{(d)})=S(\vec{\xi}^{(d)})\oplus R(\vec{\xi}^{(d)}),$$
where $R(\vec{\xi}^{(d)})$ is the orthogonal complement of the subspace $S(\vec{\xi}^{(d)})$ in the space $H(\vec{\xi}^{(d)})$.
\begin{ozn}
A stationary (wide sense) stochastic GM increment sequence $\chi_{\overline{\mu},\overline{s}}^{(d)}(\vec{\xi}(m))=\{\chi_{\overline{\mu},\overline{s}}^{(d)}(\xi_p(m))\}_{p=1}^{T}$ is called regular if $H(\vec{\xi}^{(d)})=R(\vec{\xi}^{(d)})$,
and it is called singular if
$H(\vec{\xi}^{(d)})=S(\vec{\xi}^{(d)})$.
\end{ozn}

\begin{thm}
A stationary stochastic GM increment sequence $\chi_{\overline{\mu},\overline{s}}^{(d)}(\vec{\xi}(m))
=\{\chi_{\overline{\mu},\overline{s}}^{(d)}(\xi_p(m))\}_{p=1}^{T}$ is uniquely represented in the form
\begin{equation} \label{rozklad}
\chi_{\overline{\mu},\overline{s}}^{(d)}(\xi_p(m))
=\chi_{\overline{\mu},\overline{s}}^{(d)}(\xi_{S,p}(m))+\chi_{\overline{\mu},\overline{s}}^{(d)}(\xi_{R,p}(m)),
\end{equation}
where  $\chi_{\overline{\mu},\overline{s}}^{(d)}(\xi_{R,p}(m)), p=1,\dots,T$ is a regular stationary GM   increment sequence and
$\chi_{\overline{\mu},\overline{s}}^{(d)}(\xi_{S,p}(m)), p=1,\dots,T$ is a singular stationary GM increment sequence.
The GM increment sequences
$\chi_{\overline{\mu},\overline{s}}^{(d)}(\xi_{R,p}(m)), p=1,\dots,T$  and
$\chi_{\overline{\mu},\overline{s}}^{(d)}(\xi_{S,p}(m)), p=1,\dots,T$
are
orthogonal for all $ m,k\in\mathbb{Z} $. They are defined by the formulas
\[ \chi_{\overline{\mu},\overline{s}}^{(d)}(\xi_{S,p}(m))=\mt E[\chi_{\overline{\mu},\overline{s}}^{(d)}(\xi_{p}(m))|S(\vec{\xi}^{(d)})], \]
\[\chi_{\overline{\mu},\overline{s}}^{(d)}(\xi_{R,p}(m))=\chi_{\overline{\mu},\overline{s}}^{(d)}(\xi_{p}(m))
-\chi_{\overline{\mu},\overline{s}}^{(d)}(\xi_{S,p}(m)), \quad p=1,\dots,T.
\]
\end{thm}

Consider an innovation sequence  ${\vec\varepsilon(u)=\{\varepsilon_k(u)\}_{k=1}^q, u \in\mathbb Z}$ for a regular stationary GM
increment $\chi_{\overline{\mu},\overline{s}}^{(d)}(\xi_{R,p}(m)), p=1,\dots,T$, namely, a sequence of uncorrelated random
variables such that $\mathsf{E} \varepsilon_k(u)\overline{\varepsilon}_j(v)=\delta_{kj}\delta_{uv}$,   $\mathsf{E} |\varepsilon_k(u)|^2=1, k,j=1,\dots,q; u \in\mathbb Z$, and $H^{r}(\vec\xi^{(d)} )=H^{r}(\vec\varepsilon)$ holds true for all $r \in \mathbb Z$, where  $H^r(\vec\varepsilon)$ is
the Hilbert space generated by elements $ \{ \varepsilon_k(u):k=1,\dots,q; u\leq r\}$,
 $\delta_{kj}$ and $\delta_{uv}$ are Kronecker symbols.

\begin{thm}\label{thm 4}
A   stationary GM increment sequence
$\chi_{\overline{\mu},\overline{s}}^{(d)}(\vec{\xi}(m))$ is regular if and only if there exists an
innovation sequence ${\vec\varepsilon(u)=\{\varepsilon_k(u)\}_{k=1}^q, u \in\mathbb Z}$
and a sequence of matrix-valued
functions $\varphi^{(d)}(k,\overline{\mu}) =\{\varphi^{(d)}_{ij}(k,\overline{\mu}) \}_{i=\overline{1,T}}^{j=\overline{1,q}}$, $k\geq0$, such that
\[\sum_{k=0}^{\infty}
\sum_{i=1}^{T}
\sum_{j=1}^{q}
|\varphi^{(d)}_{ij}(k,\overline{\mu})|^2
<\infty,\]
\begin{equation}\label{odnostRuhSer}
\chi_{\overline{\mu},\overline{s}}^{(d)}(\vec{\xi}(m))=
\sum_{k=0}^{\infty}\varphi^{(d)}(k,\overline{\mu})\vec\varepsilon (m-k).
\end{equation}
Representation (\ref{odnostRuhSer}) is called the canonical
moving average representation of the stochastic stationary GM increment
sequence $\chi_{\overline{\mu},\overline{s}}^{(d)}(\vec{\xi}(m))$.
\end{thm}

The spectral function $F(\lambda)$ of  a stationary GM increment sequence $\chi_{\overline{\mu},\overline{s}}^{(d)}(\vec{\xi}(m))$
which admits the canonical representation
 $(\ref{odnostRuhSer})$  has the spectral density  $f(\lambda)=\{f_{ij}(\lambda)\}_{i,j=1}^T$ admitting the canonical
factorization
\begin{equation}\label{SpectrRozclad_f}
f(\lambda)=
\Phi(e^{-i\lambda})\Phi^*(e^{-i\lambda}),
\end{equation}
where the function
$\Phi(z)=\sum_{k=0}^{\infty}\varphi(k)z^k$ has
analytic in the unit circle $\{z:|z|\leq1\}$
components
$\Phi_{ij}(z)=\sum_{k=0}^{\infty}\varphi_{ij}(k)z^k; i=1,\dots,T; j=1,\dots,q$. Based on moving average  representation
$(\ref{odnostRuhSer})$ define
\[\Phi_{\overline{\mu}}(z)=\sum_{k=0}^{\infty} \varphi^{(d)}(k,\overline{\mu})z^k=\sum_{k=0}^{\infty}\varphi_{\overline{\mu}}(k)z^k.\]
 Then the following relation holds true:
\begin{equation}
\Phi_{\overline{\mu}}(e^{-i\lambda})
\Phi^*_{\overline{\mu}}(e^{-i\lambda})=
 \frac{|\chi_{\overline{\mu}}^{(d)}(e^{-i\lambda})|^2}{|\beta^{(d)}(i\lambda)|^2}f(\lambda)=
 \prod_{j=1}^r\frac{\ld|1-e^{-i\lambda\mu_js_j}\rd|^{2d_j}}{\prod_{k_j=-[s_j/2]}^{[s_j/2]}|\lambda-2\pi  k_j/s_j|^{2d_j}}f(\lambda).
\label{dd}
\end{equation}

We will use the one-sided moving average representation (\ref{odnostRuhSer}) and
relation  (\ref{dd})  for finding the mean square optimal
estimates of unobserved values of vector-valued  sequences with stationary GM increments.

\section{Hilbert space projection method of forecasting}\label{classical_extrapolation}

\subsection{Forecasting of vector-valued stationary GM increment}

Consider a vector-valued stochastic sequence with stationary GM increments $\vec{\xi}(m)$ constructed from the sequence $\zeta(m)$ with the help of transformation \eqref{PerehidXi}.
Let the stationary GM increment sequence $\chi_{\overline{\mu},\overline{s}}^{(d)}(\vec{\xi}(m))=\{\chi_{\overline{\mu},\overline{s}}^{(d)}(\xi_p(m))\}_{p=1}^{T}$
has an absolutely continuous spectral function $F(\lambda)$
and the spectral density $f(\lambda)=\{f_{ij}(\lambda)\}_{i,j=1}^{T}$. Without loss of generality we will assume that $ \mathsf{E}\chi_{\overline{\mu},\overline{s}}^{(d)}(\vec{\xi}(m))=0$ and $\overline{\mu}>\overline{0}$.

Consider the problem of mean square optimal linear estimation of the functionals
\begin{equation}
A\vec{\xi}=\sum_{k=0}^{\infty}(\vec{a}(k))^{\top}\vec{\xi}(k), \quad
A_{N}\vec{\xi}=\sum_{k=0}^{N}(\vec{a}(k))^{\top}\vec{\xi}(k),
\end{equation}
which depend on unobserved values of the stochastic sequence $\vec{\xi}(k)=\{\xi_{p}(k)\}_{p=1}^{T}$ with stationary GM
increments. Estimates are based on observations of the sequence $\vec\xi(k)$ at points $k=-1,-2,\ldots$.

We will suppose that the following conditions are satisfied:
\begin{itemize}
\item  conditions on coefficients $\vec{a}(k)=\{a_{p}(k)\}_{p=1}^{T}$, $k\geq0$, and a linear transformation $D^{\overline{\mu}}$ to be defined in Lemma \ref{lema predst A}
 \be\label{umovana a_e_d}
\sum_{k=0}^{\infty}\|\vec{a}(k)\|<\infty,\quad
\sum_{k=0}^{\infty}(k+1)\|\vec{a}(k)\|^{2}<\infty,\ee
\be\label{umovana a_mu_e_d}
\sum_{k=0}^{\infty}\|(D^{\overline{\mu}}\me a)_k\|<\infty,\quad
\sum_{k=0}^{\infty}(k+1)\|(D^{\overline{\mu}}\me a)_k\|^2<\infty,\ee

\item  \emph{the minimality condition} on the spectral density $f(\lambda)$
\be
 \ip \text{Tr}\left[ \frac{|\beta^{(d)}(i\lambda)|^2}{|\chi_{\overline{\mu}}^{(d)}(e^{-i\lambda})|^2}f^{-1}(\lambda)\right]
 d\lambda<\infty.
\label{umova11_e_d}
\ee
\end{itemize}
The latter one  is the necessary and sufficient condition under which the mean square errors of estimates of functionals $A\vec\xi$ and $A_N\vec\xi$ are not equal to $0$.

The classical Hilbert space estimation technique proposed by Kolmogorov \cite{Kolmogorov} can be described as a $3$-stage procedure:
(i) define a target element of the space $H=L_2(\Omega, \mathcal{F},\mt P)$ to be estimated,
(ii) define a subspace of $H$ generated by observations,
(iii) find an estimate of the target element as an orthogonal projection on the defined subspace.

\textbf{Stage i}. Neither  functional $A\vec{\xi}$ nor $A_{N}\vec{\xi}$   belongs to the space $H$.
With the help of the following lemma and the corresponding corollary, we  describe representations of these functionals   as sums of functionals with finite second moments  belonging to $H$ and functionals depending on observed values of the sequence $\vec\xi(k)$ (``initial values'').

\begin{lema}\label{lema predst A}
The functional $A\vec\xi$ admits a representation
\[A\vec{\xi}=B\chi\vec{\xi}-V\vec{\xi},\]
where
\[B\chi\vec{\xi}=\sum_{k=0}^{\infty}(\vec{b}(k))^{\top}\chi_{\overline{\mu},\overline{s}}^{(d)}(\vec{\xi}(k)),
\quad
V\vec{\xi}=\sum_{k=-n(\gamma)}^{-1} (\vec{v}(k))^{\top}\vec{\xi}(k),\]
\begin{eqnarray*}
\vec{v}(k)& = &\sum_{l=0}^{k+n(\gamma)} \mt{diag}_T(e_{\nu}(l-k))\vec{b}(l), \quad k=-1,-2,\dots,-n(\gamma),
\\
\vec{b}(k)& = &\sum_{m=k}^{\infty}\mt{diag}_T(d_{\overline{\mu}}(m-k))\vec{a}(m)=(D^{\overline{\mu}}{\me a})_k,  \quad k=0,1,2,\dots,
\end{eqnarray*}
 $\vec{b}(k)=(b_{1}(k),b_{2}(k),\dots,b_{T}(k))^{\top}$, $\me a=((\vec{a}(0))^{\top},(\vec{a}(1))^{\top},(\vec{a}(2))^{\top}, \ldots)^{\top}$, $\vec{v}(k)=(v_{1}(k),v_{2}(k),\dots,v_{T}(k))^{\top}$,
  $D^{\overline{\mu}}$ is a linear transformation  determined by a matrix with $T\times T$ entries
 $D^{\overline{\mu}}(k,j), k,j=0,1,2,\dots$ such that $D^{\overline{\mu}}(k,j)=\mt{diag}_T(d_{\overline{\mu}}(j-k))$ if $0\leq k\leq j $ and $D^{\overline{\mu}}(k,j)=\mt{diag}_T(0)$ for $0\leq j<k$; $\mt{diag}_T(x)$ denotes a $T\times T$ diagonal matrix with the entry $x$ on its diagonal, coefficients $\{d_{\overline{\mu}}(k):k\geq0\}$ are determined by the relationship
\[
 \sum_{k=0}^{\infty}d_{\overline{\mu}}(k)x^k
=\prod_{i=1}^r\ld(\sum_{j_i=0}^{\infty}x^{\mu_is_ij_i}\rd)^{d_i}.\]
\end{lema}
\begin{proof}
See Appendix.
\end{proof}

\begin{nas}\label{nas predst A_N}
The functional $A_N\vec\xi$ allows a representation
\[
A_N\vec\xi=B_N\chi\vec\xi-V_N\vec\xi,\]
\[
B_{N}\chi\vec{\xi}=\sum_{k=0}^{N}(\vec{b}_{N}(k))^{\top}\chi_{\overline{\mu},\overline{s}}^{(d)}(\vec{\xi}(k)),
\quad
V_{N}\vec{\xi}=\sum_{k=-n(\gamma)}^{-1} (\vec{v}_{N}(k))^{\top}\vec{\xi}(k),
\]
where the coefficients
$\vec{b}_{N}(k)=\{b_{N,p}(k)\}_{p=1}^{T}, k=0,1,\dots,N$ and
 $\vec{v}_{N}(k)=\{ v_{N,p}(k)\}_{p=1}^{T}, k=-1,-2,\dots,-n(\gamma)$
are calculated by the formulas
\begin{eqnarray}\label{koefv_N_diskr}
 \vec{v}_N(k)& = &\sum_{l=0}^{N \wedge k+n(\gamma)} \mt{diag}_T(e_{\nu}(l-k))\vec{b}_N(l), \, k=-1,-2,\dots,-n(\gamma),
\\
 \label{determ_b_N}
\vec{b}_N(k)& = &\sum_{m=k}^N\mt{diag}_T(d_{\overline{\mu}}(m-k))\vec{a}(m)=(D^{\overline{\mu}}_{N}{\me a}_{N})_k, \, k=0,1,\dots,N,
\end{eqnarray}
  $D^{\overline{\mu}}_{N}$ is the  linear transformation  determined by an infinite matrix with the entries $(D^{\overline{\mu}}_{N})(k,j)=\mt{diag}_T(d_{\overline{\mu}}(j-k))$ if
$0\leq k\leq j\leq N$, and $(D^{\overline{\mu}}_{N})(k,j)=0$ if $j<k$ or $j,k>N$; $\me a_N=((\vec{a}(0))^{\top},(\vec{a}(1))^{\top}, \ldots,(\vec{a}(N))^{\top},\vec 0 \ldots)^{\top}$.
\end{nas}

So that, Lemma \ref{lema predst A} provides a representation of the functional $A\vec\xi$ as a sum of an element $B\chi\vec\xi$ from the space $H=L_2(\Omega, \mathcal{F},\mt P)$ under conditions (\ref{umovana a_e_d}) -- (\ref{umovana a_mu_e_d}) and linear combination $V\vec\xi$ of a finite number of initial values  $\vec\xi(k)$, $k=-1,-2,\ldots,-n(\gamma)$, which are observed. Thus, the following
equality hold true
\begin{equation}\label{mainformula_e_d}
 \widehat{A}\vec\xi=\widehat{B}\chi\vec\xi-V\vec\xi.
\end{equation}
Denote by
$\Delta(f,\widehat{A}\vec\xi):=\mathsf{E} |A\vec\xi-\widehat{A}\vec\xi |^2$
 the mean square error of the optimal
estimate $\widehat{A}\vec\xi$ of the functional $A\vec\xi$ and let
$\Delta(f,\widehat{B}\chi\vec\xi):=\mathsf{E}  |B\chi\vec\xi-\widehat{B}\chi\vec\xi |^2$
denote the mean square error of the optimal estimate $\widehat{B}\chi\vec\xi$ of the functional $B\chi\vec\xi$. Then
\begin{eqnarray*}
\Delta\left(f;\widehat{A}\vec\xi\right)
&=&\mt E \ld|A\vec\xi-\widehat{A}\vec\xi\rd|^2= \mt
E\ld|B\chi\vec\xi-V\vec\xi-\widehat{B}\chi\vec\xi+V\vec\xi\rd|^2
\\
&=&\mt E\ld|B\chi\vec\xi-\widehat{B}\chi\vec\xi\rd|^2=\Delta\ld(f;\widehat{B}\chi\vec\xi\rd).
\end{eqnarray*}
Thus, we have defined the functional $B\chi\vec\xi$ to be used in finding
the optimal linear estimate of the functional $A\vec\xi$.

At \textbf{stage ii}, we recall the subspace $H^{0-}(\vec\xi^{(d)}):=H^{-1}(\vec\xi^{(d)})$ of the Hilbert space $H=L_2(\Omega,\mathcal{F},\mt P)$
defined in  Subsection \ref{moving_average}, which is generated by observations $\{\chi_{\overline{\mu},\overline{s}}^{(d)}(\xi_p(k)), p=1,\dots,T;\, k\leq -1 \}$.
Denote by $L_2^{0-}(f)$ the closed linear subspace of the Hilbert space $L_2(f)$ of vector-valued functions endowed with the inner product $\langle g_1;g_2\rangle=\ip (g_1(\lambda))^{\top}f(\lambda)\overline{g_2(\lambda)}d\lambda$ which is generated by functions
\[ e^{i\lambda k}\chi_{\overline{\mu}}^{(d)}(e^{-i\lambda})(\beta^{(d)}(i\lambda))^{-1}\boldsymbol{\delta}_{l}, \quad \boldsymbol{\delta}_l=\{\delta_{lp}\}_{p=1}^{T}, \, l=1,\dots,T; k \leqslant -1,
\]
where $\delta_{lp}$ are Kronecker symbols.
The relation
\begin{equation}\label{spectr_zobr_zeta}
 \chi_{\overline{\mu},\overline{s}}^{(d)}(\xi_p(m))=\ip e^{i\lambda
 m}\chi_{\overline{\mu}}^{(d)}(e^{-i\lambda})\dfrac{1}{\beta^{(d)}(i\lambda)}dZ_p(\lambda), \quad p=1,2,\dots,T,
\end{equation}
implies  a relation between elements
$\chi_{\overline{\mu},\overline{s}}^{(d)}(\xi_p(m))$ of the space
$H (\vec\xi^{(d)})$ and elements
$e^{i\lambda m}\chi_{\overline{\mu}}^{(d)}(e^{-i\lambda})(\beta^{(d)}(i\lambda))^{-1}\boldsymbol{\delta}_{p}$ of the space $L_2 (f)$.
The spectral representation of the functional $B\chi\vec{\xi}$ can be written in the form
\[
B\chi\vec{\xi}=\int_{-\pi}^{\pi}\left(\vec{B}_{\overline{\mu}}(e^{i\lambda})\right)^{\top}
\frac{\chi_{\overline{\mu}}^{(d)}(e^{-i\lambda})}{\beta^{(d)}(i\lambda)}d\vec{Z}_{\xi^{(d)}}(\lambda),
 \]
 where
\[
 \vec B_{\overline{\mu}}(e^{i\lambda})=\sum_{k=0}^{\infty}{\vec b}(k)e^{i\lambda k}
 =\sum_{k=0}^{\infty}(D^{\overline{\mu}}\me a)_ke^{i\lambda k}.\]
Thus, at \textbf{stage iii}, the problem is equivalent to   finding  a projection of the element $ \vec{B}_{\overline{\mu}}(e^{i\lambda})
\frac{\chi_{\overline{\mu}}^{(d)}(e^{-i\lambda})}{\beta^{(d)}(i\lambda)}$ of the Hilbert space $L_2(f)$ on the subspace $L_2^{0-}(f)$.

Relation (\ref{mainformula_e_d}) implies that every linear estimate $\widehat{A}\vec\xi$ of the functional $A\vec\xi$
can be written in the form
\be \label{otsinka A_e_d}
 \widehat{A}\vec\xi=\ip
(\vec{h}_{\overline{\mu}}(\lambda))^{\top}d\vec{Z}_{\xi^{(d)}}(\lambda)-\sum_{k=-n(\nu)}^{-1}(\vec v(k))^{\top}\vec\xi(k),
\ee
 where
$\vec{h}_{\overline{\mu}}(\lambda)=\{h_{p}(\lambda)\}_{p=1}^{T}$ is the spectral characteristic of the estimate $\widehat{B}\vec\xi$, which is a projection of the element
$\vec B_{\overline{\mu}}(e^{i\lambda })\frac{\chi_{\overline{\mu}}^{(d)}(e^{-i\lambda}) }{\beta^{(d)}(i\lambda)}$ on the subspace $L_2^{0-}(f)$.
This estimate is characterized by the following conditions:
\begin{equation}
\label{umova1}
\vec{h}_{\overline{\mu}}(\lambda)\in L_2^{0-}\ld(f\rd),
\end{equation}
\begin{equation}
\label{umova2}
\ld(\vec{B}_{\overline{\mu}}(e^{i\lambda})\frac{\chi_{\overline{\mu}}^{(d)}(e^{-i\lambda}) }{\beta^{(d)}(i\lambda)}-\vec{h}_{\overline{\mu}}(\lambda)\rd)\perp
L_2^{0-}\ld(f\rd).
\end{equation}
Condition \eqref{umova2} implies the following
relation holding true for all $  k \leqslant -1$
\begin{gather}
\label{umova2_2}
 \int_{-\pi}^{\pi}\left(\vec{B}_{\overline{\mu}}(e^{i\lambda}) \frac{\chi_{\overline{\mu}}^{(d)}(e^{-i\lambda})}{\beta^{(d)}(i\lambda)}-\vec{h}_{\overline{\mu}}(\lambda)\right)^{\top} f(\lambda)e^{-i\lambda k}\frac{\chi_{\overline{\mu}}^{(d)}(e^{i\lambda})}{\overline{\beta^{(d)}(i\lambda)}} d\lambda=\vec 0.
\end{gather}
Thus, the spectral characteristic of the estimate
$\widehat{B}\chi\vec\xi$ can be represented in the form
\[ (\vec{h}_{\overline{\mu}}(\lambda))^{\top}=(\vec{B}_{\overline{\mu}}(e^{i\lambda}))^{\top}
\frac{\chi_{\overline{\mu}}^{(d)}(e^{-i\lambda})}{\beta^{(d)}(i\lambda)}-
\frac{\overline{\beta^{(d)}(i\lambda)}}{\chi_{\overline{\mu}}^{(d)}(e^{i\lambda})}(\vec{C}_{\overline{\mu}}(e^{i\lambda}))^{\top}f^{-1}(\lambda),\]
where
\[
\vec{C}_{\overline{\mu}}(e^{i \lambda})=\sum_{k=0}^{\infty}\vec{c}_{\overline{\mu}}(k)e^{i\lambda k},
\]
and $\vec{c}(k)=\{c_p(k)\}_{p=1}^T, k\geqslant 0$  are unknown coefficients to be found.

Condition \eqref{umova1}  implies the following representation of the
spectral characteristic $\vec{h}_{\overline{\mu}}(\lambda)$
\[
 \vec{h}_{\overline{\mu}}(\lambda)=\vec{h} (\lambda)\chi_{\overline{\mu}}^{(d)}(e^{-i\lambda})
 \frac{1}{\beta^{(d)}(i\lambda)}, \quad \vec{h} (\lambda)=
 \sum_{k=1}^{\infty}\vec{s}(k)e^{-i\lambda k},\]
which allows us to write the relations
\begin{equation} \label{eq_C}
\int_{-\pi}^{\pi} \left[(\vec{B}_{\overline{\mu}}(e^{i\lambda}))^{\top}-
\frac{|\beta^{(d)}(i\lambda)|^2}{\chi_{\overline{\mu}}^{(d)}(e^{-i\lambda}) \chi_{\overline{\mu}}^{(d)}(e^{i\lambda})}(\vec{C}_{\overline{\mu}}(e^{i\lambda}))^{\top}f^{-1}(\lambda)\right]e^{-ij\lambda}d\lambda=\vec 0,  j\geqslant 0.
\end{equation}

\noindent Next we define the matrix-valued Fourier coefficients
\begin{equation} \label{Fourier_F}
 F^{\overline{\mu}}(k,j)=\frac{1}{2\pi}\ip
 e^{i\lambda(j-k)}\frac{|\beta^{(d)}(i\lambda)|^2}{| \chi_{\overline{\mu}}^{(d)}(e^{-i\lambda})|^2}f^{-1}(\lambda)d\lambda,\quad k,j\geq0,
\end{equation}
and rewrite relation \eqref{eq_C}  as a system of linear vector equations
\[
 \vec{b}(j)=\sum_{k=0}^{\infty}F^{\overline{\mu}}(j,k)\,\,{\vec c}_{\overline{\mu}}(k),\quad j\geq0,
 \]
 determining the unknown coefficients ${\vec c}_{\overline{\mu}}(k)$, $k\geq0$.
This system can be presented in the matrix form
\be\label{meq-e}
 D^{\overline{\mu}}\me a=\me F_{\overline{\mu}}\me c_{\overline{\mu}},
\ee
where
\[\me c_{\overline{\mu} }=((\vec{c}_{\overline{\mu}}(0))^{\top},(\vec{c}_{\overline{\mu}}(1))^{\top},(\vec{c}_{\overline{\mu}}(2))^{\top}, \ldots)^{\top},\, \me a=((\vec{a}(0))^{\top},(\vec{a}(1))^{\top},(\vec{a}(2))^{\top}, \ldots)^{\top},\]
$\me F_{\overline{\mu}}$ is a linear operator in the space $\ell_2$ which is determined by a matrix with the $T\times T$ matrix entries $\me F_{\overline{\mu}}(j,k)=F^{\overline{\mu}}(j,k)$, $j,k\geq0$; the linear transformation $D^{\overline{\mu}}$ is defined in Lemma \ref{lema predst A}.

To show that operator $\me F_{\overline{\mu}}$ is invertible we note that the problem of projection of the element ${B}\vec\xi$ of the Hilbert space $H$ on the closed convex set $H^{0-}(\vec\xi^{(d)}_{\overline{\mu}})$ has a unique solution for each non-zero coefficients $\{\vec{a}(0),\vec{a}(1)),\vec{a}(2), \ldots\}$, satisfying conditions
(\ref{umovana a_e_d}) -- (\ref{umovana a_mu_e_d}). Therefore, equation (\ref{meq-e}) has a unique solution for each vector $D^{\overline{\mu}} \me a$, which implies existence of the inverse operator $\me F^{-1}_{\overline{\mu}}$.

Therefore,  coefficients $\vec c_{\overline{\mu}}(k)$, $k\geq0$, which determine the spectral characteristic $\vec{h}_{\overline{\mu}}(\lambda)$, can be calculated as
\begin{equation}
 \vec c_{\mu }(k)= (\me F_{\overline{\mu}}^{-1}D^{\overline{\mu}}\me a  )_k,\quad k\geq 0,
 \end{equation}
 where $(\me F_{\overline{\mu}}^{-1}D^{\overline{\mu}}\me a)_k$, $k\geq 0$, is the
 $k$th $T$-dimension vector element of the vector $\me F_{\overline{\mu}}^{-1}D^{\overline{\mu}}\me a$.

The spectral characteristic $\vec{h}_{\overline{\mu}}(\lambda)$ of the  estimate $\widehat{B}\chi\vec\xi$ is calculated by the formula
\begin{equation}\label{spectr A_e_d}
 (\vec{h}_{\overline{\mu}}(\lambda))^{\top}=(\vec{B}_{\overline{\mu}}(e^{i\lambda}))^{\top}
\frac{\chi_{\overline{\mu}}^{(d)}(e^{-i\lambda})}{\beta^{(d)}(i\lambda)}-
\frac{\overline{\beta^{(d)}(i\lambda)}
}
{\chi_{\overline{\mu}}^{(d)}(e^{i\lambda})}\left(
\sum_{k=0}^{\infty}(\me F_{\overline{\mu}}^{-1}D^{\overline{\mu}}\me a )_k e^{i\lambda k}\right)^{\top}f^{-1}(\lambda).
\end{equation}

The value of the mean square error of the estimate $\widehat{A}\vec\xi$ is calculated by the formula
\begin{eqnarray}
\nonumber
\Delta\ld(f;\widehat{A}\vec\xi\rd)
& = & \Delta\ld(f;\widehat{B}\chi\vec\xi\,\rd) =  \mt E\ld|B\chi\vec\xi-\widehat{B}\chi\vec\xi\,\rd|^2
\\
\nonumber
& = &
\frac{1}{2\pi}\int_{-\pi}^{\pi}\frac{\overline{\beta^{(d)}(i\lambda)}
}
{\chi_{\overline{\mu}}^{(d)}(e^{i\lambda})}\left(
\sum_{k=0}^{\infty}(\me F_{\overline{\mu}}^{-1}D^{\overline{\mu}}\me a )_k e^{i\lambda k}\right)^{\top}f(\lambda)
\\\nonumber
& \quad & \times \left(\overline{
\sum_{k=0}^{\infty}(\me F_{\overline{\mu}}^{-1}D^{\overline{\mu}}\me a )_k e^{i\lambda k}}\right) \frac{\beta^{(d)}(i\lambda)}
{\chi_{\overline{\mu}}^{(d)}(e^{-i\lambda})}d\lambda
\\
& = &\ld\langle D^{\overline{\mu}}\me a,\me F_{\overline{\mu}}^{-1}D^{\overline{\mu}}\me a\rd\rangle.\label{pohybka}
\end{eqnarray}

Next consider the problem in the case where the GM  incremental sequence of the stochastic sequence $\vec\xi(m)$ admits   moving-average representation \eqref{odnostRuhSer} and its spectral density $f(\lambda)=\{f_{ij}(\lambda)\}_{i,j=1}^T$   admits the canonical factorization \eqref{SpectrRozclad_f}, \eqref{dd}, namely
\be \label{fact_f_e_d}
 f(\lambda)=\Phi(e^{-i\lambda})\Phi^*(e^{-i\lambda}),
  \quad
 \frac{ | \chi_{\overline{\mu}}^{(d)}(e^{-i\lambda}) |^2}{ |\beta^{(d)}(i\lambda) |^2}f(\lambda)=
 \Phi_{\overline{\mu}}(e^{-i\lambda})\Phi^*_{\overline{\mu}}(e^{-i\lambda}),
\ee
where
\[
\Phi(e^{-i\lambda})=\sum_{k=0}^{\infty}\varphi(k)e^{-i\lambda k},\quad \Phi_{\overline{\mu}}(e^{-i\lambda})=\sum_{k=0}^{\infty}\varphi_{\overline{\mu}}(k)e^{-i\lambda k},
\]
and
$\varphi_{\overline{\mu}}(k)=\{\varphi_{ij}(k)\}_{i=\overline{1,T}}^{j=\overline{1,q}},\, k=0,1,2,\dots$.
Let $E_q$ denote the identity $q\times q$ matrix. Define the matrix-valued function $\Psi_{\overline{\mu}}(e^{-i\lambda})= \{\Psi_{ij}(e^{-i\lambda})\}_{i=\overline{1,q}}^{j=\overline{1,T}}$ by the equation
\[\Psi_{\overline{\mu}}(e^{-i\lambda})\Phi_{\overline{\mu}}(e^{-i\lambda})=E_q.\]

Formulas for calculating the spectral characteristic $\vec h_{\overline{\mu}}(\lambda)$ and the value of the mean square error $\Delta(f;\widehat{A}\vec\xi)$ of the estimate $\widehat{A}\vec\xi$ can be presented in terms of the function $\Psi_{\overline{\mu}}(e^{-i\lambda})$  and the factorization coefficients $\varphi_{\overline{\mu}}(k)$, $ k=0,1,2,\dots$.
One can directly check that conditions (\ref{umova1}) and (\ref{umova2}) are satisfied by the function
\be\label{simple_spectr A_e_d}
 \vec h_{\overline{\mu}}(\lambda)=\frac{\chi_{\overline{\mu}}^{(d)}(e^{-i\lambda}) }{\beta^{(d)}(i\lambda)}\ld(
 \vec B_\mu(e^{i\lambda})-(\Psi_{\overline{\mu}}(e^{-i\lambda}))^{\top}
 \vec r_{\overline{\mu}}(e^{i\lambda})
  \rd),
 \ee
 where
\begin{eqnarray*}
\vec r_{\overline{\mu}}(e^{i\lambda})&=&
 \sum_{k=0}^{\infty}(D^{\overline{\mu}}\me A\varphi_{\overline{\mu}})_ke^{i\lambda k},
\\
 (D^{\overline{\mu}}\me A\varphi_{\overline{\mu}})_k &=&  \sum_{m=0}^{\infty}\sum_{l=m}^{\infty}(\varphi_{\overline{\mu}}(m))^{\top} D^{\overline{\mu}}(m,l)\vec a(l+k)
\\
&=& \sum_{m=0}^{\infty}\sum_{l=k}^{\infty}(\varphi_{\overline{\mu}}(m))^{\top}\vec a(m+l)d_{\overline{\mu}}(l-k),
\end{eqnarray*}
 and
$\me A$ is a linear symmetric operator which is determined by a matrix with the entries
$\me A(k,j)=\vec a(k+j)$, $k,j\geq0$.   The defined operators $D^{\overline{\mu}}\me A$ and $\me A$ are compact under   conditions (\ref{umovana a_e_d}) -- (\ref{umovana a_mu_e_d}).
Then the value of the mean square error is calculated by the formula
\begin{eqnarray}
\nonumber
\Delta\ld(f;\widehat{A}\vec\xi\rd)&=&
\frac{1}{2\pi}\ip
\left(
 \sum_{k=0}^{\infty}(D^{\overline{\mu}}\me A\varphi_{\overline{\mu}})_ke^{i\lambda k}
\right)^{\top}
\left(\overline{
 \sum_{k=0}^{\infty}(D^{\overline{\mu}}\me A\varphi_{\overline{\mu}})_ke^{i\lambda k}}
\right)d\lambda
\\
&=&\frac{1}{2\pi}\ip\ld\|\vec r_{\overline{\mu}}(e^{i\lambda})\rd\|^2
d\lambda=
\ld\|D^{\overline{\mu}}\me A\varphi_{\overline{\mu}}\rd\|^2.
\label{nsimple_poh A_e_d}
\end{eqnarray}

The derived results are summarized in the following theorem.

\begin{thm}
\label{thm_est_A}
Let a vector-valued stochastic sequence $\{\vec{\xi}(m), m\in\mathbb Z\} $ determine a
stationary stochastic GM increment sequence
$\chi_{\overline{\mu},\overline{s}}^{(n)}(\vec{\xi}(m))$ with the spectral density matrix $f(\lambda)=\{f_{ij}(\lambda)\}_{i,j=1}^{T}$ which satisfies the minimality condition (\ref{umova11_e_d}).
Let coefficients $\vec {a}(j), j\geqslant 0$ satisfy conditions  (\ref{umovana a_e_d}) -- (\ref{umovana a_mu_e_d}).
\\
Then the optimal linear estimate $\widehat{A}\vec\xi$ of the functional $A\vec\xi$ based on observations of the sequence
$\vec\xi(m)$ at points $m=-1,-2,\ldots$ is calculated by formula (\ref{otsinka A_e_d}).
The spectral characteristic
$\vec h_{\overline{\mu}}(\lambda)=\{h_{p}(\lambda)\}_{p=1}^{T}$ and the value of the mean square error $\Delta(f;\widehat{A}\vec\xi)$ of the   estimate $\widehat{A}\vec\xi$ are calculated by formulas
(\ref{spectr A_e_d}) and (\ref{pohybka}) respectively.
\\
In the case where the spectral density $f(\lambda)$ admits the canonical factorization (\ref{fact_f_e_d}) the spectral characteristic and the value of the mean square error of the optimal estimate $\widehat{A}\xi$ can be calculated by formulas (\ref{simple_spectr A_e_d}) and (\ref{nsimple_poh A_e_d}) respectively.
\end{thm}

\subsection{Estimates of functional $A_N\vec\xi$ and value $\xi_p(N)$}

Theorem \ref{thm_est_A} allows us to find the optimal estimate
$\widehat{A}_N\vec\xi$ of the functional $A_N\vec\xi$ which depends on the unobserved values $\vec\xi(m)$, $m=0,1,2,\ldots,N$,
based on observations of the sequence $\vec\xi(m)$ at points $m=-1,-2,\ldots$.
Put $\vec a(k)=0$ for $k>N$.
Then we get that the spectral characteristic $\vec h_{\mu, N}(\lambda)$ of the optimal estimate

\begin{equation} \label{estimation_A_N}
\widehat{A}_{N}\vec{\xi}=\int_{-\pi}^{\pi}
(\vec{h}_{\overline{\mu},N}(\lambda))^{\top}d\vec{Z}_{\xi^{(d)}}(\lambda) - \sum_{k=-n(\gamma)}^{-1} (\vec{v}_{N}(k))^{\top}\vec{\xi}(k),
\end{equation}
is calculated by the formula
\begin{eqnarray}
\nonumber
(\vec{h}_{\overline{\mu},N}(\lambda))^{\top}
 & = & (\vec{B}_{\overline{\mu},N}(e^{i\lambda}))^{\top}
\frac{\chi_{\overline{\mu}}^{(d)}(e^{-i\lambda})}{\beta^{(d)}(i\lambda)}
\\
 & \quad & -
\frac{\overline{\beta^{(d)}(i\lambda)}
}
{\chi_{\overline{\mu}}^{(d)}(e^{i\lambda})}\left(
\sum_{k=0}^{\infty}(\me F_{\overline{\mu}}^{-1}D^{\overline{\mu}}_{N}{\me a}_N )_k e^{i\lambda k}\right)^{\top}f^{-1}(\lambda).\label{est_h_N}
\end{eqnarray}
where
\[
B_{\overline{\mu},N}(e^{i\lambda })=\sum_{k=0}^{N}(D^{\overline{\mu}}_N{\me a}_N)_ke^{i\lambda k},
\]
and   $D^{\overline{\mu}}_N$ is defined in Corollary \ref{nas predst A_N}.
 The value of the mean square error of the  estimate $\widehat{A}_N\xi$ is
\begin{eqnarray}
\nonumber
\Delta \ld(f, \widehat{A}_{N}\vec{\xi}\rd)& = &\Delta \ld(f, \widehat{B}_{N}\chi\vec{\xi}\rd)={\mathsf E}\ld|B_{N}\chi\vec{\xi}-\widehat{B}_{N}\chi\vec{\xi}\rd|^2
\\
\nonumber
& = &\frac{1}{2\pi}\int_{-\pi}^{\pi}\frac{\overline{\beta^{(d)}(i\lambda)}
}
{\chi_{\overline{\mu}}^{(d)}(e^{i\lambda})}\left(
\sum_{k=0}^{\infty}(\me F_{\overline{\mu}}^{-1}D^{\overline{\mu}}_{N}{\me a}_N )_k e^{i\lambda k}\right)^{\top}f(\lambda)
\\\nonumber
& \quad &\times\left(\overline{
\sum_{k=0}^{\infty}(\me F_{\overline{\mu}}^{-1}D^{\overline{\mu}}_{N}{\me a}_N )_k e^{i\lambda k}}\right)
\frac{\beta^{(d)}(i\lambda)}
{\chi_{\overline{\mu}}^{(d)}(e^{-i\lambda})}d\lambda
\\
& = &\ld\langle D^{\overline{\mu}}_{N}{\me a}_N,\me F_{\overline{\mu}}^{-1}D^{\overline{\mu}}_{N}{\me a}_N\rd\rangle.\label{pohybka_N}
\end{eqnarray}

In the case where the spectral density $f(\lambda)$ admits the canonical factorization (\ref{fact_f_e_d}) the spectral characteristic can be calculated as

\be\label{simple_spectr A_N_e_d}
 \vec h_{\overline{\mu},N}(\lambda)=\frac{\chi_{\overline{\mu}}^{(d)}(e^{-i\lambda})}{\beta^{(d)}(i\lambda)}
 \left( \vec B_{\overline{\mu},N}(e^{i\lambda})-(\Psi_{\overline{\mu}}(e^{-i\lambda}))^{\top}
 \vec r_{\overline{\mu},N}(e^{i\lambda})\right)
 \ee
where
\begin{eqnarray*}
 \vec r_{\overline{\mu},N}(e^{i\lambda}) & = &
 \sum_{k=0}^{N}(\widetilde{D}_N^{\overline{\mu}}\me A_N\varphi_{\overline{\mu},N})_ke^{i\lambda k},
\\
(\widetilde{D}_N^{\overline{\mu}}\me A_N\varphi_{\overline{\mu},N})_k & = & \sum_{m=0}^{N}\sum_{l=k}^{N}(\varphi_{\overline{\mu}}(m))^{\top}\vec a(m+l)d_{\overline{\mu}}(l-k),
\end{eqnarray*}
 and $\varphi_{\overline{\mu},N}=(\varphi_{\overline{\mu}}(0),\varphi_{\overline{\mu}}(1),\ldots,\varphi_{\overline{\mu}}(N))$;
 $\me A_N$ is a linear operator determined by the coefficients $\vec a(k)$, $k=0,1,\ldots, N$, as follows: $(\me A_N)(k,j)=\vec a(k+j)$, $0\leq k+j\leq N$, $(\me A_N)(k,j)=0$, $k+j>N$, $0\leq k,j\leq N$;
 $\widetilde{D}^{\overline{\mu}}_{N}$ is a matrix of the dimension $(N+1)\times(N+1)$ determined by the $T\times T$ entries $\widetilde{D}^{\overline{\mu}}_{N}(k,j)=\mt{diag}_T(d_{\overline{\mu}}(j-k))$ if
$0\leq k\leq j\leq N$ and $\widetilde{D}^{\overline{\mu}}_{N}(k,j)=\mt{diag}_T(0)$ if $0\leq j<k\leq N$.

The value of the mean square error is calculated by the formula
\begin{eqnarray}
\Delta\ld(f;\widehat{A}_N\vec \xi\rd)
\nonumber
& = &
\frac{1}{2\pi}\ip
\left(
\sum_{k=0}^{N}(\widetilde{D}_N^{\overline{\mu}}\me A_N\varphi_{\overline{\mu},N})_ke^{i\lambda k}
\right)^{\top}
\\\nonumber
& \quad &\times
\left(\overline{
\sum_{k=0}^{N}(\widetilde{D}_N^{\overline{\mu}}\me A_N\varphi_{\overline{\mu},N})_ke^{i\lambda k}}
\right)d\lambda
\\
& = &\frac{1}{2\pi}\ip\ld\|\vec r_{\overline{\mu},N}(e^{i\lambda})\rd\|^2
d\lambda=
\ld\|\widetilde{D}^{\overline{\mu}}_N\me A_N\varphi_{\overline{\mu},N}\rd\|^2.
\label{simple_poh A_N_e_d}
\end{eqnarray}

Thus, the following theorem holds true.

\begin{thm}
\label{thm_est_A_N}
Let $\{\vec{\xi}(m), m\in\mathbb Z\}$ be a stochastic sequence which determine a stationary stochastic GM increment sequence
$\chi_{\overline{\mu},\overline{s}}^{(n)}(\vec{\xi}(m))$ with the spectral density matrix $f(\lambda)$ which satisfies the minimality condition (\ref{umova11_e_d}).
The optimal linear estimate $\widehat{A}_N\vec\xi$ of the functional $A_N\vec\xi$ based on observations of the sequence
$\vec\xi(m)$ at points $m=-1,-2,\ldots$ is calculated by formula \eqref{estimation_A_N}.
The spectral characteristic $\vec{h}_{\overline{\mu},N}(\lambda)=\{h_{\mu,N,p}(\lambda)\}_{p=1}^{T}$ and the value of the mean square error $\Delta(f;\widehat{A}_N\vec \xi)$
are calculated by formulas  \eqref{est_h_N} and \eqref{pohybka_N} respectively.
In the case where the spectral density $f(\lambda)$ admits the canonical factorization (\ref{fact_f_e_d}) the spectral characteristic
$\vec h_{\overline{\mu},N}(\lambda)$ and the value of the mean square error of the  estimate $\widehat{A}_N\vec\xi$ can be calculated by formulas
(\ref{simple_spectr A_N_e_d}) and (\ref{simple_poh A_N_e_d}) respectively.
\end{thm}

For the problem of the mean square optimal estimate of the unobserved value
$A_{N,p}\vec\xi=\xi_p(N)=\vec\xi(N)\boldsymbol{\delta}_p$, $p=1,2,\dots,T$, $N\geq0$
of the stochastic sequence $\vec\xi(m)$ with GM stationary increments based on its observations  at points $m=-1,-2,\ldots$
we have the following corollary from Theorem \ref{thm_est_A_N}.

\begin{nas}\label{nas xi_e_d}
The optimal linear estimate $\widehat{\xi}_p(N)$ of the   value
$\xi_p(N)$, $p=1,2,\dots,T$, $N\geq0$, of the stochastic sequence with GM stationary increments from observations of the sequence $\vec\xi(m)$ at points $m=-1,-2,\ldots$ is calculated by formula
\begin{equation} \label{est_xi_N}
 \widehat{\xi}_p(N)=\ip
\left(\vec h_{\overline{\mu},N,p}(\lambda)\right)^{\top}d\vec Z_{\xi^{(d)}}(\lambda)- \sum_{k=-n(\gamma)}^{-1} (\vec{v}_{N,p}(k))^{\top}\vec{\xi}(k).
\end{equation}
 The spectral characteristic $\vec h_{\overline{\mu},N,p}(\lambda)$ of the estimate is calculated by the formula
\begin{eqnarray}
\nonumber
\ld(\vec{h}_{\overline{\mu},N,p}(\lambda)\rd)^{\top}
& = &\frac{\chi_{\overline{\mu}}^{(d)}(e^{-i\lambda})}{\beta^{(d)}(i\lambda)}\left(\boldsymbol{\delta}_p\sum_{k=0}^Nd_{\overline{\mu}}(N-k)e^{i\lambda k}\right)^{\top}
\\
& \quad &-
\frac{\overline{\beta^{(d)}(i\lambda)}
}
{\chi_{\overline{\mu}}^{(d)}(e^{i\lambda})}
\left(
\sum_{k=0}^{\infty}
(\me F_{\overline{\mu}}^{-1}
\me d_{\overline{\mu},N,p}
)_k e^{i\lambda k}\right)^{\top}f^{-1}(\lambda).\label{sph_est_xi_N}
\end{eqnarray}
where
$\me d_{\overline{\mu},N,p}=(d_{\overline{\mu}}(N)\boldsymbol{\delta}_p^{\top},d_{\overline{\mu}}(N-1)\boldsymbol{\delta}_p^{\top},\ldots,d_{\overline{\mu}}(0)\boldsymbol{\delta}_p^{\top},0,\ldots)^{\top}$.
The value of the mean square error of the   estimate $\widehat{\xi}_p(N)$ is calculated by the formula
\begin{eqnarray}
\nonumber
 \Delta\ld(f;\widehat{\xi}_p(N)\rd)
 & = & \Delta\ld(f;\chi_{\overline{\mu},\overline{s}}^{(n)}(\widehat{\xi}_p(m))\rd)
 =\mt E\ld|\chi_{\overline{\mu},\overline{s}}^{(n)}(\xi_p(m))-\chi_{\overline{\mu},\overline{s}}^{(n)}(\widehat{\xi}_p(m))\rd|^2 \\
\nonumber
& = &\frac{1}{2\pi}\int_{-\pi}^{\pi}\frac{\overline{\beta^{(d)}(i\lambda)}
}
{\chi_{\overline{\mu}}^{(d)}(e^{i\lambda})}\left(
\sum_{k=0}^{\infty}(\me F_{\overline{\mu}}^{-1}
\me d_{\overline{\mu},N,p} )_k e^{i\lambda k}\right)^{\top}f(\lambda)
\\\nonumber
& \quad &\times\left(\overline{
\sum_{k=0}^{\infty}(\me F_{\overline{\mu}}^{-1}
\me d_{\overline{\mu},N,p} )_k e^{i\lambda k}}\right)
\frac{\beta^{(d)}(i\lambda)}
{\chi_{\overline{\mu}}^{(d)}(e^{-i\lambda})}d\lambda
\\
\label{delta_est_xi_N}
& = &\ld\langle
\me d_{\overline{\mu},N,p},\me F_{\overline{\mu}}^{-1}
\me d_{\overline{\mu},N,p}\rd\rangle.
\end{eqnarray}
In the case where the spectral density $f(\lambda)$ admits   canonical factorization (\ref{fact_f_e_d}), and the condition  $\min_{i=\overline{1,r}}\mu_is_i>N$ is satisfied, the spectral characteristic
and the value of the mean square error of the  estimate $\widehat{\xi}_p(N)$ can be calculated by the formulas
\be
\label{sph_f_est_xi_N}
\vec{h}_{\overline{\mu},N,p}(\lambda)=
\frac{\chi_{\overline{\mu}}^{(d)}(e^{-i\lambda})}{\beta^{(d)}(i\lambda)}
 e^{iN\lambda}
 \left[\boldsymbol{\delta}_p
-(\Psi_{\overline{\mu}}(e^{-i\lambda}))^{\top}
\left(\sum_{k=0}^{N}\varphi_{\overline{\mu}}(k)e^{-i\lambda k}\right)^{\top}\boldsymbol{\delta}_p\right]
\ee
and
\begin{eqnarray}
\nonumber
 \Delta\ld(f;\widehat{\xi}_p(N)\rd) & = &
 \frac{1}{2\pi}\ip
  \left[
  \left(\boldsymbol{\delta}_p\right)^{\top}
 \sum_{k=0}^{N}\varphi_{\overline{\mu}}(k)e^{-i\lambda k}\right]
   \left[
  \left(\boldsymbol{\delta}_p\right)^{\top}
 \sum_{k=0}^{N}\varphi_{\overline{\mu}}(k)e^{-i\lambda k}\right]^{*}
  d\lambda
  \\
  & = &\sum_{k=0}^{N}\sum_{j=1}^{q}
  |\varphi_{\mu,p,j}(k)|^2.\label{delta_f_est_xi_N}
\end{eqnarray}
\end{nas}

\begin{zau}\label{zauv_operator_W}
Since for all $\overline{d}\geq\overline{1}$ and $\overline{\mu}\geq\overline{1}$ the condition
\[\ip\left |\ln \frac{| \chi_{\overline{\mu}}^{(d)}(e^{-i\lambda})|^2}{|\beta^{(d)}(i\lambda)|^2}\right |d\lambda <\infty\]
holds true, there exists a function
\[w_{\overline{\mu}}(z)=\sum_{k=0}^{\infty}w_{\overline{\mu}}(k)z^k,\quad \sum_{k=0}^{\infty}|w_{\overline{\mu}}(k)|^2<\infty
\]
such that \cite{Hannan}
\[|w_{\overline{\mu}}(e^{-i\lambda})|^2=\frac{| \chi_{\overline{\mu}}^{(d)}(e^{-i\lambda})|^2}{|\beta^{(d)}(i\lambda)|^2},\]
which can be calculated by the formula
\begin{equation}\label{function_w}
w_{\overline{\mu}}(z)=\exp\left \{\frac{1}{4\pi}\ip\frac{e^{i\lambda}+z}{e^{i\lambda}-z}\ln
\frac{ | \chi_{\overline{\mu}}^{(d)}(e^{-i\lambda}) |^2}{|\beta^{(d)}(i\lambda)|^2}d\lambda\right \}.
\end{equation}
For this reason the following relation holds true:
\begin{equation}
 \Phi_{\overline{\mu}}(e^{-i\lambda}) =w_{\overline{\mu}}(e^{-i\lambda})\Phi (e^{-i\lambda}).
\label{dddd}\end{equation}
which implies
\[
 \varphi_{\overline{\mu}}(k)= \sum_{m=0}^kw_{\overline{\mu}}(k-m)\varphi (m),\quad k=0,1,\dots
 \]
 that is
 \[
 \varphi_{\mu,ij}(k)= \sum_{m=0}^kw_{\overline{\mu}}(k-m)\varphi_{ij} (m), \quad i=1,\dots,T;\,j=1,\dots,q;\,k=0,1,\dots.
 \]
This relation can be represented in the form
 \begin{equation}\boldsymbol{\varphi}_{\overline{\mu}} =\me W^{\overline{\mu}}\boldsymbol{\varphi}, \label{dd1}
 \end{equation}
where $\boldsymbol{\varphi}_{\overline{\mu}}=(\varphi_{\overline{\mu}}(0),\varphi_{\overline{\mu}}(1),\varphi_{\overline{\mu}}(2),\ldots)^{\top}$
and $\boldsymbol{\varphi} =(\varphi (0),\varphi (1),\varphi (2),\ldots)^{\top}$ are vectors composed from matrices
$\varphi_{\overline{\mu}}(k)=\{\varphi_{\mu,ij}(k)\}_{i=\overline{1,T}}^{j=\overline{1,q}},\, k=0,1,2,\dots$, and
$\varphi(k)=\{\varphi_{ij}(k)\}_{i=\overline{1,T}}^{j=\overline{1,q}},\, k=0,1,2,\dots$, and where
 $\me W^{\overline{\mu}}$ is a linear operator  with the entries
 $(\me W^{\overline{\mu}})_{j,k}=w_{\overline{\mu}}(j-k)$ if $0\leq k\leq j $, and $(\me W^{\overline{\mu}})_{j,k}=0$ if $0\leq j<k$.
\end{zau}

\subsection{Forecasting of periodically stationary GM increment}

Consider the problem of mean square optimal linear estimation of the functionals
\begin{equation}
A{\zeta}=\sum_{k=0}^{\infty}{a}^{(\zeta)}(k)\zeta(k), \quad
A_{M}{\zeta}=\sum_{k=0}^{N}{a}^{(\zeta)}(k)\zeta(k)
\end{equation}
which depend on unobserved values of the stochastic sequence ${\zeta}(m)$ with periodically stationary
increments. Estimates are based on observations of the sequence $\zeta(m)$ at points $m=-1,-2,\ldots$.

The functional $A{\zeta}$ can be represented in the form
\begin{eqnarray}
\nonumber
A{\zeta}& = &\sum_{k=0}^{\infty}{a}^{(\zeta)}(k)\zeta(k)=\sum_{m=0}^{\infty}\sum_{p=1}^{T}
{a}^{(\zeta)}(mT+p-1)\zeta(mT+p-1)
\\\nonumber
& = & \sum_{m=0}^{\infty}\sum_{p=1}^{T}a_p(m)\xi_p(m)=\sum_{m=0}^{\infty}(\vec{a}(m))^{\top}\vec{\xi}(m)=A\vec{\xi},
\end{eqnarray}
where
\be \label{zeta}
\vec{\xi}(m)=({\xi}_1(m),{\xi}_2(m),\dots,{\xi}_T(m))^{\top},\,
 {\xi}_p(m)=\zeta(mT+p-1);\,p=1,2,\dots,T;
\ee
\be \label{azeta}
 \vec{a}(m) =({a}_1(m),{a}_2(m),\dots,{a}_T(m))^{\top},\,
 {a}_p(m)=a^{(\zeta)}(mT+p-1);\,p=1,2,\dots,T.
\ee

Making use of the introduced notations and statements of Theorem \ref{thm_est_A} we can claim that the following theorem holds true.

\begin{thm}
\label{thm_est_Azeta}
Let a stochastic sequence ${\zeta}(k)$ with periodically stationary increments generate by formula \eqref{zeta}
 a vector-valued stochastic sequence $\vec{\xi}(m) $ which determine a
stationary stochastic GM increment sequence
$\chi_{\overline{\mu},\overline{s}}^{(n)}(\vec{\xi}(m))$ with the spectral density matrix $f(\lambda)=\{f_{ij}(\lambda)\}_{i,j=1}^{T}$ that satisfy the minimality condition (\ref{umova11_e_d}).
Let the coefficients $\vec {a}(k), k\geqslant 0$ determined by formula \eqref{azeta}  satisfy conditions  (\ref{umovana a_e_d}) -- (\ref{umovana a_mu_e_d}).
\\
Then the optimal linear estimate $\widehat{A}\zeta$ of the functional $A\zeta$ based on observations of the sequence
$\zeta(m)$ at points $m=-1,-2,\ldots$ is calculated by formula (\ref{otsinka A_e_d}).
The spectral characteristic
$\vec h_{\overline{\mu}}(\lambda)=\{h_{p}(\lambda)\}_{p=1}^{T}$ and the value of the mean square error $\Delta(f;\widehat{A}\zeta)$ of the   estimate $\widehat{A}\zeta$ are calculated by formulas
(\ref{spectr A_e_d}) and (\ref{pohybka}) respectively.
\\
In the case where the spectral density matrix $f(\lambda)$ admits the canonical factorization (\ref{fact_f_e_d}), the spectral characteristic and the value of the mean square error of the   estimate $\widehat{A}\xi$ can be calculated by formulas (\ref{simple_spectr A_e_d}) and (\ref{nsimple_poh A_e_d}) respectively.
\end{thm}

The functional $A_M{\zeta}$ can be represented in the form
\begin{eqnarray}
\nonumber
A_M{\zeta}& = &\sum_{k=0}^{M}{a}^{(\zeta)}(k)\zeta(k)=\sum_{m=0}^{N}\sum_{p=1}^{T}
{a}^{(\zeta)}(mT+p-1)\zeta(mT+p-1)
\\\nonumber
& = &\sum_{m=0}^{N}\sum_{p=1}^{T}a_p(m)\xi_p(m)=\sum_{m=0}^{N}(\vec{a}(m))^{\top}\vec{\xi}(m)=A_N\vec{\xi},
\end{eqnarray}
where $N=[\frac{M}{T}]$, the sequence $\vec{\xi}(m) $ is determined by formula \eqref{zeta},
\begin{eqnarray}
\nonumber
(\vec{a}(m))^{\top}& = &({a}_1(m),{a}_2(m),\dots,{a}_T(m))^{\top},
\\\nonumber
 {a}_p(m)& = &a^{\zeta}(mT+p-1);\,0\leq m\leq N; 1\leq p\leq T;\,mT+p-1\leq M;
\\  {a}_p(N)& = &0;\quad
M+1\leq NT+p-1\leq (N+1)T-1;1\leq p\leq T. \label{aNzeta}
\end{eqnarray}

Making use of the introduced notations and statements of Theorem \ref{thm_est_A_N} we can claim that the following theorem holds true.

\begin{thm}
\label{thm_est_A_Nzeta}
Let a stochastic sequence ${\zeta}(k)$ with periodically stationary GM increments generate by formula \eqref{zeta}
 a vector-valued stochastic sequence $\vec{\xi}(m) $ which determine a
stationary  GM increment sequence
$\chi_{\overline{\mu},\overline{s}}^{(n)}(\vec{\xi}(m))$ with the spectral density matrix $f(\lambda)=\{f_{ij}(\lambda)\}_{i,j=1}^{T}$ that satisfy the minimality condition (\ref{umova11_e_d}).
Let coefficients $\vec {a}(k), k\geqslant 0$ be determined by formula \eqref{aNzeta}.
The optimal linear estimate $\widehat{A}_M\zeta$ of the functional $A_M\zeta=A_N\vec{\xi}$ based on observations of the sequence
$\zeta(m)$ at points $m=-1,-2,\ldots$ is calculated by formula \eqref{estimation_A_N}.
The spectral characteristic $\vec{h}_{\overline{\mu},N}(\lambda)=\{h_{\overline{\mu},N,p}(\lambda)\}_{p=1}^{T}$ and the value of the mean square error $\Delta(f;\widehat{A}_M\zeta)$
are calculated by formulas  \eqref{est_h_N} and \eqref{pohybka_N} respectively.
In the case where the spectral density matrix $f(\lambda)$ admits the canonical factorization (\ref{fact_f_e_d}), then the spectral characteristic
$\vec h_{\overline{\mu},N}(\lambda)$ and the value of the mean square error of the   estimate $\widehat{A}_M\zeta$ can be calculated by formulas
(\ref{simple_spectr A_N_e_d}) and (\ref{nsimple_poh A_N_e_d}) respectively.
\end{thm}

As a corollary from Theorem \ref{thm_est_A_Nzeta}, one can obtain the mean square optimal estimate of the unobserved value
$\zeta(M)$, $M\geq0$ of a stochastic sequence ${\zeta}(m)$ with periodically stationary GM increments
based on observations of the sequence ${\zeta}(m)$ at points $m=-1,-2,\ldots$
Making use of the notations
$\zeta(M)=\xi_p(N)=(\vec\xi(N))^{\top}\boldsymbol{\delta}_p$, $N=[\frac{M}{T}]$, $p=M+1-NT$,
and the obtained results we can conclude that the following corollary holds true.

\begin{nas}\label{nas zeta_e_d}
Let a stochastic sequence ${\zeta}(m)$ with periodically stationary GM increments generate by formula \eqref{zeta}
 a vector-valued stochastic sequence $\vec{\xi}(m) $ which determine a
stationary  GM increment sequence
$\chi_{\overline{\mu},\overline{s}}^{(n)}(\vec{\xi}(m))$ with the spectral density matrix $f(\lambda)=\{f_{ij}(\lambda)\}_{i,j=1}^{T}$ that satisfy the minimality condition (\ref{umova11_e_d}).
The optimal linear estimate $\widehat{\zeta}(M)$ of the unobserved value
$\zeta(M)$, $M\geq0$,
based on observations of the sequence ${\zeta}(m)$ at points $m=-1,-2,\ldots$
is calculated by formula \eqref{est_xi_N}.
The spectral characteristic $\vec h_{\overline{\mu},N,p}(\lambda)$ of the estimate is calculated by the formula
 \eqref{sph_est_xi_N}.
The value of the mean square error of the   estimate $\widehat{\zeta}(M)$ is calculated by the formula
 \eqref{delta_est_xi_N}.
If the spectral density $f(\lambda)$ admits the canonical factorization (\ref{fact_f_e_d}), then the spectral characteristic
and the value of the mean square error of the   estimate $\widehat{\zeta}(M)$ can be calculated by the formulas
 \eqref{sph_f_est_xi_N}, \eqref{delta_f_est_xi_N}.
\end{nas}

\section{Forecasting of GM fractional increments}\label{fractional_extrapolation}

In the previous section, we solved  the forecasting problem for the increment sequence $\chi_{\overline{\mu},\overline{s}}^{(d)}(\vec{\xi}(m))$ of  the positive integer orders $(d_1,\ldots,d_r)$. Here we consider the forecasting problem in the case of fractional increment orders $d_i$.

Within the section, we consider the step $\overline{\mu}=(1,1,\ldots,1)$ and represent the increment operator $\chi_{\overline{s}}^{(d)}(B)$  in the form
\be\label{FM_increment}
\chi_{\overline{s}}^{(R+D)}(B)=(1-B)^{R_0+D_0}\prod_{j=1}^r(1-B^{s_j})^{R_j+D_j},
\ee
where $(1-B)^{R_0+D_0}$ is the integrating component, $R_j$, $j=0,1,\ldots, r$, are non-negative integer numbers, $1<s_1<\ldots<s_r$.  
The goal  is to find representations $d_j=R_j+D_j$, $j=0,1,\ldots, r$, of the increment orders 
under some  conditions on the fractional parts $D_j$, such that the increment sequence $\vec y(m):=(1-B)^{R_0}\prod_{j=1}^r(1-B^{s_j})^{R_i}\vec{\xi}(m)$ to be  a stationary  fractionally integrated seasonal stochastic  sequence.
For example, in case of single  increment pattern $(1-B^{s^*})^{R^*+D^*}$ this condition is $|D^*|<1/2$.

We will call a  sequence $\chi_{\overline{s}}^{(R+D)}(\vec \xi(m))$ \emph{a fractional multiple (FM) increment sequence}.

\begin{lema}\label{frac_incr_1}
The increment operator $\chi_{\overline{s}}^{(D)}(B):=(1-B)^{D_0}\prod_{j=1}^r(1-B^{s_j})^{D_j}$ admits a representation
\begin{eqnarray*}
\chi_{\overline{s}}^{(D)}(B)
& = &\prod_{j=0}^r\prod_{k_j=0}^{[s_j/2]}(1-2\cos \nu_{k_j} B+B^2)^{D_{k_j}}
\\& = &(1-B)^{D_0+D_1+\ldots+D_r}\prod_{j=1}^r\prod_{k_j=1}^{[s_j/2]}(1-2\cos \nu_{k_j} B+B^2)^{D_{k_j}},
\end{eqnarray*}
where $s_0=1$, $\nu_{k_j}=2\pi k_j/s_j$, $k_j=0,1,\ldots, [s_j/2]$, $D_{k_j}=D_0/2$ for $k_j=0$, $D_{k_j}=D_j$ for $k_j=1,2, \ldots,[s_j/2]-1$, $D_{[s_j/2]}=D_j$ for odd $s_j$ and $D_{[s_j/2]}=D_j/2$ for even $s_j$,
\end{lema}

Note that Lemma \ref{frac_incr_1} follows from the representation
\[
(1-B^{s_j})^{D_r}=\prod_{k_j=0}^{[s_j/2]}(z_{k_j}-B)^{D_{k_j}}(z_{-k_j}-B)^{D_{k_j}},
\]
where $z_{k_j}=\exp(\nu_{k_j}i)$, $k_j=0,1,\ldots,s_j-1$, are solutions of the equation $1-B^{s_j}=0$.

Lemma \ref{frac_incr_1} implies the following statement.

\begin{lema}\label{frac_incr_2}
Define the sets $\md M_j=\{\nu_{k_j}=2\pi k_j/s_j: k_j=0,1,\ldots, [s_j/2]\}$, $j=0,1,\ldots, r$, and the set $\md M=\bigcup_{j=0}^r \md M_j$. Then
\begin{eqnarray*}
\chi_{\overline{s}}^{(D)}(B)
& = &\prod_{\nu \in \md M}(1-2\cos \nu B+B^2)^{\widetilde{D}_{\nu}}
\\
& = &(1-B)^{D_0+D_1+\ldots+D_r}(1+B)^{D_{\pi}}\prod_{\nu \in \md M\setminus\{0,\pi\}}(1-2\cos \nu B+B^2)^{D_{\nu}},
\end{eqnarray*}
where $D_{\nu}=\sum_{j=0}^rD_j \mr I \{\nu\in \md M_j\}$, $\widetilde{D}_{\nu}=D_{\nu}$ for $\nu \in \md M\setminus\{0,\pi\}$, $\widetilde{D}_{\nu}=D_{\nu}/2$ for $\nu=0$ and $\nu=\pi$.
\end{lema}


Lemma \ref{frac_incr_2} shows that a multiple seasonal increment sequence can be represented as the following $k$-factor Gegenbauer sequence 
\be \label{Gegenbauer}
\prod_{i=1}^k(1-2 u_i B+B^2)^{d_i} x(m)=\xi(m).
\ee
In the case where $\xi(m)$ is an $ARMA(p,q)$ sequence, the model $x(m)$ defined by  (\ref{Gegenbauer})   is called $k$-factor $GARMA(p,d_i,u_i,q)$ sequence.
It is stationary and invertible if $|d_i|<1/2$ for  $|u_i|<1$ and $|d_i|<1/4$ for  $|u_i|=1$.
If additionally $d_i>0$, then the model exhibits a long memory behavior \cite{Woodward}.
The function $(1-2 u_i B+B^2)^{-d_i}$ is a generating function of the Gegenbauer polynomial:
\[
(1-2 u B+B^2)^{-d}=\sum_{n=0}^{\infty}C_n^{(d)}(u)B^n,
\]
where
\[
C_n^{(d)}(u)=\sum_{k=0}^{[n/2]}\frac{(-1)^k(2u)^{n-2k}\Gamma(d-k+n)}{k!(n-2k)!\Gamma(d)}.
\]
Thus, denoting $k^*=|\md M|$, we obtain
\begin{eqnarray*}
(\chi_{\overline{s}}^{(D)}(B))^{-1}& = &\prod_{\nu \in \md M}(1-2\cos \nu B+B^2)^{-\widetilde{D}_{\nu}}
\\
& = &\sum_{m=0}^{\infty}G^+_{k^*}(m)B^m=\ld(\sum_{m=0}^{\infty}G^-_{k^*}(m)B^m\rd)^{-1},
\end{eqnarray*}
where
\begin{eqnarray}
\label{Gegenbauer_GI+}
G^+_{k^*}(m)& = &\sum_{0\leq n_1,\ldots,n_{k^*}\leq m, n_1+\ldots+n_{k^*}=m}\prod_{\nu \in \md M}C_{n_{\nu}}^{(\widetilde{D}_{\nu})}(\cos\nu),
\\
\label{Gegenbauer_GI-}
G^-_{k^*}(m)& = &\sum_{0\leq n_1,\ldots,n_{k^*}\leq m, n_1+\ldots+n_{k^*}=m}\prod_{\nu \in \md M}C_{n_{\nu}}^{(-\widetilde{D}_{\nu})}(\cos\nu).
\end{eqnarray}

The derived representations of the increment operator $\chi_{\overline{s}}^{(D)}(B)$ imply the following theorem.

\begin{thm}\label{thm_frac}
Assume that for a stochastic vector sequence $\vec \xi(m)$ and fractional differencing orders $d_j=R_j+D_j$, $j=0,1,\ldots, r$, the FM increment sequence $\chi_{\overline{1},\overline{s}}^{(R+D)}(\vec \xi(m))$ generated by increment operator (\ref{FM_increment})  is a stationary sequence with a bounded from zero and infinity spectral density $\widetilde{f}_{\overline{1}}(\lambda)$. Then for the non-negative integer numbers $R_j$, $j=0,1,\ldots, r$, the GM increment sequence $\chi_{\overline{1},\overline{s}}^{(R)}(\vec \xi(m))$    is stationary if $-1/2< D_{\nu}<1/2$ for all $\nu\in \md M$, where $D_{\nu}$ are defined by real numbers $D_j$, $j=0,1,\ldots, r$,
in Lemma \ref{frac_incr_2}, and it is long memory if $0< D_{\nu}<1/2$ for at least one $\nu\in \md M$, and invertible if $-1/2< D_{\nu}<0$. The spectral density $f(\lambda)$ of the stationary GM increment sequence $\chi_{\overline{1},\overline{s}}^{(R)}(\vec \xi(m))$ admits a representation
\[
 f(\lambda)=|\beta^{(R)}(i\lambda)|^2 \ld|\chi_{\overline{1}}^{(R)}(e^{-i\lambda})\rd|^{-2}\ld|\chi_{\overline{1}}^{(D)}(e^{-i\lambda})\rd|^{-2} \widetilde{f}_{\overline{1}}(\lambda)=:\ld|\chi_{\overline{1}}^{(D)}(e^{-i\lambda})\rd|^{-2} \widetilde{f} (\lambda),
  \]
  where
  \begin{eqnarray*}
  \ld|\chi_{\overline{1}}^{(D)}(e^{-i\lambda})\rd|^{-2}& = &\ld|\sum_{m=0}^{\infty}G^+_{k^*}(m) e^{-i\lambda m}\rd|^2=\ld|\sum_{m=0}^{\infty}G^-_{k^*}(m) e^{-i\lambda m}\rd|^{-2}
\\
& = &\prod_{\nu \in \md M}\ld|(e^{-i\nu}-e^{i\lambda})(e^{i\nu}-e^{i\lambda})\rd|^{-2\widetilde{D}_{\nu}},
  \end{eqnarray*}
 coefficients $G^+_{k^*}(m)$, $G^-_{k^*}(m)$ are defined by (\ref{Gegenbauer_GI+}), (\ref{Gegenbauer_GI-}).
\end{thm}

The spectral density $f(\lambda)$ and the structural function $D^{(R)}_{ \overline{s}}(m,\overline{1},\overline{1})$ of a stationary GM increment sequence $\chi_{\overline{1},\overline{s}}^{(R)}(\vec \xi(m))$ exhibit the following behavior in the case of constant matrices $C$ and $K$:
\begin{itemize}
\item   $|\beta^{(R)}(i\lambda)|^{-2} |\chi_{\overline{1}}^{(R)}(e^{-i\lambda})|^2f(\lambda)\sim C|\nu-\lambda|^{-2\widetilde{D}_{\nu}}$
as $\lambda\to \nu$, $\nu\in \mathcal{M}$, thus, the minimality condition (\ref{umova11_e_d}) is satisfied (for properties of eigenvalues of generalized fractional process, we refer to Palma and Bondon \cite{Palma-Bondon})

\item $D^{(R)}_{ \overline{s}}(m,\overline{1},\overline{1})\sim K\sum_{\nu\in \mathcal{M}:\widetilde{D}_{\nu}>0}|m|^{2\widetilde{D}_{\nu}-1}\cos (m\nu)$,   as $m\to\infty$ (see Giraitis and Leipus \cite{Giraitis}).
\end{itemize}

\begin{pry}
1. Consider an increment operator $(1-B)^{R_0+D_0}(1-B^2)^{R_1+D_1}$ which represents a fractional integrated component and a fractional seasonal components.
In this case $\md M_0=\{0\}$, $\md M_1=\{0,\pi\}$, $\md M=\{0,\pi\}$. The Gegenbauer representation of the increment is $(1-B)^{D_0+D_1}(1+B)^{D_1}$.
Stationarity conditions are the following: $|D|=|D_0+D_1|<1/2$, $|D_{\pi}|=|D_1|<1/2$.

2. Consider an increment operator $(1-B^2)^{R_1+D_1}(1-B^3)^{R_2+D_2}$ which represents two fractional  seasonal components.
In this case  $\md M_0=\{0,\pi\}$, $\md M_1=\{0,2\pi/3\}$, $\md M=\{0,2\pi/3,\pi\}$. The Gegenbauer representation of the increment is $(1-B)^{D_1+D_2}(1-2\cos(2\pi/3)B+B^2)^{D_2}(1+B)^{D_1}$.
Stationarity conditions are the following: $|D|=|D_1+D_2|<1/2$, $|D_{2\pi/3}|=|D_2|<1/2$, $|D_{\pi}|=|D_1|<1/2$.

3. Consider an increment operator $(1-B^2)^{R_1+D_1}(1-B^4)^{R_2+D_2}$.
In this case $\md M_0=\{0,\pi\}$, $\md M_1=\{0,\pi/2,\pi\}$, $\md M=\{0,\pi/2,\pi\}$.
The Gegenbauer representation of the increment is $(1-B)^{D_1+D_2}(1+B^2)^{D_2}(1+B)^{D_1+D_2}$.
Stationarity conditions are the following: $|D|=|D_{\pi}|=|D_1+D_2|<1/2$, $|D_{\pi/2}|=|D_2|<1/2$.
\end{pry}

In the following remarks we provide some additional details with the help of which we can use theorems proposed in the previous section in finding solution of the forecasting problem for stochastic
sequences with periodically stationary (periodically correlated) FM increments.

\begin{zau}
Theorem \ref{thm_frac} implies that the Fourier coefficients \eqref{Fourier_F} of the function
\[
|\beta^{(R)}(i\lambda)|^2| \chi_{\overline{1}}^{(R)}(e^{-i\lambda})|^{-2}f^{-1}(\lambda)\]
 are calculated by the formula
\[
 F^{\overline{1}}(k,j)=\frac{1}{2\pi}\ip
 e^{i\lambda(j-k)}\ld|\chi_{\overline{1}}^{(D)}(e^{-i\lambda})\rd|^2\widetilde{f}^{-1}_{\overline{1}}(\lambda)d\lambda,\quad k,j\geq0.
\]
\end{zau}

\begin{zau}
Assume that the spectral density $\widetilde{f}_{\overline{1}}(\lambda)$ admits a factorization
\[
\widetilde{f}_{\overline{1}}(\lambda)=\ld|\widetilde{\Phi}_{\overline{1}}(e^{-i\lambda})\rd|^2
=\ld|\sum_{k=0}^{\infty}\widetilde{\varphi}_{\overline{1}}(k)(e^{-i\lambda})\rd|^2,
\]
where
$\widetilde{\varphi}_{\overline{1}}(k)=\{\varphi_{\overline{1},ij}(k)\}_{i=\overline{1,T}}^{j=\overline{1,q}}$, $ k=0,1,2,\dots$.
Then coefficients $\{\varphi_{\overline{1},ij}\}_{i=\overline{1,T}}^{j=\overline{1,q}}$, $ k=0,1,2,\dots$ from factorization (\ref{fact_f_e_d}) are calculated by the formula
 \[
 \varphi_{\overline{1},ij}(k)= \sum_{m=0}^kG^+_{k^*}(k-m)\widetilde{\varphi}_{\overline{1},ij} (m)=(G^+_{k^*}*\widetilde{\varphi}_{\overline{1},ij})(k).
 \]
 \end{zau}

 \begin{zau}
Define a matrix-valued function $\widetilde{\Psi}_{\overline{1}}(e^{-i\lambda})= \{\widetilde{\Psi}_{\overline{1},ij}(e^{-i\lambda})\}_{i=\overline{1,q}}^{j=\overline{1,T}}$ by the equation
$\widetilde{\Psi}_{\overline{1}}(e^{-i\lambda})\widetilde{\Phi}_{\overline{1}}(e^{-i\lambda})=E_q$,
where $E_q$ is the identity $q\times q$ matrix. Then $\Psi_{\overline{1}}(e^{-i\lambda})=\chi_{\overline{1}}^{(D)}(e^{-i\lambda})\widetilde{\Psi}_{\overline{1}}(e^{-i\lambda})$.
\end{zau}

\section{Examples of forecasting for some special models}\label{examples}

\begin{pry}
Basawa et al. \cite{Basawa} consider a so-called first-order seasonal periodic autoregressive process (SPAR(1,1)) defined by the difference equation
\be\label{SPAR_model}
    X_{nT+\nu}=\phi(\nu)X_{mT+\nu-1}+\alpha(\nu)X_{(m-1)T+\nu}-\phi(\nu)\alpha(\nu)X_{(m-1)T+\nu-1}+\eps_{mT+\nu},
\ee
where $\eps_{mT+\nu}$ is an uncorrelated periodic white noise process with $\mt E (\eps_{mT+\nu})=0$ and $\mt{Var}(\eps_{mT+\nu})=\sigma^2(\nu)$, $1\leq \nu\leq T$ (we follow  notations from \cite{Basawa}). Depending on coefficients $\phi(\nu)$, $\alpha(\nu)$  model \eqref{SPAR_model} has the following properties:
\begin{itemize}
\item if $\phi(\nu)\equiv\phi$, $\alpha(\nu)\equiv\alpha$, $\sigma^2(\nu)\equiv\sigma^2$ for $1\leq \nu\leq T$, then model (\ref{SPAR_model}) reduces to Box-Jenkins SAR(1,1) model,

\item if  $\alpha(\nu)\equiv 0$  for $1\leq \nu\leq T$, then model (\ref{SPAR_model}) reduces to  PAR(1) model,

\item if $\prod_{\nu=1}^T|\phi(\nu)|<1$ and $|\alpha(\nu)|<1$ for all $\nu$, then model (\ref{SPAR_model}) admits a causal and stationary $T$-dimensional VAR representation
\[
    \me\Phi_0 \me X_m=\me \Phi_1 \me X_{m-1}+\me \Phi_2 \me X_{m-2}+ \boldsymbol{\eps}_n,
\]
where
\begin{eqnarray}
\nonumber
\me X_m & = & (X_{mT+1},X_{mT+2},\ldots,X_{mT+T})^{\top},
\\\nonumber
 \boldsymbol{\eps}_m & = & (\eps_{mT+1},\eps_{mT+2},\ldots,\eps_{mT+T})^{\top}.
\end{eqnarray}
\end{itemize}

Here we consider another case where $\prod_{\nu=1}^T|\phi(\nu)|<1$ and $\alpha(\nu)\equiv 1$ for all $\nu$.
Without loss of generality assume that $\sigma^2(\nu)\equiv 1$ for all $\nu$.
Then, taking into account $B^TX_{mT+\nu}= X_{(m-1)T+\nu}$, model (\ref{SPAR_model}) reduces to integrated PAR, or PARIMA, model
\[
    (1-B^T)(X_{mT+\nu}-\phi(\nu)X_{mT+\nu-1})=\eps_{mT+\nu},
\]
which admits a VARIMA(1,1,0) representation
\[
\me \Psi_0 \Delta\me X_m+\me \Psi_1 \Delta\me X_{m-1}=\boldsymbol{\eps}_m,
\]
where $\me \Psi_0(r,s)=1$ for $r=s$, $\me \Psi_0(r,s)=-\phi(s)$ for $r=s+1$ and $\me \Psi_0(r,s)=0$ otherwise; $\me \Psi_1(1,T)=-\phi(T)$ and $\Psi_1(r,s)=0$ otherwise, $1\leq r,s\leq T$.  Note that $\Delta\me X_m=(1-B)\me X_m=\chi_{1,1}^{(1)}(\me X_m)$ in terms of GM increments. The spectral density of the one-step increment sequence $\Delta\me X_m$ is the following:
\[
f(\lambda)=\frac{\lambda^2}{|1-e^{-i\lambda}|^2}\ld|\Psi_0+\Psi_1e^{-i\lambda}\rd|^{-2}
.
\]

Consider the following estimation problem. Let us assume that we observe a time series $X_{mT+\nu}$ at points $m\leq-1$, $1\leq\nu\leq T$.
It is necessary to find an estimate $\widehat{A}X$ of the functional which depends on future values of $X_{mT+\nu}$, $m\geq0$, $1\leq\nu\leq T$, with a discount factor $\rho, 0<\rho<1$:
\[
AX=\sum_{m=0}^{\infty}\sum_{\nu=1}^T\rho^{mT+\nu}X_{mT+\nu}=\sum_{m=0}^{\infty}\me a_m^{\top}\me X_m=:
A\me X,
\]
where $\me a_m=\rho^{mT}(\rho,\rho^2,\ldots,\rho^T)^{\top}=\rho^{mT} \me a$, $\me a=(\rho,\rho^2,\ldots,\rho^T)^{\top}$. Coefficients $\me b_m$, $m\geq0$ and $\me v_{-1}$ from the representation $A\me X=B\Delta\me X-V\me X$ are the following: $\me b_m=\frac{\rho^{mT}}{1-\rho^T}\me a$, $\me v_{-1}=\me b_0=-\frac{1}{1-\rho^T}\me a$. Since $\|\me a_m\|^2=c_1(T)\rho^{mT}$, $\|\me b_m\|^2=c_2(T)\rho^{mT}$, conditions (\ref{umovana a_e_d}) and (\ref{umovana a_mu_e_d}) are satisfied. Thus, we apply Theorem \ref{thm_est_Azeta} to find the spectral characteristic of the estimate $\widehat{A} X$.
We have
\[
\Phi_{(1)}(e^{-i\lambda})=\Psi_0^{-1}\sum_{k=0}^{\infty}(-1)^k(\Psi_1\Psi_0^{-1})^ke^{-i\lambda k},
\quad
\Psi_{(1)}(e^{-i\lambda})=\Psi_0+\Psi_1e^{-i\lambda},
\]
and
\[
\vec r_{(1)}(e^{i\lambda})=\frac{1}{1-\rho^T}\Theta^{\top}\me a\sum_{k=0}^{\infty}\rho^{kT}e^{i\lambda k},
\]
where $\Theta:=(\Psi_0 +\rho^T\Psi_1 )^{-1}$,
\[
(\Psi_{(1)}(e^{-i\lambda}))^{\top}
 \vec r_{(1)}(e^{i\lambda})=\frac{1}{1-\rho^T}\Psi_1^{\top}\Theta^{\top}\me a e^{-i\lambda }+\frac{1}{1-\rho^T}\me a\sum_{k=0}^{\infty}\rho^{kT}e^{i\lambda k},
\]
Then, the optimal estimate of the value of the functional $AX$ is calculated by the formula
\begin{eqnarray*}
\widehat{A}X & = & -\frac{1}{1-\rho^T}\me a^{\top} \Theta\Psi_1\Delta\me X_{-1}+\frac{1}{1-\rho^T}\me a^{\top}\me X_{-1}
\\
 & = & \frac{1}{1-\rho^T}\me a^{\top}\ld((E_T-\Theta\Psi_1)\me X_{-1}+\Theta\Psi_1\me X_{-2}\rd).
\end{eqnarray*}
The value of the mean square error of the estimate is calculated by the formula
\[
\Delta\ld(f;\widehat{A}X\rd)= \frac{1}{(1-\rho^T)^3(1+\rho^T)}\|\Theta^{\top}\me a\|^2.
\]
\end{pry}

\begin{pry}\label{example_MA}
To illustrates a  forecasting technique developed in Chapters \ref{classical_extrapolation} and \ref{fractional_extrapolation} we consider a seasonal time series $x(t)$, $t\in \mr Z$, exhibiting two fractional seasonal patterns and a periodic covariance behavior
\[
(1-B^s)^{d_0}(1-B^{us})^{d_1}\xi(t)=\eps(t)-a_0\eps(t-1)-a_{i(t)}\eps(t-s),\quad i(t)=(t\mod s)+1,
\]
where $d_0=1+D_0$, $d_1=1+D_1$, $\eps(t)$, $t\in \mr Z$, are i.i.d. random variables with $\mt E\eps(t)=0$, $\mt E|\eps(t)|^2=1$.
The first cycle $s$ may refer to $7$ days within a week, and this pattern shows different correlation structure for each `season'', namely, day of a week.
The second seasonal pattern $us$ may describe a year period assuming that  $u=52$ corresponds to weeks within a year.
 Under the conditions stated below, the increment $w(t)=(1-B^s)(1-B^{us})\xi(t)$ is cyclostationary since coefficients $a_{i(t)}$ are periodic with the period $T=s$.

Define the vector-valued sequences $\vec \xi(m)=(\xi_1(m), \xi_2(m),\ldots,\xi_s(m))^{\top}$, where $\xi_p(m)=\xi(sm+p-1)$, and  $\vec \eps(m)=(\eps_1(m), \eps_2(m),\ldots,\eps_s(m))^{\top}$, where $\eps_p(m)=\eps(sm+p-1)$. Consider  an increment function $\chi_{(1,1),(1,u)}^{(2)}(B)=(1-B)(1-B^u)$ with the step $\overline{\mu}=(1,1)$. The GM increment sequence $\chi_{(1,1),(1,u)}^{(2)}(\vec \xi(m))$ admits the representation
\[
(1-B)^{D_0}(1-B^u)^{D_1}\chi_{(1,1),(1,u)}^{(2)}(\vec \xi(m))=\Phi_0
\vec \eps(m)
+
\Phi_1\vec \eps(m-1),
\]

\[
\Phi_0=\left(
\begin{array}{ccccc}
1 & 0 & 0& \ldots & 0\\
-a_0 & 1 & 0& \ldots & 0\\
0 & -a_0 & 1& \ddots & 0\\
\vdots & \vdots & \ddots & \ddots & \vdots\\
0 & 0 & 0 & \ldots & 1
\end{array}
\right),
\,
\Phi_1=\left(
\begin{array}{ccccc}
-a_1 &  0 & 0& \ldots & -a_0\\
0 & -a_2 & 0 & \ldots & 0\\
0 & 0 & -a_3& \ddots & 0\\
\vdots & \vdots & \ddots & \ddots & \vdots\\
0 & 0 & 0 & \ldots & -a_s
\end{array}
\right).
\]
It is stationary under conditions  $|D_0+D_1|<1/2$, $|D_1|<1/2$. For instance, if $d_0=0.7$, $d_1=1.2$, then $D_0=-0.3$, $D_1=0.2$ and the process exhibits a long-memory behavior.

The spectral density $f(\lambda)$ of the GM increment sequence $\chi_{(1,1),(1,d)}^{(2)}\vec \xi(m)$ is
\begin{eqnarray*}
\nonumber
f(\lambda)
 & = & \frac{|\beta^{(2)}(i\lambda)|^2}{|\chi_{(1,1)}^{(2)}(e^{-i\lambda})|^2|\chi_{(1,1)}^{(D)}(e^{-i\lambda})|^2}\ld|\Phi_{(1,1)}(e^{-i\lambda})\rd|^2
\\\nonumber
 & = & \frac{\lambda^2\prod_{k=-[u/2]}^{[u/2]}(\lambda-2\pi k/u)^{2}}{|1-e^{-i\lambda}|^2|1-e^{-iu\lambda}|^2|1-e^{-i\lambda}|^{2D_0}|1-e^{-iu\lambda}|^{2D_1}}\ld|\Phi_0+\Phi_1e^{-i\lambda}\rd|^2,
\end{eqnarray*}
where
\begin{eqnarray*}(\chi_{(1,1)}^{(D)}(e^{-i\lambda}))^{-1} & = & (1-e^{-i\lambda})^{-D_0 }(1-e^{-iu\lambda})^{ -D_1}
\\ & = &
(1-e^{-i\lambda})^{-D_0-D_1}\prod_{k=1}^{[u/2]}(1-2\cos (2\pi k/u) e^{-i\lambda}+e^{-2i\lambda})^{-\widetilde{D}_k}
\\ & = & \sum_{k=0}^{\infty}G^+_{k^*}(k)e^{-i\lambda k}=\ld(\sum_{k=0}^{\infty}G^-_{k^*}(k)e^{-i\lambda k}\rd)^{-1},
\end{eqnarray*}
$k^*=[u/2]+1$, $G^+_{k^*}(m)$ and $G^-_{k^*}(m)$, $m\geq0$, are defined by (\ref{Gegenbauer_GI+}) and (\ref{Gegenbauer_GI-}) respectively, $\widetilde{D}_k=D_1$ for $1\leq k\leq[u/2]-1$, $\widetilde{D}_{[u/2]}=D_1$ for odd $u$ and $\widetilde{D}_{[u/2]}=D_1/2$ for even $u$. Note,  $G^+_{k^*}(0)=G^-_{k^*}(0)=1$.

Let us find an estimate of  a weighted sum of two average weekly values of the time series $\xi(t)$
\[
A_{2s}\xi=\alpha\ld(\frac{1}{s}\sum_{k=0}^{s-1}\xi(k)\rd)+(1-\alpha)\ld(\frac{1}{s}\sum_{k=s}^{2s-1}\xi(k)\rd)
\]
based on observations of $\xi(t)$ at points $t=-1,-2,\dots$.

In terms of the sequence $\vec \xi(m)$,  the functional $A_{2s}\xi$ is rewritten as
\[
A_2\vec \xi=(\vec{a}(0))^{\top}\vec{\xi}(0)+(\vec{a}(1))^{\top}\vec{\xi}(1),
\]
where $\vec{a}(0)=(\alpha s^{-1}, \alpha s^{-1}, \ldots,\alpha s^{-1})^{\top}$, $\vec{a}(1)=((1-\alpha) s^{-1}, (1-\alpha) s^{-1}, \ldots,(1-\alpha) s^{-1})^{\top}$, and admits a representation
\[
A_2\vec \xi=B_2\chi \vec \xi  -V_2\vec \xi=(\vec{b}(0))^{\top}\chi(\vec{\xi}(0))+(\vec{b}(1))^{\top}\chi(\vec{\xi}(1))-\sum_{k=-u-1}^{-1} (\vec{v}_2(k))^{\top}\vec{\xi}(k).
\]
Here $\vec{b}(0)=\vec{a}(0)+\vec{a}(1)=(  s^{-1},   s^{-1}, \ldots, s^{-1})^{\top}=s^{-1}\overline{\me 1}$, $\vec{b}(1)=\vec{a}(1)=((1- \alpha)  s^{-1},  (1- \alpha)  s^{-1}, \ldots, (1- \alpha)  s^{-1})^{\top}=(1- \alpha) s^{-1}\overline{\me 1}$; and further $\vec{v}_2(k)=0$ for $k=-1,-2, \ldots,-u+2$,
$\vec{v}_2(-u+1)=-\vec{b}(1)=-(1- \alpha) s^{-1}\overline{\me 1}$, $\vec{v}_2(-u)=-\vec{b}(0)+\vec{b}(1)=-\alpha s^{-1}\overline{\me 1}$, $\vec{v}_2(-u-1)= \vec{b}(0)=s^{-1}\overline{\me 1}$. By $\overline{\me 1}$ we denote   a vector $(1,1,\ldots,1)^{\top}$ of  dimension $s$.
Note  that
\[
(1-B)^{-1}(1-B^u)^{-1}=\sum_{k=0}^{\infty}d(k)B^k=\sum_{k=0}^{\infty}\ld(1+\ld[\frac{k}{u}\rd]\rd)B^k.
\]
We find the spectral characteristic $\vec h_{(1,1),2}(\lambda)$ of the estimate $\widehat{A}_2\vec \xi$ using Theorem~\ref{thm_est_A_N} as well as remarks to Theorem~\ref{thm_frac}.
First we obtain
\begin{eqnarray*}\Phi_{(1,1)}(e^{-i\lambda}) & = & \sum_{k=0}^{\infty}(G^+_{k^*}*\Phi)(k)e^{-i\lambda k}
\\
 & = & G^+_{k^*}(0)\Phi_0
+\sum_{k=1}^{\infty}(G^+_{k^*}(k)\Phi_0+G^+_{k^*}(k-1)\Phi_1)e^{-i\lambda k},
\\
\Psi_{(1,1)}(e^{-i\lambda}) & = & \sum_{k=0}^{\infty}(G^-_{k^*}*\Psi)(k)e^{-i\lambda k}=\Phi_0^{-1}\sum_{k=0}^{\infty}(G^-_{k^*}*\widetilde{\Psi})(k)e^{-i\lambda k},
\end{eqnarray*}
where $(G^-_{k^*}*\Psi)(k)$, $k\geq0$, is a convolution of two sequences $G^-_{k^*}(k)$ and $\Psi_k$, $k\geq0$,  $\Psi_k=(-1)^k\Phi_0^{-1}(\Phi_1\Phi_0^{-1})^k$, $\widetilde{\Psi}_k=\Phi_0^{-1}\Psi_k=(-1)^k (\Phi_1\Phi_0^{-1})^k$,
\[
\Phi_1\Phi_0^{-1}=\left(
\begin{array}{ccccc}
-a_1-a_0^s &  -a_0^{s-1} &-a_0^{s-2} & \ldots &  -a_0\\
-a_2a_0  &-a_2   &  0 & \ldots &  0\\
-a_3a_0^2  &-a_3a_0 & -a_3 & \ldots &  0\\
\vdots  & \vdots & \vdots& \ddots  & \vdots\\
-a_sa_0^{s-1}  & -a_sa_0^{s-2} & -a_sa_0^{s-3}& \ldots &  - a_s
\end{array}
\right).
\]
Then
\begin{eqnarray*}
 &  & (\Psi_{(1,1)}(e^{-i\lambda}))^{\top}
 \vec r_{(1,1),2}(e^{i\lambda})=\vec b(1)e^{i\lambda}+\vec b(0)
\\
 & \quad & + \sum_{k=1}^{\infty}\ld((G^-_{k^*}*\widetilde{\Psi}^{\top})(k)( \vec b(0)+G^+_{k^*}(1) \vec b(1))+ G^-_{k^*}(k+1)\vec b(1)\rd)e^{-i\lambda k}
\\
 & = &  \chi_{_{(1,u)}}^{(D)}(e^{-i\lambda})(\widetilde{\Psi}(e^{-i\lambda}))^{\top}
\ld( \vec b(0)+D_0 \vec b(1)\rd) +  \vec b(1)e^{i\lambda}\chi_{_{(1,u)}}^{(D)}(e^{-i\lambda}),
\end{eqnarray*}
and
\begin{eqnarray*}
\vec h_{(1,1),2}(\lambda) & = & -\frac{\chi_{(1,1)}^{(2)}(e^{-i\lambda})}{\beta^{(2)}(i\lambda)} \sum_{k=1}^{\infty}s^{-1}\ld((1  + (1-\alpha) D_0)(G^-_{k^*}*\widetilde{\Psi}^{\top})(k)\rd.
\\
 & \quad & \Big.+(1-\alpha) G^-_{k^*}(k+1)E_T\Big)    \overline{\me 1} e^{-i\lambda k}
=\frac{\chi_{(1,1)}^{(2)}(e^{-i\lambda})}{\beta^{(2)}(i\lambda)} \sum_{k=1}^{\infty}h(k).
\end{eqnarray*}
The optimal estimate of the functional $A_{2s} \xi$ is calculated by the formula
\begin{eqnarray*}
\widehat{A}_2\vec \xi & = &-s^{-1}(\overline{\me 1})^{\top}\vec{\xi}(-u-1)+\alpha s^{-1}(\overline{\me 1})^{\top}\vec{\xi}(-u)+(1-\alpha) s^{-1}(\overline{\me 1})^{\top}\vec{\xi}(-u+1)
\\
 & \quad &+\sum_{k=1}^{\infty}h^{\top}(k)\ld(\vec \xi(-k)-\vec \xi(-k-1)-\vec\xi(-k-u)+\vec\xi(-k-u-1)\rd).
\end{eqnarray*}

The value of the mean square error of the estimate is calculated by the formula
\[
\Delta\ld(f;\widehat{A}_2\vec\xi\rd)= s^{-2}\|(\Phi_0+(1-\alpha)(D_0\Phi_0+\Phi_1))^{\top}\overline{\me 1}\|^2+(1-\alpha)^2s^{-2}\|\Phi_0^{\top} \overline{\me 1}\|^2.
\]

In a particular case   $d_0=1$, $d_1=1$ and, respectively, $D_0=0$, $D_1=0$, we have $\chi_{_{(1,u)}}^{(D)}(e^{-i\lambda})\equiv 1$,  and $G^{\pm}_{k^*}(k)=0$ for $k\geq1$.
In this case the estimate of the functional $A_{2s} \xi$ and the value of the its mean square error  are calculated by the formulas
\begin{eqnarray*}
\widehat{A}_2\vec \xi& = &-s^{-1}(\overline{\me 1})^{\top}\vec{\xi}(-u-1)+\alpha s^{-1}(\overline{\me 1})^{\top}\vec{\xi}(-u)+(1-\alpha) s^{-1}(\overline{\me 1})^{\top}\vec{\xi}(-u+1)
\\
& \quad &+s^{-1}\sum_{k=1}^{\infty}(-1)^{k+1}(\overline{\me 1})^{\top}(\Phi_1\Phi_0^{-1})^k
\\
& \quad &\quad\times
\ld(\vec \xi(-k)-\vec \xi(-k-1)-\vec\xi(-k-u)+\vec\xi(-k-u-1)\rd),
\end{eqnarray*}
and
\begin{eqnarray*}
\Delta\ld(f;\widehat{A}_2\vec\xi\rd)& = & s^{-2}\Bigg(\sum_{k=1}^s\ld(1-(1-\alpha)^{\delta_{ks}}a_0-(1-\alpha) a_k\rd)^2\Bigg.
\\& \quad &\Bigg.+(1-\alpha)^2(s-1)(1-a_0)^2+(1-\alpha)^2 \Bigg).
\end{eqnarray*}
\end{pry}

\section{Minimax (robust) method of forecasting}\label{minimax_extrapolation}

Values of the mean square errors and spectral characteristics of the optimal estimates
of functionals ${A}\vec\xi$ and ${A}_N\vec\xi$
constructed from unobserved values of stochastic sequence $\vec{\xi}(m)$ which determine a stationary stochastic GM increment sequence
$\chi_{\overline{\mu},\overline{s}}^{(d)}(\vec{\xi}(m))$ with the spectral density matrix $f(\lambda)$
based on its observations
$\vec\xi(m)$ at points $m=-1,-2,\dots$ can be calculated by formulas
  \eqref{pohybka},  \eqref{spectr A_e_d}  and
 \eqref{pohybka_N}, \eqref{est_h_N}
respectively, in the case where the spectral density matrix
$f(\lambda)$  is exactly known.
In the case where the spectral density $f(\lambda)$ admits the canonical factorization (\ref{fact_f_e_d}), formulas
 (\ref{simple_poh A_e_d}), (\ref{simple_spectr A_e_d}) and
 (\ref{simple_poh A_N_e_d}), (\ref{simple_spectr A_N_e_d})
are derived for calculating values of the mean square errors and spectral characteristics, respectively.

In many practical cases, however, complete information about the spectral density matrix is impossible while some sets $\md D$ of admissible spectral densities can be defined.
In this case the minimax method of estimation of functionals from  unobserved values of stochastic sequences is reasonable.
This method consists in finding an estimate that minimizes
the maximal values of the mean square errors for all spectral densities
from a given class $\mathcal D$ of admissible spectral densities simultaneously.

\begin{ozn}
For a given class of spectral densities $\mathcal{D}$ a
spectral density $f_0(\lambda)\in\mathcal{D}$ is called the least
favourable in $\mathcal{D}$ for the optimal linear estimation of the
functional $A\vec \xi$ if the following relation holds true:
\[\Delta(f_0)=\Delta(h_{\overline{\mu}}(f_0);f_0)=\max_{f\in\mathcal{D}}\Delta(h_{\overline{\mu}}(f);f).
\]
\end{ozn}

\begin{ozn}
 For a given class of spectral densities $\mathcal{D}$
a spectral characteristic $h^0(\lambda)$ of the optimal linear estimate of the functional
$A \xi$ is called minimax-robust if the following conditions are satisfied
\[h^0(\lambda)\in
H_{\mathcal{D}}=\bigcap_{f\in\mathcal{D}}L_2^{0-}(f),
\quad
\min_{h\in H_{\mathcal{D}}}\max_{f\in
\mathcal{D}}\Delta(h;f)=\sup_{f\in\mathcal{D}}\Delta(h^0;f).
\]
\end{ozn}

Taking into account the introduced definitions and the relations derived in the previous sections  we can verify that the following lemmas hold true.

\begin{lema}
A spectral density $f_0(\lambda)\in\mathcal{D}$ satisfying the minimality condition (\ref{umova11_e_d}) is the least favourable density in the class $\mathcal{D}$ for the optimal linear estimation  of the functional $A\vec\xi$ based on observations of the sequence $\vec\xi(m)$ at points $m=-1,-2,\ldots$ if the operator $\me F_{\overline{\mu}}^0$ defined by
the Fourier coefficients of the function
\be\label{functions_for_lemma_e_d}
 |\beta^{(d)}(i\lambda)|^2|\chi_{\overline{\mu}}^{(d)}(e^{-i\lambda})|^{-2}f_0^{-1}(\lambda),
 \ee
determines a solution to the constrained optimization problem
\be
 \max_{f\in \mathcal{D}}\ld(\ld\langle D^{\overline{\mu}}\me a,\me F_{\overline{\mu}}^{-1}D^{\overline{\mu}}\me a\rd\rangle\rd)= \ld\langle D^{\overline{\mu}}\me a,(\me F^0_{\overline{\mu}})^{-1}D^{\overline{\mu}}\me a\rd\rangle.
\label{minimax1_e_d}
\ee
The minimax spectral characteristic $h^0=h_{\overline{\mu}}(f^0)$ is calculated by formula (\ref{spectr A_e_d}) if
$h_{\overline{\mu}}(f^0)\in H_{\mathcal{D}}$.
\end{lema}

\begin{lema}
A spectral density $f_0(\lambda)\in\mathcal{D}$ which admits the
canonical factorization (\ref{fact_f_e_d}) is the least favourable density in the class
$\mathcal{D}$ for the optimal linear estimation  of the functional $A\vec\xi$ based on observations of the sequence $\vec\xi(m)$ at points $m=-1,-2,\ldots$
if coefficients $\{\varphi^0(k):k\geq0\}$
of the canonical factorization
\begin{equation}\label{minmax f_e_d}
 f_0(\lambda)=\ld(\sum_{k=0}^{\infty}\varphi^0(k)e^{-i\lambda k}\rd)\ld(\sum_{k=0}^{\infty}\varphi^0(k)e^{-i\lambda k}\rd)^*
  \end{equation}
of the spectral density $f^0(\lambda)$ determine a solution to the constrained optimization problem
\begin{equation}
\ld\|D^{\overline{\mu}}\me A\varphi_{\overline{\mu}}\rd\|^2\to\max,\,
f(\lambda)=
\ld(\sum_{k=0}^{\infty}\varphi(k)e^{-i\lambda k}\rd)\ld(\sum_{k=0}^{\infty}\varphi(k)e^{-i\lambda k}\rd)^*
\in\mathcal{D}.
\label{zadumextpextr_A_e_d}\end{equation}
The minimax spectral characteristic $h^0=h_{\overline{\mu}}(f_0)$ is calculated by formula (\ref{simple_spectr A_e_d}) if
$h_{\overline{\mu}}(f^0)\in H_{\mathcal{D}}$.
\end{lema}

\begin{lema}
A spectral density $f_0(\lambda)\in\mathcal{D}$ which admits
the canonical factorization (\ref{fact_f_e_d}) is the least favourable density in the class $\mathcal{D}$ for the optimal linear extrapolation of the functional $A_N\vec\xi$ based on observations of the sequence $\vec\xi(m)$ at points $m=-1,-2,\ldots$ if coefficients $\{\varphi^0(k):k=0,1,\ldots,N\}$
from the canonical factorization
\begin{equation}\label{minmax f_N_e_d}
f_0(\lambda)=\ld(\sum_{k=0}^{N}\varphi^0(k)e^{-i\lambda k}\rd)\ld(\sum_{k=0}^{N}\varphi^0(k)e^{-i\lambda k}\rd)^*
 \end{equation}
of the spectral density $f_0(\lambda)$ determine a solution to the constrained optimization problem
\begin{equation}
\ld\|D^{\overline{\mu}}_N\me A_N\varphi_{\mu,N}\rd\|^2\to\max,\,
f(\lambda)=
\ld(\sum_{k=0}^{N}\varphi(k)e^{-i\lambda k}\rd)\ld(\sum_{k=0}^{N}\varphi(k)e^{-i\lambda k}\rd)^*
\in\mathcal{D}.
\label{zadumextpextr_A_N_e_d}\end{equation}
The minimax spectral characteristic $h^0=h_{\overline{\mu}}(f_0)$ is calculated by formula (\ref{simple_spectr A_N_e_d}) if
$h_{\mu,N}(f^0)\in H_{\mathcal{D}}$.
\end{lema}

For more detailed analysis of properties of the least favorable spectral densities and the minimax-robust spectral characteristics we observe that the minimax spectral characteristic $h^0$ and the least favourable spectral density $f_0$ form a saddle
point of the function $\Delta(h;f)$ on the set
$H_{\mathcal{D}}\times\mathcal{D}$.
The saddle point inequalities
\[
 \Delta(h;f_0)\geq\Delta(h^0;f_0)\geq\Delta(h^0;f)\quad\forall f\in
 \mathcal{D},\forall h\in H_{\mathcal{D}}\]
hold true if $h^0=h_{\overline{\mu}}(f_0)$,
$h_{\overline{\mu}}(f_0)\in H_{\mathcal{D}}$ and $f_0$ is a solution of the constrained optimization problem
\be
 \widetilde{\Delta}(f)=-\Delta(h_{\overline{\mu}}(f_0);f)\to
 \inf,\quad f\in \mathcal{D},\label{zad_um_extr_e_d}
 \ee
where the functional $\Delta(h_{\overline{\mu}}(f_0);f)$ is calculated by the formula

\begin{eqnarray}
 \Delta(h_{\overline{\mu}}(f_0);f)
 & = & \frac{1}{2\pi}\int_{-\pi}^{\pi}
\frac{\overline{\beta^{(d)}(i\lambda)}
}
{\overline{\chi_{\overline{\mu}}^{(d)}(e^{-i\lambda})}}
\left(
\sum_{k=0}^{\infty}((\me F^0_{\overline{\mu}})^{-1}D^{\overline{\mu}}\me a )_k e^{i\lambda k}\right)^{\top}f_0^{-1}(\lambda)
f(\lambda)
\nonumber
\\
 & \quad & \times f_0^{-1}(\lambda)\left(\overline{
\sum_{k=0}^{\infty}((\me F^0_{\overline{\mu}})^{-1}D^{\overline{\mu}}\me a )_k e^{i\lambda k}}\right)
\frac{\beta^{(d)}(i\lambda)}
{\chi_{\overline{\mu}}^{(d)}(e^{-i\lambda})}d\lambda
\end{eqnarray}
or by the formula
\begin{eqnarray}
\nonumber
\Delta(h_{\overline{\mu}}(f_0);f)
 & = & \frac{1}{2\pi}\ip
\frac{\chi_{\overline{\mu}}^{(d)}(e^{-i\lambda})}{\beta^{(d)}(i\lambda)}
\left(
 \sum_{k=0}^{\infty}(D^{\overline{\mu}}\me A\varphi^0_{\overline{\mu}})_ke^{i\lambda k}
\right)^{\top}
(\Psi^0_{\overline{\mu}}(e^{-i\lambda}))f(\lambda)
\nonumber
\\
& \quad & \times
(\Psi^0_{\overline{\mu}}(e^{-i\lambda}))^{*}
\left(\overline{
 \sum_{k=0}^{\infty}(D^{\overline{\mu}}\me A\varphi^0_{\overline{\mu}})_ke^{i\lambda k}}
\right)\frac{\overline{\chi_{\overline{\mu}}^{(d)}(e^{-i\lambda})}}{\overline{\beta^{(d)}(i\lambda)}}d\lambda\label{simple_poh A_e_d}
\end{eqnarray}
in the case where the spectral density admits the canonical factorization (\ref{fact_f_e_d}).

The constrained optimization problem (\ref{zad_um_extr_e_d}) is equivalent to the unconstrained optimization problem
\[
 \Delta_{\mathcal{D}}(f)=\widetilde{\Delta}(f)+ \delta(f|\mathcal{D})\to\inf,\]
where $\delta(f|\mathcal{D})$ is the indicator function of the set
$\mathcal{D}$, namely $\delta(f|\mathcal{D})=0$ if $f\in \mathcal{D}$ and $\delta(f|\mathcal{D})=+\infty$ if $f\notin \mathcal{D}$.
 A solution $f_0$ of this unconstrained optimization problem is characterized by the condition
 $0\in\partial\Delta_{\mathcal{D}}(f_0)$, which is the necessary and sufficient condition under which a point $f_0$ belongs to the set of minimums of the convex functional $\Delta_{\mathcal{D}}(f)$ \cite{Franke1985,Mokl_Mas_extr,Moklyachuk,Rockafellar}.
 This condition makes it possible to find the least favourable spectral densities in some special classes of spectral densities $D$.

The form of the functional $\widetilde{\Delta}(f)$ allows us to apply the Lagrange method of indefinite
multipliers for investigating the constrained optimization problem (\ref{zad_um_extr_e_d}).
The complexity of optimization problem is determined by the complexity of calculating the subdifferentials of the indicator functions of sets of admissible spectral densities.

\subsection{Least favorable spectral density in classes $\md D_0$}

Consider the forecasting problem for the functional $A\vec{\xi}$
 which depends on unobserved values of a sequence $\vec\xi(m)$ with stationary GM increments based on observations of the sequence at points $m=-1,-2,\ldots$ under the condition that sets of admissible spectral densities $\md D_0^k,k=1,2,3,4$ are defined as follows:
$$\md D_{0}^{1} =\bigg\{f(\lambda )\left|\frac{1}{2\pi} \int
_{-\pi}^{\pi}
\frac{|\chi_{\overline{\mu}}^{(d)}(e^{-i\lambda})|^2}{|\beta^{(d)}(i\lambda)|^2}
f(\lambda )d\lambda  =P\right.\bigg\},$$
$$\md D_{0}^{2} =\bigg\{f(\lambda )\left|\frac{1}{2\pi }
\int _{-\pi }^{\pi}
\frac{|\chi_{\overline{\mu}}^{(d)}(e^{-i\lambda})|^2}{|\beta^{(d)}(i\lambda)|^2}
{\rm{Tr}}\,[ f(\lambda )]d\lambda =p\right.\bigg\},$$
$$\md D_{0}^{3} =\bigg\{f(\lambda )\left|\frac{1}{2\pi }
\int _{-\pi}^{\pi}
\frac{|\chi_{\overline{\mu}}^{(d)}(e^{-i\lambda})|^2}{|\beta^{(d)}(i\lambda)|^2}
f_{kk} (\lambda )d\lambda =p_{k}, k=\overline{1,T}\right.\bigg\},$$
$$\md D_{0}^{4} =\bigg\{f(\lambda )\left|\frac{1}{2\pi} \int _{-\pi}^{\pi}
\frac{|\chi_{\overline{\mu}}^{(d)}(e^{-i\lambda})|^2}{|\beta^{(d)}(i\lambda)|^2}
\left\langle B_{1} ,f(\lambda )\right\rangle d\lambda  =p\right.\bigg\},$$

\noindent
where  $p, p_k, k=\overline{1,T}$ are given numbers, $P, B_1,$ are given positive-definite Hermitian matrices.

From the condition $0\in\partial\Delta_{\mathcal{D}}(f_0)$
we find the following equations which determine the least favourable spectral densities for these given sets of admissible spectral densities.

For the first set of admissible spectral densities $\md D_0^1$ we have equation
\begin{multline} \label{eq_4_1}
\left(
\sum_{k=0}^{\infty}((\me F^0_{\overline{\mu}})^{-1}D^{\overline{\mu}}\me a )_k e^{i\lambda k}
\right)
\left({
\sum_{k=0}^{\infty}((\me F^0_{\overline{\mu}})^{-1}D^{\overline{\mu}}\me a )_k e^{i\lambda k}}
\right)^{*}=
\\
=\left(\frac{|\chi_{\overline{\mu}}^{(d)}(e^{-i\lambda})|^2}{|\beta^{(d)}(i\lambda)|^2} f_{0} (\lambda )\right)
\vec{\alpha}\cdot \vec{\alpha}^{*}\left(\frac{|\chi_{\overline{\mu}}^{(d)}(e^{-i\lambda})|^2}{|\beta^{(d)}(i\lambda)|^2} f_{0} (\lambda )\right),
\end{multline}
where $\vec{\alpha}$ is a vector of Lagrange multipliers.

For the second set of admissible spectral densities $\md D_0^2$ we have equation
\begin{multline}  \label{eq_4_2}
\left(
\sum_{k=0}^{\infty}((\me F^0_{\overline{\mu}})^{-1}D^{\overline{\mu}}\me a )_k e^{i\lambda k}
\right)
\left({
\sum_{k=0}^{\infty}((\me F^0_{\overline{\mu}})^{-1}D^{\overline{\mu}}\me a )_k e^{i\lambda k}}
\right)^{*}=
\\=
\alpha^{2} \left(\frac{|\chi_{\overline{\mu}}^{(d)}(e^{-i\lambda})|^2}{|\beta^{(d)}(i\lambda)|^2} f_{0} (\lambda )\right)^{2},
\end{multline}
where $\alpha^{2}$ is a Lagrange multiplier.

For the third set of admissible spectral densities $\md D_0^3$ we have equation

\begin{multline}   \label{eq_4_3}
\left(
\sum_{k=0}^{\infty}((\me F^0_{\overline{\mu}})^{-1}D^{\overline{\mu}}\me a )_k e^{i\lambda k}
\right)
\left({
\sum_{k=0}^{\infty}((\me F^0_{\overline{\mu}})^{-1}D^{\overline{\mu}}\me a )_k e^{i\lambda k}}
\right)^{*}=
\\
=
\left(\frac{|\chi_{\overline{\mu}}^{(d)}(e^{-i\lambda})|^2}{|\beta^{(d)}(i\lambda)|^2} f_{0} (\lambda )\right)\left\{\alpha _{k}^{2} \delta _{kl} \right\}_{k,l=1}^{T} \left(\frac{|\chi_{\overline{\mu}}^{(d)}(e^{-i\lambda})|^2}{|\beta^{(d)}(i\lambda)|^2} f_{0} (\lambda )\right),
\end{multline}
where  $\alpha _{k}^{2}$ are Lagrange multipliers, $\delta _{kl}$ are Kronecker symbols.

For the fourth set of admissible spectral densities $\md D_0^4$  we have equation
\begin{multline}   \label{eq_4_4}
\left(
\sum_{k=0}^{\infty}((\me F^0_{\overline{\mu}})^{-1}D^{\overline{\mu}}\me a )_k e^{i\lambda k}
\right)
\left({
\sum_{k=0}^{\infty}((\me F^0_{\overline{\mu}})^{-1}D^{\overline{\mu}}\me a )_k e^{i\lambda k}}
\right)^{*}=
\\
=
\alpha^{2} \left(\frac{|\chi_{\overline{\mu}}^{(d)}(e^{-i\lambda})|^2}{|\beta^{(d)}(i\lambda)|^2} f_{0} (\lambda )\right)B_{1}^{\top} \left(\frac{|\chi_{\overline{\mu}}^{(d)}(e^{-i\lambda})|^2}{|\beta^{(d)}(i\lambda)|^2} f_{0} (\lambda )\right),
\end{multline}
where $\alpha^{2}$ is a Lagrange multiplier.

In the case where the spectral density admits the canonical factorization (\ref{fact_f_e_d}) we have the following equations, correspondingly
\begin{equation}  \label{eq_4_11}
\left(
 \sum_{k=0}^{\infty}(D^{\overline{\mu}}\me A\varphi^0_{\overline{\mu}})_ke^{i\lambda k}
\right)
\left(
 \sum_{k=0}^{\infty}(D^{\overline{\mu}}\me A\varphi^0_{\overline{\mu}})_ke^{i\lambda k}
\right)^{*}
=(\Phi^0_{\overline{\mu}}(e^{-i\lambda}))^{\top}
\vec{\alpha}\cdot \vec{\alpha}^{*}
\overline{\Phi^0_{\overline{\mu}}(e^{-i\lambda})},
\end{equation}
\begin{equation}   \label{eq_4_21}
\left(
 \sum_{k=0}^{\infty}(D^{\overline{\mu}}\me A\varphi^0_{\overline{\mu}})_ke^{i\lambda k}
\right)
\left(
 \sum_{k=0}^{\infty}(D^{\overline{\mu}}\me A\varphi^0_{\overline{\mu}})_ke^{i\lambda k}
\right)^{*}
=
\alpha^{2}(\Phi^0_{\overline{\mu}}(e^{-i\lambda}))^{\top}\overline{\Phi^0_{\overline{\mu}}(e^{-i\lambda})}.
\end{equation}
\begin{multline}   \label{eq_4_31}
\left(
 \sum_{k=0}^{\infty}(D^{\overline{\mu}}\me A\varphi^0_{\overline{\mu}})_ke^{i\lambda k}
\right)
\left(
 \sum_{k=0}^{\infty}(D^{\overline{\mu}}\me A\varphi^0_{\overline{\mu}})_ke^{i\lambda k}
\right)^{*}=
\\=(\Phi^0_{\overline{\mu}}(e^{-i\lambda}))^{\top}\left\{\alpha _{k}^{2} \delta _{kl} \right\}_{k,l=1}^{T}
\overline{\Phi^0_{\overline{\mu}}(e^{-i\lambda})},
\end{multline}
\begin{equation}    \label{eq_4_41}
\left(
 \sum_{k=0}^{\infty}(D^{\overline{\mu}}\me A\varphi^0_{\overline{\mu}})_ke^{i\lambda k}
\right)
\left(
 \sum_{k=0}^{\infty}(D^{\overline{\mu}}\me A\varphi^0_{\overline{\mu}})_ke^{i\lambda k}
\right)^{*}
=
\alpha^{2} (\Phi^0_{\overline{\mu}}(e^{-i\lambda}))^{\top}B_{1}^{\top}  \overline{\Phi^0_{\overline{\mu}}(e^{-i\lambda})},
\end{equation}

The following theorem  holds true.

\begin{thm}
Let the minimality condition (\ref{umova11_e_d}) hold true. The least favorable spectral densities $f_{0}(\lambda)$ in the classes $ \md  D_0^{k}$, $k=1,2,3,4$, for the optimal linear estimation  of the functional  $A\vec{\xi}$ from observations of the sequence $\vec{\xi}(m)$ at points  $m=-1,-2,\ldots$  are determined by  equations
\eqref{eq_4_1}, \eqref{eq_4_2}, \eqref{eq_4_3}, \eqref{eq_4_4},
(or equations \eqref{eq_4_11}, \eqref{eq_4_21}, \eqref{eq_4_31}, \eqref{eq_4_41}
in the case where the spectral densities admit the canonical factorization (\ref{fact_f_e_d}), respectively),
the constrained optimization problem (\ref{minimax1_e_d}) and restrictions  on densities from the corresponding classes $ \md  D_0^{k}$, $k=1,2,3,4$.  The minimax-robust spectral characteristic of the optimal estimate of the functional $A\vec{\xi}$ is determined by the formula (\ref{spectr A_e_d}).
\end{thm}

\subsection{Least favorable spectral density in classes $\md D_V^U$}

Consider the forecasting problem for the functional $A\vec{\xi}$
 which depends on unobserved values of a sequence $\vec\xi(m)$ with stationary GM increments based on observations of the sequence at points $m=-1,-2,\ldots$ under the condition that sets of admissible spectral densities $ {\md D_{V}^{U}} ^{k},k=1,2,3,4$ are defined as follows:

\begin{equation*}
 {\md D_{V}^{U}} ^{1}=\left\{f(\lambda )\bigg|V(\lambda )\le f(\lambda
)\le U(\lambda ), \frac{1}{2\pi } \int _{-\pi}^{\pi}
\frac{|\chi_{\overline{\mu}}^{(d)}(e^{-i\lambda})|^2}{|\beta^{(d)}(i\lambda)|^2}
f(\lambda )d\lambda=Q\right\},
\end{equation*}
\begin{multline*}
  {\md D_{V}^{U}} ^{2}  =\bigg\{f(\lambda )\bigg|{\mathrm{Tr}}\, [V(\lambda
)]\le {\mathrm{Tr}}\,[ f(\lambda )]\le {\mathrm{Tr}}\, [U(\lambda )],
\\
\frac{1}{2\pi } \int _{-\pi}^{\pi}
\frac{|\chi_{\overline{\mu}}^{(d)}(e^{-i\lambda})|^2}{|\beta^{(d)}(i\lambda)|^2}
{\mathrm{Tr}}\,  [f(\lambda)]d\lambda  =q \bigg\},
\end{multline*}
\begin{multline*}
{\md D_{V}^{U}} ^{3}  =\bigg\{f(\lambda )\bigg|v_{kk} (\lambda )  \le
f_{kk} (\lambda )\le u_{kk} (\lambda ),
\\
\frac{1}{2\pi} \int _{-\pi}^{\pi}
\frac{|\chi_{\overline{\mu}}^{(d)}(e^{-i\lambda})|^2}{|\beta^{(d)}(i\lambda)|^2}
f_{kk} (\lambda
)d\lambda  =q_{k} , k=\overline{1,T}\bigg\},
\end{multline*}
\begin{multline*}
{\md D_{V}^{U}} ^{4}  =\bigg\{f(\lambda )\bigg|\left\langle B_{2}
,V(\lambda )\right\rangle \le \left\langle B_{2},f(\lambda
)\right\rangle \le \left\langle B_{2} ,U(\lambda)\right\rangle,
\\
\frac{1}{2\pi }
\int _{-\pi}^{\pi}
\frac{|\chi_{\overline{\mu}}^{(d)}(e^{-i\lambda})|^2}{|\beta^{(d)}(i\lambda)|^2}
\left\langle B_{2},f(\lambda)\right\rangle d\lambda  =q\bigg\}.
\end{multline*}

\noindent
Here the spectral densities $V( \lambda ),U( \lambda )$ are known and fixed, $ q,  q_k, k=\overline{1,T}$ are given numbers, $ Q, B_2$ are given positive definite Hermitian matrices.

From the condition $0\in\partial\Delta_{\mathcal{D}}(f_0)$
we find the following equations which determine the least favourable spectral densities for these given sets of admissible spectral densities.

For the first set of admissible spectral densities ${\md D_{V}^{U}} ^{1}$ we have equation
\begin{multline} \label{eq_5_1}
\left(
\sum_{k=0}^{\infty}((\me F^0_{\overline{\mu}})^{-1}D^{\overline{\mu}}\me a )_k e^{i\lambda k}
\right)
\left({
\sum_{k=0}^{\infty}((\me F^0_{\overline{\mu}})^{-1}D^{\overline{\mu}}\me a )_k e^{i\lambda k}}
\right)^{*}=
\\
=\left(\frac{|\chi_{\overline{\mu}}^{(d)}(e^{-i\lambda})|^2}{|\beta^{(d)}(i\lambda)|^2} f_{0} (\lambda )\right)(\vec{\beta}\cdot \vec{\beta}^{*}+\Gamma _{1} (\lambda )+\Gamma _{2} (\lambda ))\left(\frac{|\chi_{\overline{\mu}}^{(d)}(e^{-i\lambda})|^2}{|\beta^{(d)}(i\lambda)|^2} f_{0} (\lambda )\right),
\end{multline}
where $ \vec{\beta}$ is a vector of Lagrange multipliers, $\Gamma _{1} (\lambda )\le 0$ and $\Gamma _{1} (\lambda )=0$ if $f_{0}(\lambda )>V(\lambda ),$ $
\Gamma _{2} (\lambda )\ge 0$ and $\Gamma _{2} (\lambda )=0$ if $f_{0}(\lambda )<U(\lambda ).$

For the second set of admissible spectral densities ${\md D_{V}^{U}} ^{2}$ we have equation
\begin{multline}   \label{eq_5_2}
\left(
\sum_{k=0}^{\infty}((\me F^0_{\overline{\mu}})^{-1}D^{\overline{\mu}}\me a )_k e^{i\lambda k}
\right)
\left({
\sum_{k=0}^{\infty}((\me F^0_{\overline{\mu}})^{-1}D^{\overline{\mu}}\me a )_k e^{i\lambda k}}
\right)^{*}=
\\=
(\beta^{2} +\gamma _{1} (\lambda )+\gamma _{2} (\lambda )) \left(\frac{|\chi_{\overline{\mu}}^{(d)}(e^{-i\lambda})|^2}{|\beta^{(d)}(i\lambda)|^2} f_{0} (\lambda )\right)^{2},
\end{multline}
where $ \beta^{2}$ is Lagrange multiplier,  $\gamma _{1} (\lambda )\le 0$ and $\gamma _{1} (\lambda )=0$ if ${\mathrm{Tr}}\,
[f_{0} (\lambda )]> {\mathrm{Tr}}\,  [V(\lambda )],$ $\gamma _{2} (\lambda )\ge 0$ and $\gamma _{2} (\lambda )=0$ if $ {\mathrm{Tr}}\,[f_{0}(\lambda )]< {\mathrm{Tr}}\, [ U(\lambda)].$

For the third set of admissible spectral densities ${\md D_{V}^{U}} ^{3}$ we have equation
\begin{multline}   \label{eq_5_3}
\left(
\sum_{k=0}^{\infty}((\me F^0_{\overline{\mu}})^{-1}D^{\overline{\mu}}\me a )_k e^{i\lambda k}
\right)
\left({
\sum_{k=0}^{\infty}((\me F^0_{\overline{\mu}})^{-1}D^{\overline{\mu}}\me a )_k e^{i\lambda k}}
\right)^{*}=\left(\frac{|\chi_{\overline{\mu}}^{(d)}(e^{-i\lambda})|^2}{|\beta^{(d)}(i\lambda)|^2} f_{0} (\lambda )\right)
\\
\times
\left\{(\beta_{k}^{2} +\gamma _{1k} (\lambda )+\gamma _{2k} (\lambda ))\delta _{kl}\right\}_{k,l=1}^{T}
\left(\frac{|\chi_{\overline{\mu}}^{(d)}(e^{-i\lambda})|^2}{|\beta^{(d)}(i\lambda)|^2} f_{0} (\lambda )\right),
\end{multline}
where  $ \beta_{k}^{2}$ are Lagrange multipliers, $\delta _{kl}$ are Kronecker symbols, $\gamma _{1k} (\lambda )\le 0$ and $\gamma _{1k} (\lambda )=0$ if $f_{kk}^{0} (\lambda )>v_{kk} (\lambda ),$ $\gamma _{2k} (\lambda )\ge 0$ and $\gamma _{2k} (\lambda )=0$ if $f_{kk}^{0} (\lambda )<u_{kk} (\lambda).$

For the fourth set of admissible spectral densities ${\md D_{V}^{U}} ^{4}$ we have equation
\begin{multline}  \label{eq_5_4}
\left(
\sum_{k=0}^{\infty}((\me F^0_{\overline{\mu}})^{-1}D^{\overline{\mu}}\me a )_k e^{i\lambda k}
\right)
\left({
\sum_{k=0}^{\infty}((\me F^0_{\overline{\mu}})^{-1}D^{\overline{\mu}}\me a )_k e^{i\lambda k}}
\right)^{*}=
\\
=
(\beta^{2} +\gamma'_{1}(\lambda )+\gamma'_{2}(\lambda )) \left(\frac{|\chi_{\overline{\mu}}^{(d)}(e^{-i\lambda})|^2}{|\beta^{(d)}(i\lambda)|^2} f_{0} (\lambda )\right) B_2^\top \left(\frac{|\chi_{\overline{\mu}}^{(d)}(e^{-i\lambda})|^2}{|\beta^{(d)}(i\lambda)|^2} f_{0} (\lambda )\right)
\end{multline}
where $ \beta^{2}$ is Lagrange multiplier, $\gamma'_{1}( \lambda )\le 0$ and $\gamma'_{1} ( \lambda )=0$ if $\langle B_{2},f_{0} ( \lambda) \rangle > \langle B_{2},V( \lambda ) \rangle,$ $\gamma'_{2}( \lambda )\ge 0$ and $\gamma'_{2} ( \lambda )=0$ if $\langle
B_{2} ,f_{0} ( \lambda) \rangle < \langle B_{2} ,U( \lambda ) \rangle.$

In the case where the spectral density admits the canonical factorization (\ref{fact_f_e_d}) we have the following equations, correspondingly
\begin{multline} \label{eq_5_11}
\left(
 \sum_{k=0}^{\infty}(D^{\overline{\mu}}\me A\varphi^0_{\overline{\mu}})_ke^{i\lambda k}
\right)
\left(
 \sum_{k=0}^{\infty}(D^{\overline{\mu}}\me A\varphi^0_{\overline{\mu}})_ke^{i\lambda k}
\right)^{*}=
\\=(\Phi^0_{\overline{\mu}}(e^{-i\lambda}))^{\top}(\vec{\beta}\cdot \vec{\beta}^{*}+\Gamma _{1} (\lambda )+\Gamma _{2} (\lambda ))\,\,\overline{\Phi^0_{\overline{\mu}}(e^{-i\lambda})},
\end{multline}
\begin{multline}
\label{eq_5_21}
\left(
 \sum_{k=0}^{\infty}(D^{\overline{\mu}}\me A\varphi^0_{\overline{\mu}})_ke^{i\lambda k}
\right)
\left(
 \sum_{k=0}^{\infty}(D^{\overline{\mu}}\me A\varphi^0_{\overline{\mu}})_ke^{i\lambda k}
\right)^{*}=
\\=
(\beta^{2} +\gamma _{1} (\lambda )+\gamma _{2} (\lambda )) (\Phi^0_{\overline{\mu}}(e^{-i\lambda}))^{\top}\,\,\overline{\Phi^0_{\overline{\mu}}(e^{-i\lambda})},
\end{multline}
\begin{multline}  \label{eq_5_31}
\left(
 \sum_{k=0}^{\infty}(D^{\overline{\mu}}\me A\varphi^0_{\overline{\mu}})_ke^{i\lambda k}
\right)
\left(
 \sum_{k=0}^{\infty}(D^{\overline{\mu}}\me A\varphi^0_{\overline{\mu}})_ke^{i\lambda k}
\right)^{*}=
\\=
(\Phi^0_{\overline{\mu}}(e^{-i\lambda}))^{\top}\left\{(\beta_{k}^{2} +\gamma _{1k} (\lambda )+\gamma _{2k} (\lambda ))\delta _{kl}\right\}_{k,l=1}^{T}\overline{\Phi^0_{\overline{\mu}}(e^{-i\lambda})},
\end{multline}
\begin{multline} \label{eq_5_41}
\left(
 \sum_{k=0}^{\infty}(D^{\overline{\mu}}\me A\varphi^0_{\overline{\mu}})_ke^{i\lambda k}
\right)
\left(
 \sum_{k=0}^{\infty}(D^{\overline{\mu}}\me A\varphi^0_{\overline{\mu}})_ke^{i\lambda k}
\right)^{*}=
\\=
(\beta^{2} +\gamma'_{1}(\lambda )+\gamma'_{2}(\lambda )) (\Phi^0_{\overline{\mu}}(e^{-i\lambda}))^{\top}B_2^\top\,\, \overline{\Phi^0_{\overline{\mu}}(e^{-i\lambda})},
\end{multline}

The following theorem  holds true.

\begin{thm}
Let the minimality condition (\ref{umova11_e_d}) hold true. The least favorable spectral densities $f_{0}(\lambda)$ in the classes ${\md D_{V}^{U}} ^{k}$, $k=1,2,3,4$, for the optimal linear estimation  of the functional  $A\vec{\xi}$ from observations of the sequence $\vec{\xi}(m)$ at points  $m=-1,-2,\ldots$  are determined by  equations
\eqref{eq_5_1}, \eqref{eq_5_2}, \eqref{eq_5_3}, \eqref{eq_5_4}
(or equations \eqref{eq_5_11}, \eqref{eq_5_21}, \eqref{eq_5_31}, \eqref{eq_5_41}
in the case where the spectral densities admit the canonical factorization (\ref{fact_f_e_d}), respectively),
the constrained optimization problem (\ref{minimax1_e_d}) and restrictions  on densities from the corresponding classes ${\md D_{V}^{U}} ^{k}$, $k=1,2,3,4$.  The minimax-robust spectral characteristic of the optimal estimate of the functional $A\vec{\xi}$ is determined by the formula (\ref{spectr A_e_d}).
\end{thm}

\subsection{Least favorable spectral density in classes $\md D_{\varepsilon}$}

Consider the forecasting problem for the functional $A\vec{\xi}$
 which depends on unobserved values of a sequence $\vec\xi(m)$ with stationary GM increments based on observations of the sequence at points $m=-1,-2,\ldots$ under the condition that sets of admissible spectral densities $ \md D_{\varepsilon}^{k},k=1,2,3,4$ are defined as follows:

\begin{multline*}
\md D_{\varepsilon }^{1}  =\bigg\{f(\lambda )\bigg|{\mathrm{Tr}}\,
[f(\lambda )]=(1-\varepsilon ) {\mathrm{Tr}}\,  [f_{1} (\lambda
)]+\varepsilon {\mathrm{Tr}}\,  [W(\lambda )],
\\
\frac{1}{2\pi} \int _{-\pi}^{\pi}
\frac{|\chi_{\overline{\mu}}^{(d)}(e^{-i\lambda})|^2}{|\beta^{(d)}(i\lambda)|^2}
{\mathrm{Tr}}\,
[f(\lambda )]d\lambda =p \bigg\};
\end{multline*}
\begin{multline*}
\md D_{\varepsilon }^{2}  =\bigg\{f(\lambda )\bigg|f_{kk} (\lambda)
=(1-\varepsilon )f_{kk}^{1} (\lambda )+\varepsilon w_{kk}(\lambda),
\\
\frac{1}{2\pi} \int _{-\pi}^{\pi}
\frac{|\chi_{\overline{\mu}}^{(d)}(e^{-i\lambda})|^2}{|\beta^{(d)}(i\lambda)|^2}
f_{kk} (\lambda)d\lambda  =p_{k} , k=\overline{1,T}\bigg\};
\end{multline*}
\begin{multline*}
\md D_{\varepsilon }^{3} =\bigg\{f(\lambda )\bigg|\left\langle B_{1},f(\lambda )\right\rangle =(1-\varepsilon )\left\langle B_{1},f_{1} (\lambda )\right\rangle+\varepsilon \left\langle B_{1},W(\lambda )\right\rangle,
\\
\frac{1}{2\pi}\int _{-\pi}^{\pi}
\frac{|\chi_{\overline{\mu}}^{(d)}(e^{-i\lambda})|^2}{|\beta^{(d)}(i\lambda)|^2}
\left\langle B_{1} ,f(\lambda )\right\rangle d\lambda =p\bigg\};
\end{multline*}
\begin{multline*}
\md D_{\varepsilon }^{4}=\bigg\{f(\lambda )\bigg|f(\lambda)=(1-\varepsilon )f_{1} (\lambda )+\varepsilon W(\lambda ),
\\
\frac{1}{2\pi } \int _{-\pi}^{\pi}
\frac{|\chi_{\overline{\mu}}^{(d)}(e^{-i\lambda})|^2}{|\beta^{(d)}(i\lambda)|^2}
f(\lambda )d\lambda=P\bigg\}.
\end{multline*}

\noindent
Here  $f_{1} ( \lambda )$ is a fixed spectral density, $W(\lambda)$ is an unknown spectral density, $p,  p_k, k=\overline{1,T}$, are given numbers, $P$ is a given positive-definite Hermitian matrices.

From the condition $0\in\partial\Delta_{\mathcal{D}}(f_0)$
we find the following equations which determine the least favourable spectral densities for these given sets of admissible spectral densities.

For the first set of admissible spectral densities $ \md D_{\varepsilon}^{1}$ we have equation
\begin{multline}\label{eq_6_1}
\left(
\sum_{k=0}^{\infty}((\me F^0_{\overline{\mu}})^{-1}D^{\overline{\mu}}\me a )_k e^{i\lambda k}
\right)
\left({
\sum_{k=0}^{\infty}((\me F^0_{\overline{\mu}})^{-1}D^{\overline{\mu}}\me a )_k e^{i\lambda k}}
\right)^{*}=
\\=
(\alpha^{2} +\gamma_1(\lambda ))\left(\frac{|\chi_{\overline{\mu}}^{(d)}(e^{-i\lambda})|^2}{|\beta^{(d)}(i\lambda)|^2} f_{0} (\lambda )\right)^{2},
\end{multline}
where $\alpha^{2}$ is Lagrange multiplier, $\gamma_1(\lambda )\le 0$ and $\gamma_1(\lambda )=0$ if ${\mathrm{Tr}}\,[f_{0} (\lambda )]>(1-\varepsilon ) {\mathrm{Tr}}\, [f_{1} (\lambda )]$.

For the second set of admissible spectral densities $ \md D_{\varepsilon}^{2}$ we have equation
\begin{multline}   \label{eq_6_2}
\left(
\sum_{k=0}^{\infty}((\me F^0_{\overline{\mu}})^{-1}D^{\overline{\mu}}\me a )_k e^{i\lambda k}
\right)
\left({
\sum_{k=0}^{\infty}((\me F^0_{\overline{\mu}})^{-1}D^{\overline{\mu}}\me a )_k e^{i\lambda k}}
\right)^{*}=
\\
=
\left(\frac{|\chi_{\overline{\mu}}^{(d)}(e^{-i\lambda})|^2}{|\beta^{(d)}(i\lambda)|^2} f_{0} (\lambda )\right)\left\{(\alpha_{k}^{2} +\gamma_{k}^1 (\lambda ))\delta _{kl} \right\}_{k,l=1}^{T} \left(\frac{|\chi_{\overline{\mu}}^{(d)}(e^{-i\lambda})|^2}{|\beta^{(d)}(i\lambda)|^2} f_{0} (\lambda )\right),
\end{multline}
where  $\alpha _{k}^{2}$ are Lagrange multipliers, $\gamma_{k}^1(\lambda )\le 0$ and $\gamma_{k}^1 (\lambda )=0$ if $f_{kk}^{0}(\lambda )>(1-\varepsilon )f_{kk}^{1} (\lambda )$.

For the third set of admissible spectral densities $ \md D_{\varepsilon}^{3}$ we have equation
\begin{multline}   \label{eq_6_3}
\left(
\sum_{k=0}^{\infty}((\me F^0_{\overline{\mu}})^{-1}D^{\overline{\mu}}\me a )_k e^{i\lambda k}
\right)
\left({
\sum_{k=0}^{\infty}((\me F^0_{\overline{\mu}})^{-1}D^{\overline{\mu}}\me a )_k e^{i\lambda k}}
\right)^{*}=
\\
=
(\alpha^{2} +\gamma_1'(\lambda ))\left(\frac{|\chi_{\overline{\mu}}^{(d)}(e^{-i\lambda})|^2}{|\beta^{(d)}(i\lambda)|^2} f_{0} (\lambda )\right) B_{1}^{\top} \left(\frac{|\chi_{\overline{\mu}}^{(d)}(e^{-i\lambda})|^2}{|\beta^{(d)}(i\lambda)|^2} f_{0} (\lambda )\right),
\end{multline}
where $\alpha^{2}$ is a Lagrange multiplier, $\gamma_1' ( \lambda )\le 0$ and $\gamma_1' ( \lambda )=0$ if $\langle B_{1} ,f_{0} ( \lambda ) \rangle>(1- \varepsilon ) \langle B_{1} ,f_{1} ( \lambda ) \rangle$.

For the fourth set of admissible spectral densities $ \md D_{\varepsilon}^{4}$ we have equation
\begin{multline}  \label{eq_6_4}
\left(
\sum_{k=0}^{\infty}((\me F^0_{\overline{\mu}})^{-1}D^{\overline{\mu}}\me a )_k e^{i\lambda k}
\right)
\left({
\sum_{k=0}^{\infty}((\me F^0_{\overline{\mu}})^{-1}D^{\overline{\mu}}\me a )_k e^{i\lambda k}}
\right)^{*}=
\\
=
\left(\frac{|\chi_{\overline{\mu}}^{(d)}(e^{-i\lambda})|^2}{|\beta^{(d)}(i\lambda)|^2} f_{0} (\lambda )\right) (\vec{\alpha}\cdot \vec{\alpha}^{*}+\Gamma(\lambda)) \left(\frac{|\chi_{\overline{\mu}}^{(d)}(e^{-i\lambda})|^2}{|\beta^{(d)}(i\lambda)|^2} f_{0} (\lambda )\right),
\end{multline}
where  $\vec{\alpha}$ is a vector of Lagrange multipliers, $\Gamma(\lambda )\le 0$ and $\Gamma(\lambda )=0$ if $f_{0}(\lambda )>(1-\varepsilon )f_{1} (\lambda )$.

In the case where the spectral density admits the canonical factorization (\ref{fact_f_e_d}) we have the following equations, correspondingly
\begin{multline} \label{eq_6_11}
\left(
 \sum_{k=0}^{\infty}(D^{\overline{\mu}}\me A\varphi^0_{\overline{\mu}})_ke^{i\lambda k}
\right)
\left(
 \sum_{k=0}^{\infty}(D^{\overline{\mu}}\me A\varphi^0_{\overline{\mu}})_ke^{i\lambda k}
\right)^{*}=
\\=
(\alpha^{2} +\gamma_1(\lambda ))(\Phi^0_{\overline{\mu}}(e^{-i\lambda}))^{\top}\overline{\Phi^0_{\overline{\mu}}(e^{-i\lambda})},
\end{multline}
\begin{multline}   \label{eq_6_21}
\left(
 \sum_{k=0}^{\infty}(D^{\overline{\mu}}\me A\varphi^0_{\overline{\mu}})_ke^{i\lambda k}
\right)
\left(
 \sum_{k=0}^{\infty}(D^{\overline{\mu}}\me A\varphi^0_{\overline{\mu}})_ke^{i\lambda k}
\right)^{*}=
\\=
(\Phi^0_{\overline{\mu}}(e^{-i\lambda}))^{\top}\left\{(\alpha_{k}^{2} +\gamma_{k}^1 (\lambda ))\delta _{kl} \right\}_{k,l=1}^{T} \overline{\Phi^0_{\overline{\mu}}(e^{-i\lambda})},
\end{multline}
\begin{multline}   \label{eq_6_31}
\left(
 \sum_{k=0}^{\infty}(D^{\overline{\mu}}\me A\varphi^0_{\overline{\mu}})_ke^{i\lambda k}
\right)
\left(
 \sum_{k=0}^{\infty}(D^{\overline{\mu}}\me A\varphi^0_{\overline{\mu}})_ke^{i\lambda k}
\right)^{*}=
\\=
(\alpha^{2} +\gamma_1'(\lambda ))(\Phi^0_{\overline{\mu}}(e^{-i\lambda}))^{\top}B_{1}^{\top} \overline{\Phi^0_{\overline{\mu}}(e^{-i\lambda})},
\end{multline}
\begin{multline}  \label{eq_6_41}
\left(
 \sum_{k=0}^{\infty}(D^{\overline{\mu}}\me A\varphi^0_{\overline{\mu}})_ke^{i\lambda k}
\right)
\left(
 \sum_{k=0}^{\infty}(D^{\overline{\mu}}\me A\varphi^0_{\overline{\mu}})_ke^{i\lambda k}
\right)^{*}=
\\=
(\Phi^0_{\overline{\mu}}(e^{-i\lambda}))^{\top}(\vec{\alpha}\cdot \vec{\alpha}^{*}+\Gamma(\lambda)) \overline{\Phi^0_{\overline{\mu}}(e^{-i\lambda})},
\end{multline}

The following theorem  holds true.

\begin{thm}
Let the minimality condition (\ref{umova11_e_d}) hold true. The least favorable spectral densities $f_{0}(\lambda)$ in the classes $ \md D_{\varepsilon}^{k},k=1,2,3,4$ for the optimal linear estimation  of the functional  $A\vec{\xi}$ from observations of the sequence $\vec{\xi}(m)$ at points  $m=-1,-2,\ldots$  are determined by the equations
\eqref{eq_6_1}, \eqref{eq_6_2}, \eqref{eq_6_3}, \eqref{eq_6_4}
(or equations \eqref{eq_6_11}, \eqref{eq_6_21}, \eqref{eq_6_31}, \eqref{eq_6_41}
in the case where the spectral densities admit the canonical factorization (\ref{fact_f_e_d}), respectively),
the constrained optimization problem (\ref{minimax1_e_d}) and restrictions  on densities from the corresponding classes $ \md D_{\varepsilon}^{k},k=1,2,3,4$.  The minimax-robust spectral characteristic of the optimal estimate of the functional $A\vec{\xi}$ is determined by the formula (\ref{spectr A_e_d}).
\end{thm}

\subsection{Least favorable spectral density in classes $\md D_{1\delta}$}

Consider the forecasting problem for the  functional $A\vec{\xi}$
 which depends on unobserved values of a sequence $\vec\xi(m)$ with stationary GM increments based on observations of the sequence at points $m=-1,-2,\ldots$ under the condition that the sets of admissible spectral densities $ \md D_{1\delta}^{k},k=1,2,3,4$ are defined as follows:
\begin{equation*}
\md D_{1\delta}^{1}=\left\{f(\lambda )\biggl|\frac{1}{2\pi} \int_{-\pi}^{\pi}
\frac{|\chi_{\overline{\mu}}^{(d)}(e^{-i\lambda})|^2}{|\beta^{(d)}(i\lambda)|^2}
\left|{\rm{Tr}}(f(\lambda )-f_{1} (\lambda))\right|d\lambda \le \delta\right\};
\end{equation*}
\begin{equation*}
\md D_{1\delta}^{2}=\left\{f(\lambda )\biggl|\frac{1}{2\pi } \int_{-\pi}^{\pi}
\frac{|\chi_{\overline{\mu}}^{(d)}(e^{-i\lambda})|^2}{|\beta^{(d)}(i\lambda)|^2}
\left|f_{kk} (\lambda )-f_{kk}^{1} (\lambda)\right|d\lambda  \le \delta_{k}, k=\overline{1,T}\right\};
\end{equation*}
\begin{equation*}
\md D_{1\delta}^{3}=\left\{f(\lambda )\biggl|\frac{1}{2\pi } \int_{-\pi}^{\pi}
\frac{|\chi_{\overline{\mu}}^{(d)}(e^{-i\lambda})|^2}{|\beta^{(d)}(i\lambda)|^2}
\left|\left\langle B_{2} ,f(\lambda )-f_{1}(\lambda )\right\rangle \right|d\lambda  \le \delta\right\};
\end{equation*}
\begin{equation*}
\md D_{1\delta}^{4}=\left\{f(\lambda )\biggl|\frac{1}{2\pi} \int_{-\pi}^{\pi}
\frac{|\chi_{\overline{\mu}}^{(d)}(e^{-i\lambda})|^2}{|\beta^{(d)}(i\lambda)|^2}
\left|f_{ij} (\lambda )-f_{ij}^{1} (\lambda)\right|d\lambda  \le \delta_{i}^j, i,j=\overline{1,T}\right\}.
\end{equation*}

\noindent
Here  $f_{1} ( \lambda )$ is a fixed spectral density,  $\delta,\delta_{k},k=\overline{1,T}$, $\delta_{i}^{j}, i,j=\overline{1,T}$, are given numbers.

From the condition $0\in\partial\Delta_{\mathcal{D}}(f_0)$
we find the following equations which determine the least favourable spectral densities for these given sets of admissible spectral densities.

For the first set of admissible spectral densities $ \md D_{1\delta}^{1}$ we have equation
\begin{multline} \label{eq_7_1}
\left(
\sum_{k=0}^{\infty}((\me F^0_{\overline{\mu}})^{-1}D^{\overline{\mu}}\me a )_k e^{i\lambda k}
\right)
\left({
\sum_{k=0}^{\infty}((\me F^0_{\overline{\mu}})^{-1}D^{\overline{\mu}}\me a )_k e^{i\lambda k}}
\right)^{*}=
\\=
\beta^{2} \gamma_2( \lambda )\left(\frac{|\chi_{\overline{\mu}}^{(d)}(e^{-i\lambda})|^2}{|\beta^{(d)}(i\lambda)|^2} f_{0} (\lambda )\right)^{2},
\end{multline}
\begin{equation} \label{eq_7_11}
\frac{1}{2 \pi} \int_{-\pi}^{ \pi}
\frac{|\chi_{\overline{\mu}}^{(d)}(e^{-i\lambda})|^2}{|\beta^{(d)}(i\lambda)|^2}
\left|{\mathrm{Tr}}\, (f_0( \lambda )-f_{1}(\lambda )) \right|d\lambda =\delta,
\end{equation}
where $ \beta^{2}$ is Lagrange multiplier,
$ \left| \gamma_2( \lambda ) \right| \le 1$ and
\[\gamma_2( \lambda )={ \mathrm{sign}}\; ({\mathrm{Tr}}\, (f_{0} ( \lambda )-f_{1} ( \lambda ))): \; {\mathrm{Tr}}\, (f_{0} ( \lambda )-f_{1} ( \lambda )) \ne 0.\]

For the second set of admissible spectral densities $ \md D_{1\delta}^{2}$ we have equation
\begin{multline}   \label{eq_7_2}
\left(
\sum_{k=0}^{\infty}((\me F^0_{\overline{\mu}})^{-1}D^{\overline{\mu}}\me a )_k e^{i\lambda k}
\right)
\left({
\sum_{k=0}^{\infty}((\me F^0_{\overline{\mu}})^{-1}D^{\overline{\mu}}\me a )_k e^{i\lambda k}}
\right)^{*}=
\\
=
\left(\frac{|\chi_{\overline{\mu}}^{(d)}(e^{-i\lambda})|^2}{|\beta^{(d)}(i\lambda)|^2} f_{0} (\lambda )\right) \left \{ \beta_{k}^{2} \gamma^2_{k} ( \lambda ) \delta_{kl} \right \}_{k,l=1}^{T} \left(\frac{|\chi_{\overline{\mu}}^{(d)}(e^{-i\lambda})|^2}{|\beta^{(d)}(i\lambda)|^2} f_{0} (\lambda )\right),
\end{multline}
\begin{equation} \label{eq_7_21}
\frac{1}{2 \pi} \int_{- \pi}^{ \pi} \frac{|\chi_{\overline{\mu}}^{(d)}(e^{-i\lambda})|^2}{|\beta^{(d)}(i\lambda)|^2} \left|f^0_{kk} ( \lambda)-f_{kk}^{1} ( \lambda ) \right| d\lambda =\delta_{k},
\end{equation}
where  $\beta_{k}^{2}$ are Lagrange multipliers,  $\left| \gamma^2_{k} ( \lambda ) \right| \le 1$ and
\[\gamma_{k}^2( \lambda )={ \mathrm{sign}}\;(f_{kk}^{0}( \lambda)-f_{kk}^{1} ( \lambda )): \; f_{kk}^{0} ( \lambda )-f_{kk}^{1}(\lambda ) \ne 0, \; k= \overline{1,T}.\]

For the third set of admissible spectral densities $ \md D_{1\delta}^{3}$ we have equation
\begin{multline}   \label{eq_7_3}
\left(
\sum_{k=0}^{\infty}((\me F^0_{\overline{\mu}})^{-1}D^{\overline{\mu}}\me a )_k e^{i\lambda k}
\right)
\left({
\sum_{k=0}^{\infty}((\me F^0_{\overline{\mu}})^{-1}D^{\overline{\mu}}\me a )_k e^{i\lambda k}}
\right)^{*}=
\\
=
\beta^{2} \gamma_2'( \lambda )\left(\frac{|\chi_{\overline{\mu}}^{(d)}(e^{-i\lambda})|^2}{|\beta^{(d)}(i\lambda)|^2} f_{0} (\lambda )\right)B_{2}^{ \top}\,\left(\frac{|\chi_{\overline{\mu}}^{(d)}(e^{-i\lambda})|^2}{|\beta^{(d)}(i\lambda)|^2} f_{0} (\lambda )\right),
\end{multline}
\begin{equation} \label{eq_7_31}
\frac{1}{2 \pi} \int_{- \pi}^{ \pi} \frac{|\chi_{\overline{\mu}}^{(d)}(e^{-i\lambda})|^2}{|\beta^{(d)}(i\lambda)|^2} \left| \left \langle B_{2}, f_0( \lambda )-f_{1} ( \lambda ) \right \rangle \right|d\lambda
= \delta,
\end{equation}
where $\beta^{2}$ is a Lagrange multiplier,  $\left| \gamma_2' ( \lambda ) \right| \le 1$ and
\[\gamma_2' ( \lambda )={ \mathrm{sign}}\; \left \langle B_{2} ,f_{0} ( \lambda )-f_{1} ( \lambda ) \right \rangle : \; \left \langle B_{2} ,f_{0} ( \lambda )-f_{1} ( \lambda ) \right \rangle \ne 0.\]

For the fourth set of admissible spectral densities $ \md D_{1\delta}^{4}$ we have equation
\begin{multline}  \label{eq_7_4}
\left(
\sum_{k=0}^{\infty}((\me F^0_{\overline{\mu}})^{-1}D^{\overline{\mu}}\me a )_k e^{i\lambda k}
\right)
\left({
\sum_{k=0}^{\infty}((\me F^0_{\overline{\mu}})^{-1}D^{\overline{\mu}}\me a )_k e^{i\lambda k}}
\right)^{*}=
\\
=
\left(\frac{|\chi_{\overline{\mu}}^{(d)}(e^{-i\lambda})|^2}{|\beta^{(d)}(i\lambda)|^2} f_{0} (\lambda )\right)\left \{ \beta_{ij}( \lambda ) \gamma_{ij} ( \lambda ) \right \}_{i,j=1}^{T} \left(\frac{|\chi_{\overline{\mu}}^{(d)}(e^{-i\lambda})|^2}{|\beta^{(d)}(i\lambda)|^2} f_{0} (\lambda )\right),
\end{multline}
\begin{equation} \label{eq_7_41}
\frac{1}{2 \pi} \int_{- \pi}^{ \pi} \frac{|\chi_{\overline{\mu}}^{(d)}(e^{-i\lambda})|^2}{|\beta^{(d)}(i\lambda)|^2} \left|f^0_{ij}(\lambda)-f_{ij}^{1}( \lambda ) \right|d\lambda = \delta_{i}^{j},
\end{equation}
where  $ \beta_{ij}$ are Lagrange multipliers,  $\left| \gamma_{ij} ( \lambda ) \right| \le 1$ and
\[\gamma_{ij} ( \lambda )= \frac{f_{ij}^{0} ( \lambda )-f_{ij}^{1} (\lambda )}{ \left|f_{ij}^{0} ( \lambda )-f_{ij}^{1}(\lambda) \right|} : \; f_{ij}^{0} ( \lambda )-f_{ij}^{1} ( \lambda ) \ne 0, \; i,j= \overline{1,T}.
\]

In the case where the spectral density admits the canonical factorization (\ref{fact_f_e_d}) we have the following equations, correspondingly
\begin{multline} \label{eq_7_111}
\left(
 \sum_{k=0}^{\infty}(D^{\overline{\mu}}\me A\varphi^0_{\overline{\mu}})_ke^{i\lambda k}
\right)
\left(
 \sum_{k=0}^{\infty}(D^{\overline{\mu}}\me A\varphi^0_{\overline{\mu}})_ke^{i\lambda k}
\right)^{*}=
\\=
\beta^{2} \gamma_2( \lambda )(\Phi^0_{\overline{\mu}}(e^{-i\lambda}))^{\top}\overline{\Phi^0_{\overline{\mu}}(e^{-i\lambda})},
\end{multline}
\begin{multline}   \label{eq_7_211}
\left(
 \sum_{k=0}^{\infty}(D^{\overline{\mu}}\me A\varphi^0_{\overline{\mu}})_ke^{i\lambda k}
\right)
\left(
 \sum_{k=0}^{\infty}(D^{\overline{\mu}}\me A\varphi^0_{\overline{\mu}})_ke^{i\lambda k}
\right)^{*}=
\\=
(\Phi^0_{\overline{\mu}}(e^{-i\lambda}))^{\top} \left \{ \beta_{k}^{2} \gamma^2_{k} ( \lambda ) \delta_{kl} \right \}_{k,l=1}^{T} \overline{\Phi^0_{\overline{\mu}}(e^{-i\lambda})},
\end{multline}
\begin{multline}   \label{eq_7_311}
\left(
 \sum_{k=0}^{\infty}(D^{\overline{\mu}}\me A\varphi^0_{\overline{\mu}})_ke^{i\lambda k}
\right)
\left(
 \sum_{k=0}^{\infty}(D^{\overline{\mu}}\me A\varphi^0_{\overline{\mu}})_ke^{i\lambda k}
\right)^{*}=
\\=
\beta^{2} \gamma_2'( \lambda )(\Phi^0_{\overline{\mu}}(e^{-i\lambda}))^{\top}  B_{2}^{ \top} \overline{\Phi^0_{\overline{\mu}}(e^{-i\lambda})},
\end{multline}
\begin{multline}  \label{eq_7_411}
\left(
 \sum_{k=0}^{\infty}(D^{\overline{\mu}}\me A\varphi^0_{\overline{\mu}})_ke^{i\lambda k}
\right)
\left(
 \sum_{k=0}^{\infty}(D^{\overline{\mu}}\me A\varphi^0_{\overline{\mu}})_ke^{i\lambda k}
\right)^{*}=
\\=
(\Phi^0_{\overline{\mu}}(e^{-i\lambda}))^{\top} \left \{ \beta_{ij}( \lambda ) \gamma_{ij} ( \lambda ) \right \}_{i,j=1}^{T} \overline{\Phi^0_{\overline{\mu}}(e^{-i\lambda})},
\end{multline}

The following theorem  holds true.

\begin{thm}
Let the minimality condition (\ref{umova11_e_d}) hold true. The least favorable spectral densities $f_{0}(\lambda)$ in the classes $ \md D_{1\delta}^{k},k=1,2,3,4$  for the optimal linear estimation  of the functional  $A\vec{\xi}$ from observations of the sequence $\vec{\xi}(m)$ at points  $m=-1,-2,\ldots$  are determined by  equations
\eqref{eq_7_1}, \eqref{eq_7_11}; \eqref{eq_7_2}, \eqref{eq_7_21}; \eqref{eq_7_3}, \eqref{eq_7_31}; \eqref{eq_7_4}, \eqref{eq_7_41}
(or equations \eqref{eq_7_111}, \eqref{eq_7_211}, \eqref{eq_7_311}, \eqref{eq_7_411}
in the case where the spectral densities admit the canonical factorization (\ref{fact_f_e_d}), respectively),
the constrained optimization problem (\ref{minimax1_e_d}) and restrictions  on densities from the corresponding classes $ \md D_{1\delta}^{k},k=1,2,3,4$.  The minimax-robust spectral characteristic of the optimal estimate of the functional $A\vec{\xi}$ is determined by the formula (\ref{spectr A_e_d}).
\end{thm}

\section{Conclusions}

In this article, we present results of investigation of stochastic sequences with periodically stationary long memory multiple seasonal increments.
We give definition of generalized multiple increment sequence and introduce stochastic sequences $\zeta(m)$ with periodically stationary (periodically correlated, cyclostationary) generalized multiple  increments.
These non-stationary stochastic sequences combine  periodic structure of covariation functions of sequences as well as multiple seasonal factors,
including the integrating one. A short review of the spectral theory of vector-valued generalized multiple increment sequences is presented.
We describe methods of solution of the classical forecasting problem for linear functionals which are constructed from unobserved values of a sequence with periodically stationary generalized multiple  increments in the case where the spectral structure of the sequence is exactly known.
Estimates are obtained by representing the sequence under investigation as a vector-valued sequence with stationary generalized multiple increments and applying the Hilbert space projection technique.
An approach to forecasting in the presence of non-stationary  fractional integration is discussed.
Examples of solution of the forecasting problem for particular models of time series are proposed.
The minimax-robust approach to forecasting problem is applied in the case of spectral uncertainty where densities of sequences are not exactly known while, instead, sets of admissible spectral densities are specified.
We propose a representation of the mean square error in the form of a linear functional in $L_1$ with respect to spectral densities, which allows
us to solve the corresponding constrained optimization problem and describe the minimax (robust) estimates of the functionals.
Relations are described which determine the least favourable spectral densities and the minimax spectral characteristics of the optimal estimates of linear functionals
for a collection of specific classes
of admissible spectral densities.

\newpage

\appendix

\section{Proofs}

 \emph{Proof of Lemma \ref{lem_multiplicative_Pryrist}}

We have
\begin{eqnarray*}
    \prod_{i=1}^r(1-B_{\mu_i}^{s_i})^{d_i} & = & \prod_{i=1}^r\ld(\sum_{j_i=0}^{d_i}(-1)^{j_i}{d_i \choose j_i}B^{j_i\mu_i s_i}\rd)
\\
 & = & \prod_{i=1}^r\ld(\sum_{j_i=0}^{d_i\mu_i s_i}(-1)^{[j_i/\mu_i s_i]}{d_i \choose [j_i/\mu_i s_i]}\mr I\{j_i\mod \mu_i s_i=0\}B^{j_i}\rd)
 \\
     & = & \sum_{j_1=0}^{n_1}\ldots\sum_{j_r=0}^{n_r}\ld((-1)^{\sum_{i=1}^rM_i^{j_i}}
\prod_{i=1}^rI_i^{j_i}
\prod_{i=1}^r{d_i \choose M_i^{j_i}}\rd)B^{\sum_{i=1}^rj_i}.
   \end{eqnarray*}
By replacing consequently $j_1\to k_1$, $k_1+j_2\to k_2$, $k_2+j_3\to k_3$, \ldots, $k_{r-1}+j_r\to k_r:=k$, we obtain
\[
    \prod_{i=1}^r(1-B_{\mu_i}^{s_i})^{d_i}=\sum_{k_2=0}^{n_1+n_2}\sum_{k_1=0\vee k_2-n_2}^{n_1\wedge k_2}\sum_{j_3=0}^{n_3}\ldots\sum_{j_r=0}^{n_r}(-1)^{M_2^{k_2-k_1}+M_1^{k_1}}I_2^{k_2-k_1}I_1^{k_1}
\]
\[
\times{d_2 \choose M_2^{k_2-k_1}}{d_1 \choose M_1^{k_1}}\ld((-1)^{\sum_{i=3}^rM_i^{j_i}}
\prod_{i=3}^rI_i^{j_i}
\prod_{i=3}^r{d_i \choose M_i^{j_i}}\rd)B^{k_2+\sum_{i=3}^rj_i}
    \]
\begin{eqnarray*}
     & = & \sum_{k_3=0}^{n_1+n_2+n_3}\sum_{k_2=0\vee k_3-n_3}^{n_1+n_2\wedge k_3}\sum_{k_1=0\vee k_2-n_2}^{n_1\wedge k_2}\sum_{j_4=0}^{n-4}\ldots\sum_{j_r=0}^{n_r}
(-1)^{\sum_{i=1}^3M_i^{k_i-k_{i-1}}}
\prod_{i=1}^3I_i^{k_i-k_{i-1}}
\\
 & \quad & \times\prod_{i=1}^3{d_i \choose M_i^{k_i-k_{i-1}}}\ld((-1)^{\sum_{i=4}^rM_i^{j_i}}
\prod_{i=3}^rI_i^{j_i}
\prod_{i=3}^r{d_i \choose M_i^{j_i}}\rd)B^{k_3+\sum_{i=4}^rj_i}
\\
 & = & \sum_{k=0}^{n(\gamma)}e_{\gamma}(k)B^k,
\end{eqnarray*}
where $e_{\gamma}(k)$ are coefficients from the lemma statement. $\square$

\

 \emph{Proof of Theorem \ref{thm1}}

We follow the idea proposed by Yaglom \cite{Yaglom:1955} for continuous time stationary increments.  Consider a GM increment sequence with one seasonal factor $\chi_{\mu,s}^{(d)}(\eta(m))=\eta^{(d)}_s(m,\mu)=(1-B_{\mu}^s)^d\eta(m)$.
Formula  $(\ref{tot2})$ implies
\[
    c^{(d)}_s(\mu)=\mt
E\eta^{(d)}_s(m,\mu)=(A_0+A_1+\ldots+A_{(\mu-1)d})c^{(d)}_s(1)=\mu^dc^{(d)}_s(1)=c\mu^d,\]
where $c=c^{(d)}_s(1)$ does not depend on $\mu$.

Since $D^{(n)}_s(m;\mu,\mu)$ is a positive-definite function with respect to variable  $m$, one can define a function $F_{\mu,s}(\lambda)$ depending on the parameter $\mu$, which is a real bounded non-decreasing left-continuous with respect to
$\lambda\in[-\pi,\pi)$ function, such that
\begin{equation}D^{(d)}_s(m;\mu,\mu)=\int_{-\pi}^{\pi}e^{i\lambda
m}dF_{\mu,s}(\lambda).\label{zobr1}\end{equation}

Again, formula $(\ref{tot2})$ implies
\begin{eqnarray*}
 &  & \mt E\eta^{(d)}_s(m_1+m,\mu)\overline{\eta^{(d)}_s(m,\mu)}
\\ & = &
\sum_{p=0}^{(\mu-1)d}\sum_{q=0}^{(\mu-1)d}A_pA_q\mt
E\eta^{(d)}_s(m_1+m-ps,1)\overline{\eta^{(d)}_s(m_2-qs,1)}
\\
 & = & \int_{-\pi}^{\pi}\sum_{p=0}^{(\mu-1)d}\sum_{q=0}^{(\mu-1)d}A_pA_qe^{i\lambda(m-(p-q)s)}dF_{1,s}(\lambda)
\\
 & = & \int_{-\pi}^{\pi}e^{i\lambda m}\sum_{p=0}^{(k-1)d}A_pe^{-ips\lambda}\sum_{q=0}^{(\mu-1)d}A_qe^{iqs\lambda}dF_{1,s}(\lambda)
\\
 & = & \int_{\pi}^{\pi}e^{im\lambda}\frac{(1-e^{-i\mu s\lambda})^d}{(1-e^{-is\lambda})^d}
\frac{(1-e^{i\mu s\lambda})^d}{(1-e^{is\lambda})^d}dF_{1,s}(\lambda).
\end{eqnarray*}
Thus,
\begin{equation}\int_{\pi}^{\pi}e^{im\lambda}dF_{\mu, s}(\lambda)=
\int_{-\pi}^{\pi}e^{im\lambda}\frac{(1-\cos
\mu s\lambda)^d}{(1-\cos s \lambda)^d}F_{1,s}(\lambda).\label{teoremSpectr1}\end{equation}
The latter equality implies
\[
F_{\mu,s}(\lambda)=\int_0^{\lambda}\frac{(1-\cos \mu su)^d}{(1-\cos su)^d}F_{1,s}(u),\]
and
\begin{equation}
\int_0^{\lambda}\frac{|\beta^{(d)}(iu)|^2}{(1-\cos \mu su)^d}dF_{\mu,s}(u)
=\int_0^{\lambda}\frac{|\beta^{(d)}(iu)|^2}{(1-\cos su)^d}dF_{1,s}(u),\label{FormulaSpectrFunk2}\end{equation}
where the function $\beta^{(d)}(iu)$ is to be chosen in the way that the integrals are defined for $\lambda\in[-\pi,\pi)$ and converge at the neighborhoods of the points $\cos su=1$, $u\in[-\pi,\pi]$, namely, the points $u=2\pi k/s$ for $|k|\leq s/2$, $k\in \mr Z$. Thus, we choose the function
$\beta^{(d)}(iu)=  \prod_{k =-[s /2]}^{[s /2]}(iu-2\pi i k /s )^d$.

The right side of equality (\ref{FormulaSpectrFunk2}) doesn't depend on $k$, thus, put
\begin{equation}
F(\lambda)=\frac{1}{2^d}\int_0^{\lambda}\frac{|\beta^{(d)}(iu)|^2}{(1-\cos \mu su)^d}dF_{\mu,s}(u)\label{FormulaSpectrFunkZag}\end{equation}
The function $F(\lambda)$ is real-valued non-decreasing bounded function defined on  $[-\pi,\pi)$, such that $F(0)=0$.
Consider the function $F(\lambda)$ as left-continuous. This function is the one stated in the theorem.

Relation   $(\ref{strFnaR})$ for $r=1$ is obtained by considering the following equality for positive  $\mu_1$, $\mu_2$:
\begin{eqnarray*}
D^{(n)}_s(m,\mu_1,\mu_2) & = & \ip
e^{im\lambda}\frac{(1-e^{-i\mu_1s\lambda})^d}{(1-e^{-is\lambda})^d}
\frac{(1-e^{i\mu_2s\lambda})^d}{(1-e^{is\lambda})^d}dF_{1,s}(\lambda)
\\
 & = & \ip
e^{im\lambda}(1-e^{-i\mu_1s\lambda})^d(1-e^{i\mu_2s\lambda})^d\frac{1}{|\beta^{(d)}(i\lambda)|^2}dF(\lambda).
\end{eqnarray*}

For negative $\mu_1$, $\mu_2$,  equality
$(\ref{tot1})$ is applied.

By generalizing for $r>1$ the given reasonings, we obtain the relations $(\ref{serFnaR})$ and $(\ref{strFnaR})$. $\square$

\

\emph{Proof of Lemma \ref{lema predst A}}

Using Definition \ref{def_multiplicative_Pryrist} we obtain the formal equality
\[
\xi_p(k)=\frac{1}{(1-B_{\mu})^n}\chi_{\overline{\mu},\overline{s}}^{(n)}(\xi_p(k))=
\sum_{j=-\infty}^kd_{\mu}(k-j)\chi_{\overline{\mu},\overline{s}}^{(n)}(\xi_p(j)),
\]
which imply
\be
\label{sss1_diskr}
 \sum_{k=0}^{\infty}a_p(k)\xi_p(k)
  =
 -\sum_{i=-n(\gamma)}^{-1}v_p(i)\xi_p(i)
+\sum_{i=0}^{\infty}
 \ld(\sum_{k=i}^{\infty}a_p(k)d_{\mu}(k-i)\rd)\chi_{\overline{\mu},\overline{s}}^{(n)}(\xi_p(i)),
 \ee
 and
 \begin{eqnarray}
 \nonumber
 \sum_{i=0}^{\infty}b_p(i)\chi_{\overline{\mu},\overline{s}}^{(n)}(\xi_p(i))
 & = &  \sum_{k=-n(\gamma)}^{-1}\xi_p(k)\sum_{l=0}^{k+n(\gamma)}e_{\nu}(l-k)b_p(k)
\\
 & \quad & +\sum_{k=0}^{\infty}\xi_p(k)\sum_{l=k}^{k+n(\gamma)}e_{\nu}(l-k)b_p(k),
 \label{sss2_diskr}
\end{eqnarray}
Relations (\ref{sss1_diskr}) and (\ref{sss2_diskr}) imply a representation  $A\vec{\xi}=B\vec{\xi}-V\vec{\xi}$ and relations which prove the lemma:
\begin{eqnarray*}
v_p(k) & = & \sum_{l=0}^{k+n(\gamma)}  e_{\nu}(l-k)b_p(l), \quad k=-1,-2,\dots,-n(\gamma), \quad p=1, 2,\dots,T,
\\
b_p(k) & = & \sum_{m=k}^{\infty}d_{\overline{\mu}}(m-k)a_p(m),  \quad k=0,1,2,\dots, \quad p=1, 2,\dots,T.\quad \square
\end{eqnarray*}

\end{document}